\documentclass[10pt,oneside,a4paper]{amsart}

\usepackage{amsmath,amsfonts,amssymb,amsthm}
\usepackage{array}
\usepackage[all,textures]{xy}
\usepackage{color}
\usepackage[orange]{xcolor}
\definecolor{blue(munsell)}{rgb}{0.0, 0.5, 0.69}
\usepackage{hyperref}
\hypersetup{hypertexnames = false, bookmarksdepth = 2, bookmarksopen = true, colorlinks, linkcolor = black, citecolor = blue(munsell), urlcolor = blue(munsell), pdfstartview={XYZ null null 1}}
\usepackage{tikz-cd}
\usepackage{etoolbox}
\usepackage{cleveref}
\usepackage{mathtools}
\usepackage{placeins}
\usepackage{caption}
\usepackage{enumitem}

\usepackage[T1]{fontenc}
\usepackage{newtxtext,newtxmath}

\usepackage{pdfsync}

\usepackage[colorinlistoftodos]{todonotes}
\usepackage{xparse}

\pdfmapfile{=ntx-Bold-tlf-t1.map}

\DeclareDocumentCommand\issue{g}{\todo[size=\footnotesize,color = green!40]{Wendy\IfNoValueF{#1}{: #1}}}
\DeclareDocumentCommand\tobedone{g}{\todo[size=\footnotesize,color = yellow!50]{To be done\IfNoValueF{#1}{: #1}}}
\DeclareDocumentCommand\notationissue{g}{\todo[size=\footnotesize,color = red!30]{Notation?\IfNoValueF{#1}{: #1}}}
\DeclareDocumentCommand\doubt{g}{\todo[size=\footnotesize,color = blue!10]{Doubt\IfNoValueF{#1}{: #1}}}
\DeclareDocumentCommand\observation{g}{\todo[size=\footnotesize,color = orange!10]{Observation\IfNoValueF{#1}{: #1}}}

\theoremstyle{plain}
\newtheorem{theorem}{Theorem}[section]
\newtheorem{proposition}[theorem]{Proposition}
\newtheorem{lemma}[theorem]{Lemma}
\newtheorem{corollary}[theorem]{Corollary}

\theoremstyle{definition}
\newtheorem{definition}[theorem]{Definition}

\theoremstyle{remark}
\newtheorem{remarks}[theorem]{Remarks}
\newtheorem{remark}[theorem]{Remark}
\newtheorem{example}[theorem]{Example}

\newcommand{\Img}{\mathrm{Im}}

\newcommand{\Obj}{\mathrm{Obj}}

\newcommand{\Tri}{\ensuremath{\mathsf{Tr}}}

\newcommand{\Ho}{\ensuremath{\mathsf{Ho}}}

\newcommand{\qrep}{\ensuremath{\mathsf{qrep}}}

\newcommand{\op}{\ensuremath{\mathsf{op}}}

\newcommand{\id}{\text{Id}}

\newcommand{\Kern}{\mathrm{Ker}}

\newcommand{\D}{\ensuremath{\mathsf{D}}}
\newcommand{\C}{\ensuremath{\mathsf{C}}}

\renewcommand{\lim}{\mathrm{lim}}

\newcommand{\hocolim}{\mathrm{hocolim}}

\newcommand{\RHom}{\mathrm{RHom}}

\newcommand{\AAA}{\mathfrak{a}}
\newcommand{\BBB}{\mathfrak{b}}

\newcommand{\CCC}{\mathfrak{c}}
\newcommand{\DDD}{\mathfrak{d}}
\newcommand{\EEE}{\mathfrak{e}}

\newcommand{\DD}{\boldsymbol{\mathsf{D}}}

\newcommand{\Ab}{\ensuremath{\mathsf{Ab}} }

\newcommand{\dgcat}{\ensuremath{\mathsf{dgcat}}}

\newcommand{\hodgcat}{\ensuremath{\mathsf{Hqe}}}

\newcommand{\dgmod}{\ensuremath{\mathsf{dgMod}}}
\newcommand{\Ac}{\ensuremath{\mathsf{Ac}}}

\newcommand{\wg}{\ensuremath{\mathsf{wg}}}
\newcommand{\ccm}{\ensuremath{\mathsf{c}}}

\newcommand{\iso}{\ensuremath{\mathsf{iso}}}
\newcommand{\Iso}{\ensuremath{\mathrm{Iso}}}

\newcommand{\dertp}{\ensuremath{\otimes^{\mathrm{L}}}}
\newcommand{\hztp}{\ensuremath{\otimes_{H^0}}}

\newcommand{\Ind}{\ensuremath{\mathsf{Ind}}}

\newcommand{\Fun}{\ensuremath{\mathsf{Fun}}}

\newcommand{\lra}{\longrightarrow}

\newcommand{\ra}{\rightarrow}

\newcommand{\aaa}{\ensuremath{\mathcal{A}}}
\newcommand{\bbb}{\ensuremath{\mathcal{B}}}
\newcommand{\ccc}{\ensuremath{\mathcal{C}}}

\newcommand{\eee}{\ensuremath{\mathcal{E}}}

\newcommand{\GGG}{\ensuremath{\mathcal{G}}}
\newcommand{\hhh}{\ensuremath{\mathcal{H}}}

\newcommand{\LLL}{\ensuremath{\mathcal{L}}}

\newcommand{\nnn}{\ensuremath{\mathcal{N}}}

\newcommand{\ppp}{\ensuremath{\mathcal{P}}}

\newcommand{\ttt}{\ensuremath{\mathcal{T}}}

\newcommand{\www}{\ensuremath{\mathcal{W}}}
\newcommand{\xxx}{\ensuremath{\mathcal{X}}}

\newcommand{\zzz}{\ensuremath{\mathcal{Z}}}

\newcommand{\fraku}{\ensuremath{\mathfrak{U}}}
\newcommand{\frakv}{\ensuremath{\mathfrak{V}}}
\newcommand{\frakw}{\ensuremath{\mathfrak{W}}}

\title{On the tensor product of well generated dg categories}	
	
	\author{Wendy Lowen} 
\address[Wendy Lowen]{Universiteit Antwerpen, Departement Wiskunde, Middelheimcampus,
Middelheimlaan 1,
2020 Antwerp, Belgium}
\address{Laboratory of Algebraic Geometry, National Research University, Higher School of Economics, Moscow, Russia}
\email{wendy.lowen@uantwerpen.be}

\author{Julia Ramos Gonz\'alez}
\address[Julia Ramos Gonz\'alez]{Universiteit Antwerpen, Departement Wiskunde, Middelheimcampus,
	Middelheimlaan 1,
	2020 Antwerp, Belgium}
\email{julia.ramosgonzalez@uantwerpen.be}

\thanks{This project has received funding from the European Research Council (ERC) under the European Union's Horizon 2020 research and innovation programme (grant agreement No. 817762). The authors acknowledge the support of the Research Foundation - Flanders (FWO) under Grant No G.0D86.16N and of the Russian Academic Excellence Project `5-100'.\\
The second named author is a postdoctoral fellow of the Research Foundation -  Flanders (FWO)}

\usepackage{pdfsync}
\begin{document}
\begin{abstract}
We endow the homotopy category of well generated (pretriangulated) dg ca\-tegories with a tensor product satisfying a universal property. The resulting monoidal structure is symmetric and closed with respect to the cocontinuous RHom of dg categories (in the sense of To\"en \cite{homotopy-theory-dg-categories-derived-morita-theory}). We give a construction of the tensor product in terms of localisations of dg derived categories, making use of the enhanced derived Gabriel-Popescu theorem  \cite{the-popescu-gabriel-theorem-for-triangulated-categories}.
Given a regular cardinal $\alpha$, we define and construct a tensor product of homotopically $\alpha$-cocomplete dg categories and prove that the well generated tensor product of $\alpha$-continuous derived dg categories (in the sense of \cite{the-popescu-gabriel-theorem-for-triangulated-categories}) is the $\alpha$-continuous dg derived category of the homotopically $\alpha$-cocomplete tensor product. In particular, this shows that the tensor product of well generated dg categories preserves $\alpha$-com\-pact\-ness. 
\end{abstract}

\maketitle

\section{Introduction}
The main aim of this paper is the development of a suitable tensor product for \emph{well generated dg categories}, that is, pretriangulated dg categories $\aaa$ for which $H^0(\aaa)$ is well generated in the sense of Neeman \cite{triangulated-categories}. Well generated triangulated categories were introduced in loc. cit. as a natural class of triangulated categories sharing important properties like Brown representability with the subclass of compactly generated triangulated categories, while at the same time having a good localisation theory (see \cite{triangulated-categories} and \cite{localization-theory-triangulated-categories}). The derived category of a Grothendieck abelian category being well generated \cite{derived-category-sheaves-manifold}, there is a rich supply of examples of algebro-geometric origin and in the spririt of noncommutative geometry, our tensor product can be thought of as a kind of (derived) product of noncommutative spaces.

Our starting point is the homotopy category of dg categories $\hodgcat$ developed by Tabuada \cite{structure-categorie-modeles-quillen-categorie-dg-categories} and To\"en \cite{homotopy-theory-dg-categories-derived-morita-theory}. As shown in \cite{homotopy-theory-dg-categories-derived-morita-theory}, $\hodgcat$ has a monoidal structure given by the derived tensor product of dg categories $\dertp$ and this monoidal structure is closed with the internal hom (denoted by $\RHom$) given by the dg category of (cofibrant) right quasi-representable bimodules (also called quasi-functors). 

When we restrict our attention to dg categories $\aaa, \bbb$ that are (homotopically) cocomplete, it is natural to restrict to quasi-functors $F \in \RHom(\aaa, \bbb)$ whose associated underlying exact functor $H^0(F): H^0(\aaa) \lra H^0(\bbb)$ preserves coproducts. These will be called \emph{cocontinuous quasi-functors} and they form a full dg subcategory $\RHom_{\ccm}(\aaa, \bbb) \subseteq \RHom(\aaa, \bbb)$. We show (\Cref{cocompletecontinuousRHom} and \Cref{closed}):

\begin{theorem}\label{Thm0intro}
Consider pretriangulated dg categories $\aaa$ and $\bbb$.
\begin{enumerate}
\item If $\aaa$ and $\bbb$ are homotopically cocomplete, the same holds for $\RHom_{\ccm}(\aaa, \bbb)$.
\item If $\aaa$ and $\bbb$ are well generated, the same holds for $\RHom_{\ccm}(\aaa, \bbb)$.
\end{enumerate}
\end{theorem}

We define the \emph{well generated tensor product} of two well generated dg categories $\aaa$ and $\bbb$, if it exists, as the unique well generated dg category $\aaa \boxtimes \bbb$ satisfying the following universal property in $\hodgcat$ with respect to all well generated dg categories $\ccc$:
\begin{equation}\label{eqbox}
\RHom_{\ccm}(\aaa \boxtimes \bbb, \ccc) \cong \RHom_{\ccm}(\aaa, \RHom_{\ccm}(\bbb, \ccc)).
\end{equation}

Our main result is the existence of the well generated tensor product (see Theorem \ref{Thm1intro} below). In combination with Theorem \ref{Thm0intro} (2) we immediately obtain:

\begin{corollary}
The homotopy category $\hodgcat_{\wg}$ of well generated dg categories with cocontinuous quasi-functors is symmetric monoidal closed.
\end{corollary} 

Our approach to the existence of the tensor product makes use of the localisation theory of well generated dg categories. More precisely, we use the (enhanced) derived Gabriel-Popescu theorem from \cite{the-popescu-gabriel-theorem-for-triangulated-categories} which identifies the well generated dg categories in $\hodgcat$ as the dg quotients of dg derived categories $\DD(\AAA)$ by an (enhanced) localising subcategory $\www \subseteq \DD(\AAA)$ generated by a set, for small dg categories $\AAA$. We show:

\begin{theorem}\label{Thm1intro}
	Let $\aaa$, $\bbb$ be two well generated dg categories such that $\aaa \cong \DD(\AAA)/\www_{\AAA}$ and $\bbb \cong \DD(\BBB)/\www_{\BBB}$ for small dg categories $\AAA$, $\BBB$ with $\www_{\AAA}\subseteq \DD(\AAA)$ and $\www_{\BBB}\subseteq \DD(\BBB)$ (enhanced) localising subcategories generated by a set of objects. 
	There exists an (enhanced) localising subcategory $\www_{\AAA} \boxtimes \www_{\BBB} \subseteq \DD(\AAA \dertp \BBB)$ such that the well generated tensor product of $\aaa$ and $\bbb$ exists and is given  by the dg quotient
	\begin{equation}
		\aaa \boxtimes \bbb = \DD(\AAA \dertp \BBB)/\www_{\AAA} \boxtimes \www_{\BBB}.
	\end{equation}
	In particular, $\aaa \boxtimes \bbb$ is independent of the chosen realisations of $\aaa$ and $\bbb$.
\end{theorem}

In the paper, we give a description of $\www_{\AAA} \boxtimes \www_{\BBB}$ in terms of generators (Theorem \ref{thmexist}) as well as an intrinsic description (Theorem \ref{tplocalisingsubcats}).
We also give a description of the well generated tensor product in terms of Bousfield localisations (Theorem \ref{deftpstrictloc}) which is specifically applied to $\alpha$-continuous dg derived categories in the sense of \cite{the-popescu-gabriel-theorem-for-triangulated-categories} (we call them $\alpha$-cocontinuous in line with the rest of our terminology).
More precisely, we show (Theorem \ref{thmalphatp}, Proposition \ref{tpsmalltobig}, Corollary \ref{alphapreserved}):

\begin{theorem} Let $\alpha$ be a regular cardinal.
Let $\AAA$, $\BBB$ be two homotopically $\alpha$-cocomplete small dg categories. Then, we have that
	\begin{equation}
	\DD_{\alpha}(\AAA) \boxtimes \DD_{\alpha}(\BBB) \cong \DD_{\alpha}(\AAA \dertp_{\alpha} \BBB)
	\end{equation}
	in $\hodgcat_{\wg}$, where $\AAA \dertp_{\alpha} \BBB$ is the homotopically $\alpha$-cocomplete tensor product of $\AAA$ and $\BBB$.
	
	In particular, the well generated tensor product preserves $\alpha$-compactness.
\end{theorem}

\begin{remarks}\hspace{0.5cm}
\begin{enumerate}
\item In \cite{tensor-product-linear-sites-grothendieck-categories}, a tensor product of Grothendieck abelian categories was defined. The precise relationship between this tensor product and the tensor product of well generated dg categories (with t-structures) is currently under investigation in a joint project with Francesco Genovese and Michel Van den Bergh.
\item In contrast to the tensor product of well generated dg categories, the tensor product of Grothendieck categories from \cite{tensor-product-linear-sites-grothendieck-categories} is not closed (as follows for instance from \cite[Rem 6.5]{covers-envelopes-cotorsion-theories-locally-presentable-abelian-categories-contramodule-categories}). An in depth study of the nature of morphism categories between abelian categories is the topic of an ongoing joint project with Michel Van den Bergh.
\item There is well known correspondence between pretriangulated dg categories on the one hand and stable linear infinity categories on the other hand, see for instance \cite{differential-graded-categories-k-linear-stable-infinity-categories}. Since a pretriangulated dg category is well generated precisely when it is locally presented \cite[\S 2.1]{derived-azumaya-algebras-generators-twisted-derived-categories}, we expect our tensor product to correspond to a natural tensor product of presentable stable linear infinity categories. Such a tensor product can be obtained as a linear analogue of the tensor product of presentable stable infinity categories from \cite{DAGII-noncommutative-algebra, higher-algebra}. The details of such a monoidal correspondence remain to be elucidated. 
\end{enumerate}
\end{remarks}
The present work extends part of the work carried out by the second named author in her PhD thesis under the supervision of Wendy Lowen and Boris Shoikhet.

\vspace{0,3cm}
\noindent \emph{Acknowledgements.} 
The authors would like to thank an anonymous referee for the very careful
reading of the paper and the valuable comments and corrections, in particular with respect to the proof of \Cref{generatedbyaset2variables}. 
The second named author is very grateful to Francesco Genovese for explaining how the notion of dg Bousfield localisation we consider gives rise to an adjunction of quasi-functors in the sense of \cite{adjunctions-quasi-functors-dg-categories} (see \Cref{equivalentdefinitionBousfield}). The authors would also like to thank Pieter Belmans, Boris Shoikhet, Greg Stevenson and Michel Van den Bergh for interesting discussions.

\section{The homotopy category of dg categories}\label{parhomandtensor}
We fix a commutative ground ring $k$ throughout the paper. 

Let $\fraku$ be a fixed (Grothendieck) universe. Without further notice, categories are $\fraku$-categories, small categories are $\fraku$-small categories and cocomplete categories are $\fraku$-cocomplete (i.e. have all $\fraku$-small colimits) etc.
In the sequel, making use of the universe axiom, we will sometimes use additional universes $\fraku \in \frakv$ and $\frakv \in \frakw$, which will be made explicit in the terminology and notation.

In this chapter, we revise the essential aspects of the homotopy theory of dg categories that will be used further on.

\subsection{The model structure on the category of dg categories}
We denote by $C(k) = \fraku - C(k)$ the category of cochain complexes of $\fraku$-small $k$-modules with cochain morphisms. The category $\dgcat_k = \fraku - \dgcat_{k}$ of $\fraku$-small dg categories over $k$ with $k$-linear dg functors has a standard model structure with the quasi-equivalences as weak equivalences \cite{structure-categorie-modeles-quillen-categorie-dg-categories}.
This model structure has the following properties.
\begin{proposition}{\cite[Prop 2.3]{homotopy-theory-dg-categories-derived-morita-theory}}\label{cofibrantreplacement}
	Consider $\dgcat_k$ with the standard model structure. The following hold:
	\begin{enumerate}
		\item Any object in $\dgcat_k$ is fibrant;
		\item There exists a cofibrant replacement $Q: \dgcat_k \lra \dgcat_k$ such that the natural morphism $Q(\aaa) \lra \aaa$ is the identity on objects;
		\item If $\aaa$ is cofibrant in $\dgcat_k$ and $A, A' \in \aaa$ then $\aaa(A,A')$ is cofibrant in $C(k)$ for the projective model structure.
	\end{enumerate}
\end{proposition}

We denote by $\hodgcat = \fraku - \hodgcat = \Ho(\fraku - \dgcat_{k})$ the homotopy category of $\fraku$-small dg categories. Given a dg functor $F: \aaa \lra \bbb$, we denote by $[F]$ its image in $\hodgcat$ and as usual we denote by $[-,-] = \fraku - [-,-] = \fraku - \hodgcat(-,-)$ the set of morphisms in $\hodgcat$. Observe that an element $f \in [\aaa,\bbb]$ induces a functor $H^0(f): H^0(\aaa) \lra H^0(\bbb)$ between the corresponding $H^0$-categories. 

\subsection{The monoidal structure on the homotopy category of dg categories}
Let $\ccc$ be a small dg category and $\dgmod(\ccc)$ the dg category of all dg modules (that is, dg functors from $\ccc^{\op}$ to $C(k)$). We denote by $\DD(\ccc)$ the \emph{dg derived category} of $\ccc$, that is the full dg subcategory $\DD(\ccc) \subseteq \dgmod(\ccc)$ of the cofibrant dg modules for the projective model structure on $\dgmod(\ccc)$ (see for example \cite[\S 3]{homotopy-theory-dg-categories-derived-morita-theory}, where the dg derived category of $\ccc$ is denoted by $Int(\ccc)$). By construction, $H^0(\DD(\ccc))$ is equivalent to the derived category $\D(\ccc)$ of $\ccc$ \cite[Prop 3.1]{differential-graded-categories}.

The homotopy category of dg categories $\hodgcat$ can be endowed with a closed symmetric monoidal structure, described by To\"en in \cite[\S 6]{homotopy-theory-dg-categories-derived-morita-theory}. In particular, given $\aaa, \bbb, \ccc$ small dg categories, in $\hodgcat$ we have the adjunction 
\begin{equation}\label{eq1}
	[\aaa \dertp \bbb, \ccc] \cong [\aaa, \RHom(\bbb, \ccc)].
\end{equation}
between the derived tensor product $\aaa \dertp \bbb$ and To\"en's internal $\RHom(\bbb, \ccc)$, which can be constructed as follows.

Let $\aaa$ and $\bbb$ be small dg categories. A bimodule $F \in \dgmod(\bbb \dertp \aaa^{\op})$ induces a dg functor $\Phi_F: \aaa \lra \dgmod(\bbb)$, and it is called \emph{right quasi-representable} provided that 
the induced $H^0(\Phi_F): H^0(\aaa) \lra H^0(\dgmod(\bbb))$ factors through a functor $H^0(F): H^0(\aaa) \lra H^0(\bbb)$. In other words, for all $A \in \aaa$, $\Phi_F(A) \in \dgmod(\bbb)$ is \emph{quasi-representable}, that is, quasi-isomorphic to a representable dg $\bbb$-module. We will denote by $\qrep(\bbb)$ the full dg subcategory of $\dgmod(\bbb)$ with as objects the quasi-representable objects. In particular, the dg Yoneda embedding $Y_{\bbb}: \bbb \lra \dgmod(\bbb)$ induces a quasi-equivalence $\bbb \lra \qrep(\bbb)$.

We denote by $\RHom(\aaa, \bbb) \subseteq \DD(\bbb \dertp \aaa^{\op})$ the full dg subcategory of (cofibrant) right quasi-representable bimodules. This category is not small, but essentially small, and hence can still be considered as an element of $\hodgcat$ (see \cite{homotopy-theory-dg-categories-derived-morita-theory}). In the literature, the elements of the category $H^0(\RHom(\aaa,\bbb))$ are usually called quasi-functors between $\aaa$ and $\bbb$ (see, for example \cite{differential-graded-categories}). Given $F \in \RHom(\aaa,\bbb)$, we denote the same element considered in $H^0(\RHom(\aaa,\bbb))$ also by $F$ and we will refer to both objects as \emph{quasi-functors}. 

In particular, the adjunction from (\ref{eq1}) above can easily be extended (see for example \cite[Cor 4.1]{internal-homs-via-extensions-dg-functors}) to the following isomorphism in $\hodgcat$:
\begin{equation}\label{eq2}
	\RHom(\aaa \dertp \bbb, \ccc) \cong \RHom(\aaa, \RHom(\bbb, \ccc)).
\end{equation}
Concretely, the isomorphism \eqref{eq2} is given by sending $F \in \RHom(\aaa \dertp \bbb, \ccc)$ to the associated dg functor
$$\aaa \lra \dgmod(\ccc \dertp \bbb^{\op}): A \longmapsto F_A$$ with $F_A(B,C) \coloneqq F(A,B,C)$. Then $F_A$ is right quasi-representable, and the resulting $\aaa \lra \RHom(\bbb, \ccc)$ gives rise to a representable element in $\RHom(\aaa, \RHom(\bbb, \ccc))$.

In addition, we have the following result, relating the morphisms in $\hodgcat$ and the internal hom of the monoidal structure.
\begin{proposition}[{{\cite[Cor 4.8]{homotopy-theory-dg-categories-derived-morita-theory}}}]\label{isoclasses}
	Let $\aaa$,$\bbb$ be two small dg categories. There exists a functorial bijection between the set $\left[ \aaa,\bbb \right] $ of maps between $\aaa$ and $\bbb$ in $\hodgcat$ and the set $\Iso(H^0(\RHom(\aaa,\bbb)))$ of isomorphism classes of quasi-functors.
\end{proposition}

Consider \emph{small} dg categories $\aaa$ and $\bbb$ and $F \in [\aaa, \bbb]$. By Yoneda's Lemma, if $F$ induces a bijection $F \circ -: [\ccc, \aaa] \cong [\ccc, \bbb]$ for every small dg category $\ccc$, it follows that $F$ is an isomorphism in $\hodgcat$. In the sequel, we will need the following refinement:

\begin{proposition}\label{propunivyon}
Consider dg $\fraku$-categories $\aaa$ and $\bbb$ and let $\frakv$ be a universe such that $\aaa$ and $\bbb$ are $\frakv$-small. We consider the homotopy category $\frakv - \hodgcat$ of $\frakv$-small dg categories and $F \in \frakv - [\aaa, \bbb]$. If $F$ induces a bijection $F \circ -: \frakv - [\ccc, \aaa] \cong \frakv - [\ccc, \bbb]$ for every $\fraku$-small dg category $\ccc$, it follows that $F$ is an isomorphism in $\frakv - \hodgcat$.
\end{proposition}

\begin{proof}
We may suppose that $F$ is given by a dg functor $F: \aaa \lra \bbb$. Suppose that $F$ induces a bijection $F\circ -: \frakv - [\ccc, \aaa] \cong \frakv - [\ccc, \bbb]$ for every $\fraku$-small dg category $\ccc$. We are to show that $F$ is a quasi-equivalence. 

We start by showing that $F$ is quasi-essentially surjective. Consider the dg category $\underline{k}$ with a single object $\ast$ and $\underline{k}(\ast, \ast) = k$. It is readily seen that there is a natural quasi-equivalence $\aaa \cong \frakv - \RHom(\underline{k}, \aaa)$ for every $\frakv$-small dg category $\aaa$ and hence by Proposition \ref{isoclasses} a natural bijection $\frakv - [\underline{k}, \aaa] \cong \Iso (H^0(\aaa))$. Hence, by the assumption (for $\ccc = \underline{k}$) $F$ induces a bijection $\Iso (H^0(\aaa)) \lra \Iso (H^0(\bbb))$ as desired.

Next we show that $F$ is quasi-faithful. Consider $$H^n(F_{A, A'}): H^n \aaa(A, A') \lra H^n \bbb(F(A), F(A'))$$ and $f \in Z^n \aaa(A, A')$ with $H^n(F_{A, A'})([f]) = 0 \in H^n \bbb(F(A), F(A'))$.
Consider the dg category $\mathrm{Ar}_n$ with two objects $X, X'$ and $\mathrm{Ar}_n(X, X) = k1_X$, $\mathrm{Ar}_n(X', X') = k1_{X'}$, $\mathrm{Ar}_n(X, X') = kx$ for $x$ in degree $n$, $\mathrm{Ar}_n(X', X) = 0$. Consider the dg functor $\phi: \mathrm{Ar}_n \lra \aaa: x \longmapsto f$. We have $F\phi(x) = d(h)$ for some $h \in \bbb(F(A), F(A'))^{n-1}$. Consider the dg functors $\psi_1: \mathrm{Ar}_n \lra \aaa: x \longmapsto 0_{A, A'}$ and $\psi_2: \mathrm{Ar}_n \lra \bbb: x \longmapsto 0_{F(A), F(A')}$ for the zero morphisms $0_{A, A'} \in \aaa(A,A')^n$ and $0_{F(A), F(A')} \in \bbb(F(A), F(A'))^n$. We claim that $[F \phi] = [F][\phi] = [\psi_2]$ in $[\mathrm{Ar}_n, \bbb]$. Let $\ppp(\bbb)$ be the path object dg category for $\bbb$ as described in \cite[\S 2.2]{internal-homs-via-extensions-dg-functors}. Then it is readily seen that a homotopy between $F \phi$ and $\psi_2$ is given by 
$$H: \mathrm{Ar}_n \lra \ppp(\bbb)$$
with 
\begin{equation*}
	\begin{aligned}
		 H(X) &= (F(A),F(A),1_{F(A)})\\
		H(X') &= (F(A'),F(A'), 1_{F(A')})\\
		H(x) &= (F(f),0_{F(A),F(A')},(-1)^{n-1} h)
	\end{aligned}
\end{equation*} 
Since also $[F \psi_1] = [F][\psi_1] = [\psi_2]$ it follows from the assumption (for $\ccc = \mathrm{Ar}_n$) that $[\phi] = [\psi_1] \in [\mathrm{Ar}_n, \aaa]$ and consequently $[f] = 0 \in H^n \aaa(A, A')$ as desired.

Finally we show that $F$ is quasi-full. Thanks to the bijection $\Iso (H^0(\aaa)) \lra \Iso (H^0(\bbb))$, it suffices to show that for all $B, B' \in \bbb$, there exist $A, A' \in \aaa$ and isomorphisms $B \cong F(A)$ and $B' \cong F(A')$ in $H^0 \bbb$ such that $H^n\aaa(A,A') \lra H^n \bbb(F(A), F(A'))$ is an isomorphism for every $n$.
So let $B, B' \in \bbb$. Consider the full dg subcategory $\iota: \bbb_0 \subseteq \bbb$ spanned by the objects $B$ and $B'$ and let $Q: Q(\bbb_0) \lra \bbb_0$ be a cofibrant resolution which is the identity on objects. By the assumption (for $\ccc = Q(\bbb_0)$), there exists a dg functor $G: Q(\bbb_0) \lra \aaa$ with $[F][G] = [FG] = [\iota Q] \in [Q(\bbb_0), \bbb]$. It follows that $$H^n(F_{G(B), G(B')}): H^n\aaa(G(B), G(B')) \lra H^n\bbb(F(G(B)), F(G(B')))$$ is surjective as desired.
\end{proof}

\subsection{Variations upon the inner hom}

Consider dg $\fraku$-categories $\aaa$ and $\bbb$. For universes $\frakv \subseteq \frakv'$ such that $\aaa$ and $\bbb$ are $\frakv$-small, there is easily seen to be a quasi-equivalence $\frakv{-}\RHom(\aaa, \bbb) \cong \frakv'{-}\RHom(\aaa, \bbb)$. Hence, we will often omit the decoration $\frakv$ from the notation and simply write $\RHom(\aaa, \bbb)$ where it is understood that we make use of some universe for which the categories under considerations are small. If $\aaa$ is $\fraku$-small, then $\RHom(\aaa, \bbb)$ is seen to be a dg $\fraku$-category.

For $F \in \RHom(\aaa, \bbb)$, we have an induced functor $H^0(F): H^0(\aaa) \lra H^0(\bbb)$. We will consider several full subcategories of $\RHom(\aaa, \bbb)$ determined by properties of the functors $H^0(F)$.

Given a universe $\fraku$, its cardinality $|\fraku|$ is the unique inaccessible (and hence regular) cardinal such that $\fraku = V_{|\fraku|}$ where, for a cardinal $\kappa$, $V_{\kappa} =\{ X \,\,|\,\, |X| < \kappa \}$ - consisting of all the $\kappa$-small sets - denotes the $\kappa^{\mathrm{th}}$-level of the von Neumann hierarchy (see \cite{grothendieck-universes}). Observe that, for $\fraku \in \frakv$, we have that $|\fraku| < |\frakv|$ and hence $|\fraku|$ is a $|\frakv|$-small cardinal.

\begin{definition}
Let $\ccc$ be a dg $\fraku$-category.
\begin{enumerate}
\item Let $\alpha$ be a cardinal. We say that $\ccc$ is \emph{homotopically $\alpha$-cocomplete} if $H^0(\aaa)$ has all $\alpha$-small coproducts.
\item We say that $\ccc$ is \emph{homotopically cocomplete} if $\ccc$ is homotopically $|\fraku|$-cocomplete, that is, $H^0(\aaa)$ has all $\fraku$-small coproducts.
\end{enumerate}
\end{definition}

\begin{definition}
Consider dg $\fraku$-categories $\aaa$ and $\bbb$. 
\begin{enumerate}
\item Let $\alpha$ be a cardinal. A quasi-functor $F \in \RHom(\aaa, \bbb)$ is called \emph{$\alpha$-cocontinuous} if the induced functor $H^0(F): H^0(\aaa) \lra H^0(\bbb)$ preserves all $\alpha$-small coproducts. We let 
$$\RHom_{\alpha}(\aaa, \bbb) \subseteq \RHom(\aaa, \bbb)$$ denote the full dg subcategory of $\alpha$-cocontinuous quasi-functors.
\item A quasi-functor $F \in \RHom(\aaa, \bbb)$ is called \emph{cocontinuous} if it is $|\fraku |$-cocontinuous, that is if the induced functor $H^0(F): H^0(\aaa) \lra H^0(\bbb)$ preserves all $\fraku$-small coproducts. We put
$$\RHom_{\ccm}(\aaa, \bbb) =  \RHom_{|\fraku|}(\aaa, \bbb).$$ 
\end{enumerate}
\end{definition}
Next we look at annihilation of classes of objects.

\begin{definition}
Consider dg categories $\aaa, \bbb$ and let $\nnn \subseteq \mathrm{Ob}(\aaa)$ be a class of objects. We say that $F \in \RHom(\aaa,\bbb)$ \emph{annihilates $\nnn$} if the induced functor $H^0(F): H^0(\aaa) \lra H^0(\bbb)$ is such that $H^0(F)(N) = 0$ for every $N \in \nnn$.
 We denote by
$$\RHom_{\nnn}(\aaa, \bbb) \subseteq \RHom(\aaa, \bbb)$$
the full dg subcategory of quasi-functors annihilating $\nnn$. 
\end{definition}

\begin{remark}
	We will use the same terminology and notation for a full dg subcategory $\aaa' \subseteq \aaa$, where it is understood that annihilation is intended with respect to the class $\nnn = \mathrm{Ob}(\aaa')$.
\end{remark}

The \emph{dg quotient} $\bbb/\aaa$ of a dg category $\bbb$ along a full dg subcategory $\aaa \subseteq \bbb$ was introduced by Keller in \cite{on-the-cyclic-homology-exact-categories} and analysed further by Drinfeld in \cite{dg-quotients-dg-categories}. The dg quotient fulfills the following universal property in $\hodgcat$:
\begin{equation}\label{kellerquotient}
\RHom(\bbb/\aaa,\ccc) \cong \RHom_{\aaa}(\bbb,\ccc),
\end{equation}
for all $\ccc \in \hodgcat$ (see \cite{drinfeld-dg-quotient}). 

\begin{example}
Let $\ccc$ be a small dg category and let $\Ac_{\mathrm{dg}}(\ccc)$ be the full dg subcategory of $\dgmod(\ccc)$ of acyclic dg modules, that is, the dg modules which are pointwise acyclic.%
The natural composition of morphisms in $\hodgcat$
\begin{equation}\label{interpretationofdgderived}
\DD(\ccc) \lra \dgmod(\ccc) \lra \dgmod(\ccc)/ \Ac_{\mathrm{dg}}(\ccc)
\end{equation}
is an isomorphism, and hence it induces a morphism $Q \in [\dgmod(\ccc), \DD(\ccc)]$. 
\end{example}

\section{Well generated dg categories}\label{wellgenerateddgcategories}
Well-generated triangulated categories in the sense of Neeman \cite{triangulated-categories} form a very important class of triangulated categories. They enjoy very nice properties concerning for example localisations (see \cite{localization-theory-triangulated-categories}) and Brown representability (see \cite[\S 8.4]{triangulated-categories}), and they also appear naturally in many contexts.
In particular, derived categories of Grothendieck abelian categories are well generated triangulated \cite{derived-category-sheaves-manifold}.

Porta shows in \cite{the-popescu-gabriel-theorem-for-triangulated-categories} that in the triangulated world, well generated algebraic triangulated categories play  the analogous role to the one that Grothendieck categories play in the abelian world, in the sense that they fulfill a triangulated version of the well-known Gabriel-Popescu theorem for Grothendieck categories \cite{caracterisation-categories-abeliennes-generateurs-limites-inductives-exactes}. 

In this article we will focus on the pretriangulated dg version of well generated algebraic triangulated categories:
\begin{definition}\label{defwellgen}
	A pretriangulated dg category $\aaa$ is called \emph{well generated} if the homotopy category $H^0(\aaa)$ is a well generated triangulated category. It is called \emph{$\alpha$-compactly generated} for some cardinal $\alpha$ if $H^0(\aaa)$ is $\alpha$-compactly generated.
\end{definition}

Observe that in Definition \ref{defwellgen}, $H^0(\aaa)$ is automatically algebraic as it has $\aaa$ as an enhancement.
\begin{remark}
	From now on, when dealing with well generated pretriangulated dg categories, we will usually omit the term pretriangulated for the sake of brevity.
\end{remark}

In section \S \ref{equivapproaches} we discuss the localisation theory of well generated dg categories, which can be obtained as an enhancement of the localisation theory of well generated triangulated categories as described for example in \cite{localization-theory-triangulated-categories} (see \S \ref{triangulatedloc}). After recalling $\alpha$-cocontinuous (dg) derived categories in \S \ref{dgderivedcategories}, in \S \ref{parenhGP} we formulate the (enhanced) derived Gabriel-Popescu theorem due to Porta \cite{the-popescu-gabriel-theorem-for-triangulated-categories}. In \S \ref{parcocontint}, we prove that the cocontinuous internal hom between homotopically cocomplete dg categories is again homotopically cocomplete (Theorem \ref{main1}). In \S \ref{parwellint}, we prove the main result of this chapter: the cocontinuous internal hom between well generated dg categories is again well generated (Theorem \ref{closed}).   

\subsection{Localisation of well generated triangulated categories}\label{triangulatedloc}
The Verdier quotient of a triangulated category $\ttt$ with respect to a full triangulated subcategory $\www$ is given by a triangulated category $\ttt/\www$ and an exact functor $Q: \ttt \lra \ttt/\www$ annihilating $\www$ such that any exact functor $\ttt \lra \ttt'$ annihilating $\www$ factors through $Q$. In other words, we have that the Verdier quotient $\ttt/\www$ has the following universal property:
\begin{equation}
	\Fun_{\Tri}(\ttt/\www,\ttt') \underset{\cong}{\xrightarrow{-\circ Q}} \Fun_{\Tri,\www}(\ttt,\ttt'),
\end{equation}
where $\Fun_{\Tri}(\ttt/\www,\ttt')$ denotes the collection of exact functors from $\ttt/\www$ to $\ttt'$ and $\Fun_{\Tri,\www}(\ttt,\ttt')$ denotes the collection of exact functors from $\ttt$ to $\ttt'$ that annihilate $\www$.
On the other hand, a Bousfield localisation functor $L:\ttt \lra \ttt$ can be characterized as the composition of a Verdier quotient $Q: \ttt \lra \ttt/\Kern(L)$ followed by its right adjoint $\ttt/\Kern(L) \lra \ttt$ (see \cite{localization-theory-triangulated-categories}).

If we restrict to the realm of well generated triangulated categories, we have that localising subcategories of well generated categories which are generated by a set of objects are again well generated, and so are the corresponding Verdier quotients \cite[Thm 7.2.1]{localization-theory-triangulated-categories}. Then, we have two equivalent approaches to the localisation of well generated triangulated categories which produce again well generated triangulated categories and which are equivalent, namely:
\begin{itemize}
	\item Verdier quotients along localising subcategories generated by a set; 
	\item Bousfield localisations with kernel generated by a set; 
\end{itemize}
where we say that a localising subcategory $\www$ of a well generated triangulated category $\ttt$ is generated by a set if there exists a set of objects of $\ttt$ such that the smallest localising subcategory containing them is $\www$. 
The fact that these two approaches are equivalent can be directly deduced from \cite[Thm 7.2.1 \& Prop 5.2.1]{localization-theory-triangulated-categories}.

In what follows, we analyse the induced correspondence of localisation theories in the dg setting. But before we proceed, we make an observation on the universal properties of the Verdier and dg quotients when we restrict to the well generated case with cocontinuous functors. 

Let $\ttt$ be a well generated triangulated category and $\www \subseteq \ttt$ a localising subcategory generated by a set. One can easily observe that under this hypothesis the quotient functor $Q: \ttt \lra \ttt/\www$ preserves coproducts, as it is a left adjoint between well generated triangulated categories. It is then not hard to check that the Verdier quotient, restricted to well generated triangulated categories, has the following universal property. Given $\ttt$ a well generated triangulated category, and $\www \subseteq \ttt$ a localising subcategory generated by a set of objects (and hence well generated), the Verdier quotient $\ttt/\www$ is a well generated triangulated category such that for any well generated triangulated category $\ttt'$, one has that
\begin{equation}
\Fun_{\Tri,\ccm}(\ttt/\www,\ttt') \underset{\cong}{\xrightarrow{-\circ Q}} \Fun_{\Tri,\ccm,\www}(\ttt,\ttt'),
\end{equation}
where the subindex $\ccm$ indicates that we are considering the exact functors which preserve coproducts. 

In the dg realm one can check in a similar fashion, for example by means of Keller's construction, that if $\bbb$ is a well generated dg category and $\aaa \subseteq \bbb$ is a dg subcategory with $H^0(\aaa)$ localising in $H^0(\bbb)$ and generated by a set, then the dg quotient $\bbb/\aaa$ is also a well generated dg category (as it is an enhancement of the Verdier quotient $H^0(\bbb)/H^0(\aaa)$) and the canonical morphism $Q: \bbb \lra \bbb/\aaa$ in $\hodgcat$ is cocontinuous, that is, the induced $H^0(\bbb) \lra H^0(\bbb/\aaa)$ preserves coproducts. Observe then, that for all well generated dg categories $\ccc$, the universal property of the dg quotient \eqref{kellerquotient} in $\hodgcat$ restricts to a quasi-equivalence 
\begin{equation}\label{dgquotientrestricted}
\RHom_{\ccm}(\bbb/\aaa,\ccc) \cong \RHom_{\ccm,\aaa}(\bbb,\ccc).
\end{equation}
\subsection{Localisation of well generated dg categories} \label{equivapproaches}
\subsubsection{Localising subcategories generated by a set}
Let $\bbb$ be a well generated dg category. Observe that in particular $H^0(\bbb)$ is localising as a subcategory of itself and it is, as localising subcategory, generated by a set. In addition, the intersection of localising subcategories of $H^0(\bbb)$ generated by a set is again such (see \cite[Lem 3.2]{bousfield-lattice-triangulated-category-stratification}). Consequently, for every full triangulated subcategory $\hhh \subseteq H^0(\bbb)$ there is a smallest localising subcategory generated by a set containing $\hhh$. In particular, the poset of localising subcategories of $H^0(\bbb)$ generated by a set is a complete lattice with $\inf_i \www_i = \cap_i \www_i$ and $\sup_i \www_i = \langle \cup_i \www_i \rangle$, where $\langle \cup_i \www_i \rangle$ denotes the smallest localising subcategory that contains $\cup_i \www_i$. Observe that $\langle \cup_i \www_i \rangle$ is indeed generated by a set, taking for example $\cup_i \nnn_i$ where, for every $i$, $\nnn_i$ is a set such that $\langle \nnn_i \rangle = \www_i$.
\begin{definition}
	Consider $\hhh \subseteq H^0(\bbb)$ and $B \in \bbb$. A \emph{filtration} of $B$ consists of a countable collection $(X_i)_{i=0}^{\infty}$ of objects in $H^0(\bbb)$ with $X_0=0$ and maps $x_i: X_{i} \lra X_{i+1}$ for all $i \geq 0$ such that $\hocolim(X_i) = B$. A filtration $(X_i)_{i=0}^{\infty}$ of $B$ is called an $\hhh$-filtration if the cone of each $x_i: X_{i} \lra X_{i+1}$ belongs to $\hhh$ and in this case $B$ is called \emph{$\hhh$-filtered}.
\end{definition}

\begin{proposition}\label{filtration}
	Let $\www$ be a localising subcategory of $H^0(\bbb)$ generated by a set. Then, there exists a set $\nnn$ generating $\www$ (i.e. $\www = \langle \nnn \rangle$) such that $X \in H^0(\bbb)$ belongs to $\www$ if and only if it is $\overline{\nnn}$-filtered, where $\overline{\nnn}$ is the class of small coproducts of elements in $\nnn$.
	\begin{proof}
		By \cite[Thm 7.2.1]{localization-theory-triangulated-categories}, we know we can take a regular cardinal $\alpha$ such that $\www$ and $H^0(\bbb)$ are both $\alpha$-compactly generated. In particular, the class of $\alpha$-compact objects $\www^{\alpha}= \www \cap \bbb^{\alpha}$ is essentially small (see \cite[Prop 3.2.5, Lem 4.4.5]{triangulated-categories}). Take $\nnn$ to be the set of objects in $\www$ consisting of taking for each isomorphism class of $\www^{\alpha}$ a representative. We have that $\www = \langle \nnn \rangle$. By applying \cite[Lemma B.1.3]{triangulated-categories} to $\www$, we know that every $X \in \www$ is $\overline{\nnn}$-filtered. On the other hand, as $\www$ is localising, every $\overline{\nnn}$-filtered object $X$ in $H^0(\bbb)$ belongs to $\www$, which concludes the argument.  
	\end{proof}
\end{proposition}
We describe now the relation with orthogonal complements. 

Let $\ttt$ be a triangulated category. Recall that an object $X\in \ttt$ is said to be \emph{left orthogonal} to an object $Y \in \ttt$ (or $Y$ \emph{right orthogonal} to $X$) if $\ttt(X,Y) = 0$ and we denote this by $X \perp Y$. For a full subcategory $\hhh \subseteq \ttt$, we obtain the following $k$-linear subcategories of $\ttt$:
\begin{itemize}
	\item $\hhh^{\perp} = \{X \in \ttt \,\,|\,\, H \perp X \,\, \text{ for all }H \in \hhh\}$
	\item ${^{\perp}}\hhh = \{X \in \ttt \,\,|\,\, X \perp H \,\, \text{ for all }H \in \hhh\}$
\end{itemize}
\begin{remark}
	This notation for the right and left orthogonals is the most common in the literature, though it is not standard. For example, the notation in \cite{triangulated-categories} is reversed (see \cite[Def 9.1.10 \& 9.1.11]{triangulated-categories}).
\end{remark}
\begin{proposition}\label{leftorth}
	Let $\www$ be a localising subcategory of $\bbb$ generated by a set $\nnn$, i.e. $\www = \langle \nnn \rangle$. Then we have that $\www^{\perp} = \nnn^{\perp}$.
	\begin{proof}
		We have that $\nnn \subseteq \www$, hence $\www^{\perp} \subseteq \nnn^{\perp}$. On the other hand, we have that $\nnn \subseteq \,^{\perp}(\nnn ^{\perp})$ and $\,^{\perp}(\nnn ^{\perp})$ is easily seen to be a localising (hence triangulated) subcategory \cite[Lem 9.1.12]{triangulated-categories}. Hence we have that $\www = \langle \nnn \rangle \subseteq \,^{\perp}(\nnn ^{\perp})$. Then, applying right orthogonals and the fact that $\nnn^{\perp} = (\,^{\perp}(\nnn ^{\perp}))^{\perp}$, we obtain that $\nnn^{\perp} \subseteq \www^{\perp}$, which concludes the argument.
	\end{proof}
\end{proposition}

\subsubsection{Bousfield localisations}
\begin{definition}
	Given two pretriangulated dg categories $\aaa$, $\bbb$ and two right quasi-representable functors $F \in \RHom(\aaa,\bbb)$, $G \in \RHom(\bbb,\aaa)$, we say that $F$ is \emph{left quasi-adjoint} to $G$ if and only if $H^0(F) \dashv H^0(G)$. In this case we write $F \dashv_{H^0} G$.
\end{definition}

\begin{definition}\label{defbfloc}
	Let $\aaa, \bbb$ be pretriangulated dg categories and $i: \bbb \lra \aaa$ a quasi-fully faithful dg functor. We say that $i \in \RHom(\bbb, \aaa)$ is a \emph{dg Bousfield localisation} of $\aaa$ if $H^0(i): H^0(\bbb) \hookrightarrow H^0(\aaa)$ admits a left adjoint. 
\end{definition}

\begin{remark}\label{equivalentdefinitionBousfield}
	This definition is seen to be equivalent to the following definition: $i: \bbb \lra \aaa$ is a \emph{dg Bousfield localisation} if and only if there exists a right quasi-representable functor $a \in \RHom(\aaa, \bbb)$ which is left adjoint to $i \in \RHom(\bbb,\aaa)$ in the sense of adjoint pairs of quasi-functors from \cite{adjunctions-quasi-functors-dg-categories}. Obviously, this second definition implies the first.
	On the other hand, if $H^0(i)$ has a left adjoint $F: H^0(\aaa) \lra H^0(\bbb)$, then we have an isomorphism
	$$H^0(\bbb)(F(A),-) \cong H^0(i(A,-))$$
	for all $A \in \aaa$, where we consider $i$ in $\dgmod(\aaa \dertp \bbb^{\op})$. This isomorphism is, by Yoneda lemma, determined by an element $f \in H^0(i(A,F(A)))$. Consider $g$ a closed element of degree 0 in $i(A,F(A))$ lifting $f$. By dg Yoneda lemma, $g$ induces a morphism
	$$\bbb(F(A),-) \lra i(A,-)$$
	which is a quasi-isomorphism because it is a lift of the previous $0^{\text{th}}$-cohomology isomorphim, and both $\aaa$ and $\bbb$ are pretriangulated. This shows that $i$ is left quasi-representable as a bimodule and hence it admits a left adjoint $a \in \RHom(\aaa,\bbb)$ as a consequence of \cite[Prop 7.1]{adjunctions-quasi-functors-dg-categories}. In particular, by unicity of adjoints, we have that $H^0(a) \cong F$. 
	
	Observe this implies, in particular, that dg Bousfield localisations have left quasi-adjoints.
\end{remark}

\begin{remark} \label{unit}
	Fix the same notations as in \Cref{equivalentdefinitionBousfield}. As a direct consequence of the theory of adjunctions of quasi-functors from \cite[\S 6]{adjunctions-quasi-functors-dg-categories}, there exist morphisms $\id_{\aaa} \lra i \dertp_{\bbb} a$ in $H^0(\RHom(\aaa,\aaa)) \subseteq \D(\aaa \otimes \aaa^{\op})$ and $a \dertp_{\aaa} i \lra \id_{\bbb}$ in $H^0(\RHom(\bbb,\bbb)) \subseteq \D(\bbb \otimes \bbb^{\op})$, called the unit and counit of the adjunction respectively, where $\dertp$ is the composition of bimodules, which preserves right quasi-representability (see \cite[\S6.1]{deriving-dg-categories}). Observe that in our particular situation the counit $a \dertp_{\aaa} i \lra \id_{\bbb}$ is an isomorphism in $H^0(\RHom(\aaa,\aaa))$ and hence $a \dertp_{\aaa} i$ and $\id_{\bbb}$ are quasi-isomorphic in $\RHom(\bbb,\bbb)$. Moreover, notice that $a$ is cocontinuous, i.e. it belongs to $\RHom_{\ccm}(\aaa,\bbb)$.
\end{remark}

\begin{remark}\label{remqlpretriang}
	Observe that a dg Bousfield localisation induces a classical Bousfield localisation of the corresponding underlying triangulated category.
\end{remark}  

\subsubsection{Equivalent approaches to localisation}\label{posetisom}
When we restrict to the world of well generated triangulated categories, there is a nice correspondence between localising subcategories and Bousfield localisation, as we have pointed out at the beginning of \S\ref{triangulatedloc}. This result can be easily enhanced to the dg realm. In particular, for a well generated dg category $\bbb$, there is a poset isomorphism between:
\begin{enumerate}
	\item The poset $W_{\mathsf{dg}}$ of localising subcategories of $H^0(\bbb)$ generated by a set, ordered by inclusion;
	\item The opposite poset $(L_{\mathsf{dg}})^{\op}$ of the poset $L_{\mathsf{dg}}$ of Bousfield localisations of $\bbb$ with kernel of the left adjoint (at the $0^{\text{th}}$-cohomology level) generated by a set, ordered by inclusion, i.e. we write $i \subseteq i'$ if and only if $\Img(i) \subseteq \Img(i')$ as sub-dg-categories, where $\Img(i)$ denotes the quasi-essential image of $i$.
\end{enumerate} 
The poset isomorphism is described as follows:
\begin{enumerate}
	\item Let $\www$ be a localising subcategory of $H^0(\bbb)$ generated by a set. In particular, we have that $\www^{\perp} \subseteq H^0(\bbb)$ has a left adjoint and hence gives rise to a localisation functor
	$$H^0(\bbb) \lra \www^{\perp} \lra H^0(\bbb),$$ 
	such that the composition $\www^{\perp} \hookrightarrow H^0(\bbb) \rightarrow H^0(\bbb)/\www$ is an equivalence and $\www^{\perp}$ is well generated (see \cite[Prop 7.2.1, Prop 5.2.1 \& Prop. 4.9.1]{localization-theory-triangulated-categories}).
	
	Denote by $\LLL_{\www}$ the full dg subcategory of $\bbb$ obtained as an enhancement of $\www^{\perp} \subseteq H^0(\bbb)$ via the natural enhancement of $H^0(\bbb)$. We have that $\LLL_{\www}$ is a well generated dg category, and that $F: H^0(\bbb) \lra H^0(\bbb)/ \www \cong \www^{\perp}$ is a left adjoint of $H^0(i): H^0(\LLL_{\www}) \subseteq H^0(\bbb)$, where $i$ denotes the embedding $\LLL_{\www} \subseteq \bbb$. In addition, $\Kern(F)= \www$, which is generated by a set of objects. 
	
	To each $\www \in W_{\mathsf{dg}}$ we assign the so constructed $\LLL_{\www} \in L_{\mathsf{dg}}$.
	\item Let $i: \bbb \lra \aaa$ be a Bousfield localisation of a well generated dg category $\aaa$ such that the kernel of the left adjoint $F$ of $H^0(i)$ is generated by a set of objects. Observe that $\Kern(F)$ is a localising subcategory of $H^0(\bbb)$. We put $\www_{\LLL}= \Kern(F)$.
	
	We assign to $\LLL \in L_{\mathsf{dg}}$ the so constructed $\www_{\LLL} \in W_{\mathsf{dg}}$.
\end{enumerate}

\subsection{The $\alpha$-cocontinuous derived category}\label{dgderivedcategories}
In this section we recall the \emph{$\alpha$-cocontinuous derived category} of an \emph{$\alpha$-cocomplete} dg category from \cite{the-popescu-gabriel-theorem-for-triangulated-categories} (note that in loc. cit it is called the ``$\alpha$-continuous derived category''). 

\begin{definition} {\cite[\S 6]{the-popescu-gabriel-theorem-for-triangulated-categories}}
Let $\ccc$ be a homotopically $\alpha$-cocomplete small dg category. The \emph{$\alpha$-cocontinuous derived category} $\D_{\alpha}(\ccc)$ is defined as the full subcategory of $\D(\ccc)$ with objects given by the dg functors $X$ such that for every $\alpha$-small family of objects $\{A_i\}_{i \in I}$ the canonical morphism
	\begin{equation}
	H^n(X) \left( \coprod_i^{H^0(\ccc)} A_i \right) \lra \prod_i H^n(X) (A_i)
	\end{equation}
	is invertible for all $n \in \mathbb{Z}$, where $\underset{i}{\overset{H^0(\ccc)}{\coprod}} A_i $ denotes the coproduct taken in $H^0(\ccc)$.
\end{definition}
\begin{remark}
	Observe that, in particular, the representable dg modules belong to $\D_{\alpha}(\ccc)$.
\end{remark}
\begin{remark}\label{remDalpha}
	In addition, one can give an equivalent definition of $\D_{\alpha}(\ccc)$ as a Verdier quotient of $\D(\ccc)$ with respect to the localising subcategory $\nnn$ generated by the cones of the morphisms
\begin{equation}
	\{\sigma_{\lambda}: \coprod_{i\in I} h_{C_i} \rightarrow  h_{\coprod_{i\in I}^{H^0(\ccc)} C_i} \}_{\lambda},
\end{equation}
where $\lambda$ varies in the set of all $\alpha$-small families $\{C_i\}_{i\in I}$ of objects of $\ccc$.
\end{remark}
\begin{definition}
	We call the natural enhancement of $\D_{\alpha}(\ccc)$ via the enhancement $\DD(\ccc)$ of $\D(\ccc)$ the \emph{$\alpha$-cocontinuous derived dg category} of $\ccc$. We will denote it by $\DD_{\alpha} (\ccc)$. 
\end{definition}
There is an equivalent construction of $\DD_{\alpha}(\ccc)$ in $\hodgcat$ as a dg quotient. Indeed, we have that for the dg quotient $\DD(\ccc) / \nnn'$, where $\nnn'$ is the natural enhancement of $\nnn$ above via the enhancement $\DD(\ccc)$ of $\D(\ccc)$, the natural composition of morphisms in $\hodgcat$
\begin{equation}
\DD_{\alpha}(\ccc) \lra \DD(\ccc) \lra \DD(\ccc)/ \nnn'
\end{equation}
is an isomorphism. This induces a morphism $Q_{\alpha}\in [\DD(\ccc), \DD_{\alpha}(\ccc)]$.
In particular, we have the following
\begin{theorem}[{{\cite[Thm 6.4]{the-popescu-gabriel-theorem-for-triangulated-categories}}}]
	Let $\AAA$ be a homotopically $\alpha$-cocomplete small dg category. Then $\D_{\alpha}(\AAA)$ is $\alpha$-compactly generated by the images of the free dg modules $\{\AAA(-,A)\}_{A \in \AAA}$ through the localisation functor $\D(\AAA) \lra \D_{\alpha}(\AAA)$.
\end{theorem}

\subsection{Enhanced derived Gabriel-Popescu theorem}\label{parenhGP}

In \cite{the-popescu-gabriel-theorem-for-triangulated-categories}, 
Porta proved a derived version of the Gabriel-Popescu theorem, showing that a triangulated category $\ttt$ is well generated and algebraic if and only if there exists a small dg category $\AAA$ such that $\ttt$ is triangle equivalent to the Verdier quotient of $\D(\AAA)$ by a localising subcategory generated by a set. 
Further, $\ttt$ is $\alpha$-compactly generated and algebraic if and only if there exists a small homotopically $\alpha$-cocomplete dg category $\AAA$ such that $\ttt$ is triangle equivalent to $\D_{\alpha}(\AAA)$. 

We are interested in enhanced versions of these results, which can easily be deduced making use of the higher observations (see also  \cite{uniqueness-dg-enhancements-derived-category-grothendieck-category}).

\begin{theorem}\label{thmenhGP}
Let $\ccc$ be a pretriangulated dg category.
\begin{enumerate}
\item $\ccc$ is well generated if and only if there exists a small dg category $\AAA$ such that $\ccc \cong \DD(\AAA)/\www$ in $\hodgcat$, where $\www$ is the enhancement of a localising subcategory of $\D(\AAA)$ generated by a set.
\item $\ccc$ is $\alpha$-compactly generated if and only if there exists a small homotopically $\alpha$-cocomplete dg category $\AAA$ such that $\ccc \cong  \DD_{\alpha}(\AAA)$ in $\hodgcat$.
\end{enumerate}
\end{theorem}

From Theorem \ref{thmenhGP}, one deduces (see \cite[\S 3.1]{derived-azumaya-algebras-generators-twisted-derived-categories}):

\begin{corollary}\label{corwglocpres}
Let $\ccc$ be a pretriangulated dg category. Then $\ccc$ is well generated if and only if $\ccc$ is locally presentable in the sense of \cite{derived-azumaya-algebras-generators-twisted-derived-categories}.
\end{corollary}

\subsection{The cocontinuous internal hom of homotopically cocomplete dg categories}\label{parcocontint}
In this section we prove that given a $\fraku$-small dg category $\BBB$, and a well generated $\frakv$-small dg category $\ccc$ with a $\fraku$-small set of generators, the internal hom $\RHom(\BBB,\ccc)$ in $\frakv-\hodgcat$ is a well generated dg category as well. As a consequence of this result, we  prove that for any two $\fraku$-small dg categories $\AAA, \BBB$ with $\BBB$ homotopically $\fraku$-cocomplete (resp. $\alpha$-cocomplete) the internal hom $\RHom(\AAA,\BBB)$ is also homotopically $\fraku$-cocomplete (resp. $\alpha$-cocomplete) in $\fraku-\hodgcat$, while if also $\AAA$ is homotopically $\fraku$-cocomplete (resp. $\alpha$-cocomplete), then so is $\RHom_{\ccm}(\AAA,\BBB)$ (resp. $\RHom_{\alpha}(\AAA,\BBB)$) in $\fraku- \hodgcat$.

We will start first with some considerations on the two variable setting.

The fact that the cofibrant replacement $Q$ in $\dgcat_k$ can be taken to be the identity on objects, permits to define a canonical functor
$$i_B : \AAA \lra \AAA \dertp \BBB = \AAA \otimes Q(\BBB): A \longmapsto (A,B)$$
for all $B\in \BBB$ (see \cite[\S 4]{homotopy-theory-dg-categories-derived-morita-theory}). 

One can then consider the induced dg functor 
$$(i_B)^*: \dgmod(\AAA \dertp \BBB) \lra \dgmod(\AAA): F \longmapsto F\circ i_B = F(-,B),$$
sometimes called \emph{restriction of scalars}. This dg functor has a left adjoint $$(i_B)_!: \dgmod(\AAA) \lra \dgmod(\AAA \dertp \BBB),$$ sometimes also called \emph{extension of scalars}. Moreover, $(i_B)^*$ preserves acyclic dg modules, hence it induces an exact functor
$$(i_B)^*: \D(\AAA \dertp \BBB) \lra \D(\AAA)$$
In addition, the left derived functor 
$$\mathrm{L}(i_B)_{!}: \D(\AAA) \lra \D(\AAA \dertp \BBB).$$
is a left adjoint for $(i_B)^*$ (see \cite[\S 1]{uniqueness-enhancement-triangulated-categories}). Observe our notations for the restriction and extension of scalars functors follow the convention from classical topos theory as in \cite{SGA4-1} while in loc.cit. another convention is used.   

\begin{lemma}\label{adjunctiontwovariables}
	Let $\AAA$ and $\BBB$ be small dg categories and consider an object $B \in \BBB$. Then we have that the functor $\mathrm{L}(i_B)_{!} \cong - \dertp \BBB(-,B)$.
	\begin{proof}
		Since $\mathrm{L}(i_B)_{!}$ is a left adjoint between well generated triangulated categories, it preserves coproducts. Therefore, it is fully determined by its value on the representables, as they generate $\D(\AAA)$. Consider a module $F \in \D(\AAA \dertp \BBB)$. Then, for any object $A \in \AAA$ we have that
		\begin{equation*}
			\begin{aligned}
				\D(\AAA \dertp \BBB)\left( \mathrm{L}(i_B)_{!}(\AAA(-,A)), F\right)  &\cong \D(\AAA) \left( \AAA(-,A), (i_B)^*(F)\right) \\
				&\cong \D(\AAA)(\AAA(-,A), F(-,B))\\
				&\cong H^0(F(A,B)) \\
				&\cong \D(\AAA \dertp \BBB )(\AAA \dertp \BBB ((-,-),(A,B)), F)\\
				&\cong \D(\AAA \dertp \BBB) (\AAA(-,A) \dertp \BBB(-,B), F),
			\end{aligned}
		\end{equation*}
		where the first equivalence is given by the adjunction $\mathrm{L}(i_B)_{!} \dashv (i_B)^*$, the second by definition of $(i_B)^*$, the third and the fourth by definition of the morphisms in derived categories (see \cite[\S 4]{deriving-dg-categories}) and the last one can be readily seen using \Cref{cofibrantreplacement}. As this holds for all $F \in \D(\AAA \dertp \BBB)$, we conclude by Yoneda lemma.
	\end{proof}
\end{lemma}

\begin{lemma}\label{generatedbyaset2variables}
	Let $\AAA, \BBB$ be two small dg categories and $\www \subseteq \D(\AAA)$ a localising subcategory generated by a set. Then the triangulated subcategory
	$$\www' =  \{X \in \D(\AAA \dertp \BBB)\,\,|\,\, (i_B)^*(X) = X(-,B) \in \www \text{ for all } B \in \BBB \} $$
	of $\D(\AAA \dertp \BBB)$ is localising and generated by a set. In particular, if $\www$ is generated by a set $\nnn$, then we have that $\www'$ is generated by the set $\nnn' = \{ \mathrm{L}(i_B)_!(N) \,\,|\,\, N \in \nnn, B\in \BBB \}$. 
	\begin{proof}
		The fact that $(i_B)^*$ preserves small coproducts immediately shows that $\www'$ is a localising subcategory of $\D(\AAA \dertp \BBB)$. 
		
		We first prove that $\www'$ is generated by a set (and hence well generated by \cite[Thm 7.2.1]{localization-theory-triangulated-categories}). Given an object $B \in \BBB$, consider the composition
		$$\D(\AAA \dertp \BBB) \xrightarrow{(i_B)^*} \D(\AAA) \xrightarrow{\,\,\,Q\,\,\,} \D(\AAA)/\www,$$
		where $Q$ denotes the Verdier quotient functor. Observe that both $(i_B)^*$ and $Q$ preserve small coproducts, as they are left adjoint functors between well generated categories. Therefore, by \cite[Thm 7.4.1]{localization-theory-triangulated-categories} we have that
		$$\Kern(Q \circ (i_B)^*) = \{X \in \D(\AAA \dertp \BBB) \,\,|\,\,(i_B)^*X = X(-,B) \in \www \}$$
		is also well-generated and in particular generated by a set of objects. Notice now that, as 
		$$\www' = \bigcap_{B\in \BBB} \Kern(Q \circ (i_B)^*) \subseteq \D(\AAA \dertp \BBB),$$ 
		we can apply \cite[Lem 3.2]{bousfield-lattice-triangulated-category-stratification} to conclude that $\www'$ is also generated by a set of objects.
		
		We now show the second part of the statement, namely, that if $\www = \langle \nnn \rangle$ for a set $\nnn$, then $\www' = \langle \nnn' \rangle$ with $\nnn' = \{ \mathrm{L}(i_B)_!(N) \,\,|\,\, N \in \nnn, B\in \BBB \}$.
			
		We first prove that $\langle \nnn' \rangle \subseteq \www'$. As $\www'$ is localising, it suffices to show that $\nnn' \subseteq \www'$.
		Let's take $X = \mathrm{L}(i_B)_!(N) \in \nnn'$. We have that $X(-,B') = \mathrm{L}(i_B)_!(N)(-,B') = N(-) \dertp \BBB(B',B)$ by \Cref{adjunctiontwovariables}, and one can easily see that it belongs to $\www$. Indeed, we have that $N = N(-) \dertp k[0] \in \www$ where $k[0] \in \D(k)$ denotes the complex concentrated in degree 0 with $k$ in the 0-term. In addition, $k[0]$ is a compact generator of $\D(k)$, hence $\BBB(B',B) \in\D(k)$ can be written in terms of direct sums, extensions and shifts of $k[0]$. As $N(-) \dertp -: \D(k) \lra \D(\AAA)$ commutes with all these, and $\www$ is localising, we can conclude. Hence, we have that $\langle \nnn' \rangle \subseteq \www'$.	
		
		Now we prove that $\www' \subseteq \langle \nnn' \rangle$. Observe that it suffices to show that $\langle \nnn' \rangle^{\perp} \subseteq \www'^{\perp}$. Indeed, if we take left orthogonals, we have that $$\www' \subseteq \,^{\perp}(\www'^{\perp}) \subseteq \,^{\perp}(\langle \nnn' \rangle^{\perp}) = \langle \nnn' \rangle,$$
		where the last equality comes from \cite[Prop 4.9.1(6)]{localization-theory-triangulated-categories} because $\langle \nnn' \rangle$ is a localising subcategory generated by a set of a well generated category.	
		Recall from \Cref{leftorth} that $\langle \nnn' \rangle^{\perp} = \nnn'^{\perp}$. Let's consider $X \in \nnn'^{\perp}$. Then, we have that
		$$0 =\D(\AAA \dertp \BBB)(\mathrm{L}(i_B)_!(N), X),$$ for all for all $N \in \nnn$ and all $B\in \BBB$.
		Hence we have that 
		$$0 = \D(\AAA \dertp \BBB)(\mathrm{L}(i_B)_!(N), X) \cong \D(\AAA)(N, (i_B)^*(X))$$ 
		for all $N \in \nnn$ and all $B\in \BBB$. Thus $(i_B)^*(X) = X(-,B) \in \nnn^{\perp} = \www^{\perp}$ for all $B \in \BBB$. We are going to show that this is enough to conclude that $X \in \www'^\perp$.
		
		Observe that, because $\www'$ is a well generated subcategory of $\D(\AAA \dertp \BBB)$ closed under coproducts, by \cite[Thm 5.1.1]{localization-theory-triangulated-categories} we have that $\www'$ is a right admissible subcategory \cite[Def 1.2]{representable-functors-serre-functors-mutations} and $\langle \www'^\perp, \www' \rangle$ is a semiorthogonal decomposition of $\D(\AAA \dertp \BBB)$ (see for example \cite[Lem 2.5]{base-change-semiorthogonal-decompositions}). Therefore, we have a diagram of adjoint functors
		\begin{equation*}
			\begin{tikzcd}
				\www' \arrow[rr, "j'", hook, shift left] &  & \D(\AAA \dertp \BBB) \arrow[rr, "q'", shift left] \arrow[ll, "a'", shift left] &  & \www'^\perp \arrow[ll, "i'", hook, shift left]
			\end{tikzcd}
		\end{equation*}
		where $j' \dashv a'$, $q'\dashv i'$ with $a'j' \cong 1_{\www'}$, $q'i'\cong 1_{\www'^\perp}$ and furthermore $\ker(q') = \Img(j')$ and $\ker(a') = \Img(i')$. In particular, the projection functors associated to the semiorthogonal decomposition as in \cite[\S2.2]{base-change-semiorthogonal-decompositions} are precisely given by $i'q':\D(\AAA \dertp \BBB) \ra \D(\AAA \dertp \BBB)$ and $j'a':\D(\AAA \dertp \BBB) \ra \D(\AAA \dertp \BBB)$ (see, for example, the proof of \cite[Lem 3.1]{representations-associative-algebras-coherent-sheaves}). Analogously, we have that $\langle \www^\perp, \www\rangle$ is a semiorthogonal decomposition of $\D(\AAA)$ and thus we have a diagram of adjoint functors
		\begin{equation*}
			\begin{tikzcd}
				\www \arrow[rr, "j", hook, shift left] &  & \D(\AAA) \arrow[rr, "q", shift left] \arrow[ll, "a", shift left] &  & \www^\perp \arrow[ll, "i", hook, shift left]
			\end{tikzcd}
		\end{equation*}
		where $j \dashv a$, $q\dashv i$ with $aj \cong 1_{\www}$, $qi\cong 1_{\www^\perp}$ and furthermore $\ker(q) = \Img(j)$ and $\ker(a) = \Img(i)$. The projection functors associated to this semiorthogonal decomposition are $iq:\D(\AAA) \ra \D(\AAA)$ and $ja:\D(\AAA) \ra \D(\AAA)$. Now, observe that for all $B \in \BBB$ we have that the functor $(i_B)^*$ is compatible with the given semiorthogonal decompositions in the sense of \cite[\S3]{base-change-semiorthogonal-decompositions}, that is, we have that for all $B \in \BBB$:
		\begin{itemize}
			\item $(i_B)^*(\www') \subseteq \www$: This follows by definition of $\www'$;
			\item $(i_B)^*(\www'^\perp) \subseteq \www^\perp$: Let $X \in \www'^\perp$. Then, given any $Y \in \www$ we have that $\D(\AAA)(Y,(i_B)^*(X)) \cong \D(\AAA \dertp \BBB)(Y \dertp \BBB(-,B), X)$ and this latter is equal to $0$ because $Y \dertp \BBB(-,B) \in \www'$, which can be shown using the same argument as in the proof of the inclusion $\langle \nnn'\rangle \subseteq \www'$ above.
		\end{itemize}
		Hence, by applying \cite[Lem 3.1]{base-change-semiorthogonal-decompositions} we have that, with the notations above 
		\begin{equation}\label{compatibilityso}
			\begin{aligned}
				(i_B)^* j'a' &\cong j a (i_B)^*\\
				(i_B)^* i'q' &\cong i q (i_B)^*.
			\end{aligned}
		\end{equation}
		Notice now that, because $\langle \www'^\perp, \www' \rangle$ is a semiorthogonal decomposition, we have that our initial object $X \in \nnn'^\perp \subseteq \D(\AAA \dertp \BBB)$ fits in a distinguished triangle of the form
		\begin{equation*}
			j'a'(X) \lra X \lra i'q'(X) \lra j'a'(X)[1],
		\end{equation*}
		where the morphisms are induced by the counit and unit of the adjunctions above. If we now apply $(i_B)^*$, we obtain a distinguished triangle
		\begin{equation*}
			(i_B)^*j'a'(X) \lra (i_B)^*X \lra (i_B)^*i'q'(X) \lra (i_B)^*j'a'(X)[1].
		\end{equation*}
		Observe that $(i_B)^*j'a'(X) \cong ja(i_B)^*(X)$ for all $B \in \BBB$ because of \eqref{compatibilityso}. We showed above that for all $B \in \BBB$ we have that $(i_B)^*(X) = X(-,B)$ belongs to $\www^\perp$ and hence $a(i_B)^*(X) = 0$. Therefore, we have that $0 = ja(i_B)^*(X) \cong (i_B)^*j'a'(X)$ for all $B \in \BBB$. Consequently, we have that $j'a'(X) = 0$ and thus $X \cong i'q'(X) \in \www'^\perp$, which concludes the argument.
		\end{proof}
	\end{lemma}

\begin{remark}
	Observe that the presented proof is symmetric in the arguments, and hence the similar statement for $\www \subseteq \D(\BBB)$ a localising subcategory generated by a set holds as well.\end{remark}

We are now in disposition to prove:
\begin{theorem}\label{innerhomhalfway}
	Let $\BBB$ be a $\fraku$-small dg category and $\ccc$ a well generated $\frakv$-small dg category. Then, $\RHom(\BBB,\ccc)$ is a well generated $\frakv$-small dg category.
\end{theorem}
	\begin{proof}
		As $\ccc$ is pretriangulated by hypothesis, so is $\RHom(\BBB,\ccc)$ for any small dg category $\BBB$ (see for instance \cite[Rem E.2 \& E.4]{dg-quotients-dg-categories}). 
		
		As $\ccc$ is a well generated dg category, by Porta's Gabriel-Popescu theorem, there exists a small dg category $\CCC$ such that $\ccc$ is a Bousfield localisation of $\DD(\CCC)$, that is, there exists a quasi-fully faithful functor $i: \ccc \lra \DD(\CCC)$ which has a cocontinuous quasi-left adjoint $F \in \RHom_{\ccm}(\DD(\CCC),\ccc)$. 
		
		We have that, by \cite[Cor 6.6]{homotopy-theory-dg-categories-derived-morita-theory}, the morphism
		$$j: \RHom(\BBB,\ccc) \lra \RHom(\BBB,\DD(\CCC))$$
		induced by the dg functor $i$ is quasi-fully faithful as well, thus $H^0(j)$ is a fully-faithful functor. Observe that, if we consider $i$ as a right quasi-representable bimodule and we denote it by $h_{i(-)} \in \RHom(\ccc,\DD(\CCC))$, we have that $j = h_{i(-)} \dertp_{\ccc} (-)$. 
		Then, the natural bimodule $F' = F \dertp_{\DD(\CCC)} (-)  \in \RHom(\RHom(\BBB,\DD(\CCC)), \RHom(\BBB,\ccc))$ can be easily seen to be a quasi-left adjoint of $j:\RHom(\BBB,\ccc) \lra \RHom(\BBB,\DD(\CCC))$. Indeed, we have a counit $$F' \dertp_{\RHom(\bbb,\DD(\CCC))} j = F \dertp_{\DD(\CCC)} h_{i(-)} \dertp_{\ccc} (-) \lra \id_{\ccc} \dertp_{\ccc} (-) = \id_{\RHom(\BBB,\ccc)}$$
		induced by the counit of the quasi-adjunction $F \dashv_{H^0} h_{i(-)}$ and a unit
		$$\id_{\RHom(\BBB,\DD(\CCC))} = \id_{\DD(\CCC)} \dertp_{\DD(\CCC)} (-) \lra  h_{i(-)} \dertp_{\ccc} F \dertp_{\DD(\CCC)} (-) = j \dertp_{\RHom(\BBB,\ccc)} F'$$
		induced by the unit of the quasi-adjunction $F \dashv_{H^0} h_{i(-)}$.
		We thus have that $j$ is a dg Bousfield localisation of $\RHom(\BBB,\DD(\CCC))$, and hence   
		$$H^0(\RHom(\BBB,\ccc))  \xrightarrow{H^0(j)} H^0(\RHom(\BBB,\DD(\CCC)))$$
		is a Bousfield localisation of $H^0(\RHom(\BBB,\DD(\CCC)))$. 
		
		Now observe that $\RHom(\BBB,\DD(\CCC)) \cong \DD(\BBB^{\op} \dertp \CCC)$ in $\hodgcat$ as a direct consequence of the fact that $\DD(\CCC) \cong \RHom(\CCC^{\op}, \DD(k))$ in $\hodgcat$ (see \cite[\S 7]{homotopy-theory-dg-categories-derived-morita-theory}). We hence have an exact isomorphism $f:\D(\BBB^{\op} \dertp \CCC) \lra H^0(\RHom(\BBB,\DD(\CCC)))$ and it is not hard to see that every $X \in \D(\BBB^{\op} \dertp \CCC)$ is sent via $f$ to the associated quasi-respresentable bimodule $X: \BBB \lra \DD(\CCC)$ in $H^0(\RHom(\BBB,\DD(\CCC)))$. Consequently, via $f$, we have that 
		$$\D(\BBB^{\op} \dertp \CCC) \xrightarrow{H^0(F')\circ f} H^0(\RHom(\BBB,\ccc))  \xrightarrow{f^{-1} \circ H^0(j)} \D(\BBB^{\op} \dertp \CCC)$$ 
		provides a Bousfield localisation of $\D(\BBB^{\op} \dertp \CCC)$. In addition, observe that $$\Kern(H^0(F') \circ f) \cong \{X \in \D(\BBB^{\op} \dertp \CCC)\,\,|\,\, (i_B)^*(X) = X(B,-) \in \Kern(H^0(F)) \textrm{ for all }B \in \bbb \},$$ where $\Kern(H^0(F))$ is a localising subcategory of $\D(\CCC)$ generated by a set of objects. Then we can conclude by \Cref{generatedbyaset2variables} that $\Kern(H^0(F') \circ f)$ is also generated by a set of objects. Consequently, $H^0(\RHom(\BBB,\ccc))$ is also well generated, as we wanted to show. 
		\end{proof}

\begin{lemma}\label{lemimage}
Consider small dg categories $\BBB$, $\CCC$, $\CCC'$ and a quasi-fully faithful dg functor $\varphi: \CCC \lra \CCC'$. Consider the induced morphism 
$h_\varphi \dertp_{\CCC} - : \RHom(\BBB, \CCC) \lra \RHom(\BBB, \CCC')$. Consider a quasi-functor $F \in \RHom(\BBB, \CCC')$ which is such that $H^0(F): H^0(\BBB) \lra H^{0}(\CCC')$ factors through $H^0(\varphi)$.
Then, there exists $\bar{F} \in \RHom(\BBB,\CCC)$ such that $h_\varphi \dertp_{\CCC} \bar{F} \cong F$ as elements in $H^0(\RHom(\BBB,\CCC'))$.
\end{lemma}

\begin{proof}
Consider the first argument (derived) bimodule restriction functor $\varphi_1^{\ast}: \D(\BBB^{\op} \dertp \CCC') \lra \D(\BBB^{\op} \dertp \CCC)$ and the first argument derived bimodule extension functor $L((\varphi_1)_{!}): \D(\BBB^{\op} \dertp \CCC) \lra \D(\BBB^{\op} \dertp \CCC')$. Let $F$ be as in the statement of the lemma and put $\bar{F} = \varphi_1^{\ast}(F)$. 
Consider the Yoneda embeddings $y_{\CCC}: H^0(\CCC) \lra \D(\CCC)$, $y_{\CCC'}: H^0(\CCC') \lra \D(\CCC')$, the (derived) restriction functor $\varphi^{\ast}: \D(\CCC') \lra \D(\CCC)$ and the derived extension functor $L(\varphi_{!}): \D(\CCC) \lra \D(\CCC')$, for which we have 
\begin{equation}\label{lemeq1}
L(\varphi_{!}) y_{\CCC} \cong  y_{\CCC'} H^0(\varphi)
\end{equation}
in $\D(\CCC^{\op} \dertp \CCC')$. 
Consider $F \in \RHom(\BBB, \CCC')$ as in the statement of the lemma. We consider $F: \BBB \lra \DD(\CCC')$ and $H^0(F): H^0(\BBB) \lra \D(\CCC')$. Since $F$ is a quasi-functor, there exists $f: H^0(\BBB) \lra H^0(\CCC')$ with $H^0(F) \cong y_{\CCC'}f$ (note that we usually denote $H^0(F) = f$ for quasi-functors, but we refrain from doing so within this proof). 
By assumption, there exists $g: H^0(\BBB) \lra H^0(\CCC)$ with 
\begin{equation}\label{lemeq2}
H^0(\varphi) g \cong f.
\end{equation}
Using \eqref{lemeq1} and \eqref{lemeq2}, we thus have
\begin{equation}\label{lemeq3}
H^0(F) = y_{\CCC'}f \cong y_{\CCC'}H^0(\varphi) g \cong L(\varphi_{!}) y_{\CCC} g
\end{equation}
and 
\begin{equation}\label{lemeq4}
H^0(\bar{F}) = \varphi^{\ast} L(\varphi_{!}) y_{\CCC} g \cong y_{\CCC} g
\end{equation}
where in the last equation we have used that $\varphi$ is quasi-fully faithful.
Equation \eqref{lemeq4} already shows $\bar{F}$ to be a quasi-functor.
Comparing the expressions \eqref{lemeq3} and \eqref{lemeq4}, we see that $H^0(F) \cong L(\varphi_{!}) H^0(\bar{F})$ canonically.
From this, one readily deduced that the canonical natural transformation $h_\varphi \dertp_{\CCC} \bar{F} = L((\varphi_1)_{!}) \varphi_1^{\ast}(F) \lra F$ is an isomorphism in $H^0(\RHom(\BBB,\CCC'))$, as desired.
\end{proof}

\begin{corollary}\label{cocompleteRHom}\label{main0}
Consider a small dg category $\BBB$ and a dg category $\ccc$. Let $\alpha$ be a cardinal with $\alpha \leq |\fraku |$.
\begin{enumerate}
\item If $\ccc$ is homotopically $\alpha$-cocomplete, then so is $\RHom(\BBB, \ccc)$. Moreover, on the level of induced functors between the $H^0$-categories, coproducts are pointwise: for an $\alpha$-small family $(F_i)_{i \in I}$ with $F_i \in H^0(\RHom(\BBB, \ccc))$ with coproduct $\coprod_i F_i$, the functors $H^0(F_i), H^0(\coprod_i F_i): H^0(\BBB) \lra H^0(\ccc)$ are such that 
$$H^0(\coprod_i F_i)(B) = \coprod_i^{H^0(\ccc)} H^0(F_i)(B).$$
\item If $\ccc$ is homotopically cocomplete, then so is $\RHom(\BBB, \ccc)$, and the coproducts are pointwise on the level of induced functors between the $H^0$-categories.
\end{enumerate}
\end{corollary}		
		
\begin{proof}
Clearly, (2) is the case $\alpha = |\fraku |$ in (1). Suppose $\ccc$ is homotopically $\alpha$-cocomplete for $\alpha \leq |\fraku |$. Let $\frakv$ be a universe such that $\fraku \in \frakv$ and $\ccc$ is $\frakv$-small. Then $\frakv {-} \RHom(\BBB, \ccc)$ is constructed as an essentially $\frakv$-small dg category which is a $\fraku$-category, and $\alpha$ is a $\frakv$-small cardinal. For $Y_{\ccc}^{\alpha}: \ccc \lra \frakv{-}\DD_{\alpha}(\ccc)$ we have a canonical morphism
$$\tilde{Y}_{\ccc}^{\alpha}: \RHom(\BBB, \ccc) \lra \RHom(\BBB, \frakv{-}\DD_{\alpha}(\ccc))$$
which is quasi-fully faithful \cite[Cor 6.6]{homotopy-theory-dg-categories-derived-morita-theory}. Since the codomain is $\frakv$-well generated by \Cref{innerhomhalfway} and hence $\frakv$-cocomplete, it suffices to show that $H^0(\RHom(\BBB, \ccc))$ is closed under $\alpha$-small coproducts in $H^0(\RHom(\BBB, \frakv{-}\DD_{\alpha}(\ccc)))$.

Let $(F_i \in \RHom(\BBB, \ccc))_{i \in I}$ be an $\alpha$-small collection of objects. We may assume that $\BBB$ is cofibrant and that we have dg functors $f_i: \BBB \lra \ccc$ with $F_i = Y_{\ccc} f_i$ for $Y_{\ccc}: \ccc \lra {\frakv}{-}{\DD(\ccc)}$.

We will consider the functors $F_i^{\alpha} = Y^{\alpha}_{\ccc} f_i$ as representatives of the objects $\tilde{Y}_{\ccc}^{\alpha}(F_i) \in \RHom(\BBB, \frakv{-}\DD_{\alpha}(\ccc))$ (where we refrain from writing the composition with a further Yoneda embedding in order to obtain the associated bimodules). 

Consider the canonical quotient of dg $\frakv$-categories $Q: \frakv {-} \DD(\ccc) \lra \frakv {-} \DD_{\alpha}(\ccc)$ and the induced quotient
$$\tilde{Q}: \RHom(\BBB, \frakv {-} \DD(\ccc)) \lra \RHom(\BBB, \frakv {-} \DD_{\alpha}(\ccc)).$$
The coproduct of the objects $F_i^{\alpha} \in H^0(\RHom(\BBB, \frakv{-}\DD_{\alpha}(\ccc)))$ is given by $F = \tilde{Q}(\coprod_i F_i) = Q \circ \coprod_i F_i$ for $F_i \in H^0(\RHom(\BBB, \frakv {-} \DD(\ccc))) \cong \frakv{-}\D(\BBB^{\op} \otimes \ccc)$. By Lemma \ref{lemimage}, it suffices to show that $H^0(F)$ factors through $H^0(Y^{\alpha}_{\ccc})$. 
To see this, we compute
\begin{equation}\label{compcoprod}
H^0(F)(B) = Q(\coprod^{\frakv{-}\D(\ccc)}_i h_{f_i(B)}) = \coprod^{\frakv{-}\D_{\alpha}(\ccc)}_i h_{f_i(B)} \cong h_{ \coprod_{H^0(\ccc)} f_i(B)}
\end{equation}
where we have used the characterisation of the $\alpha$-cocontinuous derived category from Remark \ref{remDalpha}. The computation \eqref{compcoprod} also demonstrates the additional claim.
\end{proof}

\begin{corollary}\label{main1} \label{cocompletecontinuousRHom}
Consider dg categories $\aaa$ and $\bbb$. Let $\alpha$ be a cardinal with $\alpha \leq |\fraku |$.
\begin{enumerate}
\item If $\aaa$ and $\bbb$ are homotopically $\alpha$-cocomplete, then so is $\RHom_{\alpha}(\aaa, \bbb)$. Moreover, on the level of induced functors between the $H^0$-categories, the $\alpha$-small coproducts are pointwise. If in addition $\bbb$ is pretriangulated, then so is $\RHom_{\alpha}(\aaa, \bbb)$.
\item If $\aaa$ and $\bbb$ are homotopically cocomplete, then so is $\RHom_{\ccm}(\aaa, \bbb)$. Moreover, on the level of induced functors between the $H^0$-categories, small coproducts are pointwise. If in addition $\bbb$ is pretriangulated, then so is $\RHom_{\ccm}(\aaa, \bbb)$.
\end{enumerate}
\end{corollary}

	\begin{proof}
	Again, (2) is the case $\alpha = |\fraku |$ in (1). Let $\frakv$ be a universe with $\fraku \in \frakv$ and such that $\aaa$ and $\bbb$ are $\frakv$-small. Then $\alpha \leq |\frakv|$.
	From \Cref{cocompleteRHom} we know that $\RHom(\aaa,\bbb)$ is homotopically $\alpha$-cocomplete. In order to prove that $H^0(\RHom_{\alpha}(\aaa,\bbb))$ is $\alpha$-cocomplete it is enough to show that it is closed under $\alpha$-small coproducts in $H^0(\RHom(\aaa,\bbb))$. Consider an $\alpha$-small family $(F_i)_{i\in I}$ in $H^0(\RHom_{\alpha}(\aaa,\bbb))$. By \Cref{cocompleteRHom}, we have that 
	$$H^0\left( \coprod_i^{H^0(\RHom(\aaa,\bbb))} F_i\right)  \left( \coprod_j^{H^0(\aaa)} A_j\right)  \cong \coprod_i^{H^0(\bbb)}H^0(F_i)\left( \coprod_j^{H^0(\aaa)} A_j\right) $$ for all $\alpha$-small families $(A_j)_{j\in J}$ of elements of $H^0(\aaa)$. From the fact that the $F_i$'s are $\alpha$-cocontinuous we have that 
	$$\coprod_i^{H^0(\bbb)}H^0(F_i)\left( \coprod_j^{H^0(\aaa)} A_j\right)  = \coprod_j^{H^0(\bbb)} \coprod_i^{H^0(\bbb)} H^0(F_i)(A_j)= \coprod_j^{H^0(\bbb)} \left( H^0\left(\coprod_i^{H^0(\RHom(\aaa,\bbb))} F_i\right)(A_j) \right),$$
	which proves that $\coprod_i^{H^0(\RHom(\aaa,\bbb))} F_i$ belongs to $H^0(\RHom_{\alpha}(\aaa,\bbb))$ as desired.
	
	Now assume that $\bbb$ is pretriangulated. By \cite[Rem E.2 \& E.4]{dg-quotients-dg-categories}, we know that $\RHom(\aaa,\bbb)$ is also pretriangulated, with the triangulated structure inherited from that of $\DD(\bbb \dertp \aaa^{\op})$. It is then enough to show that $H^0(\RHom_{\alpha}(\aaa,\bbb)) \subseteq H^0(\RHom(\aaa,\bbb))$ is a triangulated subcategory. Take $F \in \RHom_{\alpha}(\aaa, \bbb)$ and consider its shift $F[1]$ when seen in the triangulated category $H^0(\RHom(\aaa, \bbb))$. We prove that $F[1]$ is $\alpha$-cocontinuous. Indeed, for any small family $(A_i)_i$ of objects of $\aaa$, we have that 
	\begin{equation*}
	\begin{aligned}
	(H^0(F)[1]) \left( \coprod_i^{H^0(\aaa)} A_i\right) &= \left( H^0(F) \left( \coprod_i^{H^0(\aaa)} A_i\right)\right)  [1] \\
	&= \left( \coprod_i^{H^0(\bbb)} H^0(F) (A_i)\right) [1] \\
	&= \coprod_i^{H^0(\bbb)}( H^0(F)(A_i) [1])\\
	&= \coprod_i^{H^0(\bbb)} (H^0(F)[1]) (A_i),
	\end{aligned}
	\end{equation*}
	where in the first and last equalities we use the fact that triangulated structure in $H^0(\RHom(\aaa,\bbb))$ is inherited from the canonical one in $\D(\bbb \dertp \aaa^{\op})$, in the second equality we use that $F$ is $\alpha$-cocontinuous and in the third equality we use that shifts commute with coproducts. Now consider an exact triangle $$F \lra F' \lra F'' \lra F[1]$$ 
	in $H^0(\RHom(\aaa,\bbb))$, where $F,F' \in H^0(\RHom_{\alpha}(\aaa,\bbb))$. Given an $\alpha$-small family $(A_i)$ of elements in $H^0(\aaa)$, for all $i$ we have the exact triangle 
	\begin{equation}\label{familytriang}
	\begin{tikzcd}
	H^0(F)(A_i) \arrow[r] &H^0(F')(A_i) \arrow[r] &H^0(F'')(A_i) \arrow[r] &H^0(F) (A_i)[1]
	\end{tikzcd}	
	\end{equation}
	in $H^0(\bbb)$. Observe now that we have the following diagram with rows exact triangles: 
	\begin{equation*}
	\begin{tikzcd}[column sep=0.29cm]
	H^0(F) \left( \overset{H^0(\aaa)}{\underset{i}{\coprod}} A_i\right) \arrow[r] \arrow[d,equal] &H^0(F') \left( \overset{H^0(\aaa)}{\underset{i}{\coprod}} A_i\right) \arrow[r] \arrow[d,equal] &H^0(F'') \left( \overset{H^0(\aaa)}{\underset{i}{\coprod}} A_i\right) \arrow[r] &H^0(F) \left( \overset{H^0(\aaa)}{\underset{i}{\coprod}} A_i\right)[1] \\
	\overset{H^0(\bbb)}{\underset{i}{\coprod}} H^0(F)(A_i) \arrow[r] &\overset{H^0(\bbb)}{\underset{i}{\coprod}}H^0(F')(A_i) \arrow[r] &\overset{H^0(\bbb)}{\underset{i}{\coprod}}H^0(F'')(A_i) \arrow[r] &\overset{H^0(\bbb)}{\underset{i}{\coprod}}H^0(F) (A_i)[1] 
	\end{tikzcd}	
	\end{equation*}
	where the exact triangle below is the coproduct of the family of exact triangles from \eqref{familytriang} above, and the vertical equalities are given because both $H^0(F)$ and $H^0(F')$ are $\alpha$-cocontinuous by hypothesis. By the axioms of triangulated categories, we have that
	$$H^0(F'') \left( \coprod_i^{H^0(\aaa)} A_i\right) \cong \coprod_i^{H^0(\bbb)}H^0(F'')(A_i).$$
	Hence $F'' \in \RHom_{\alpha}(\aaa,\bbb)$, which concludes the argument.    
	\end{proof}

\subsection{The cocontinuous internal hom of well generated dg categories}\label{parwellint}

In this section we prove that for well generated dg categories $\aaa$ and $\bbb$, the dg category $\RHom_c(\aaa, \bbb)$ is again well generated.

\begin{remark}\label{rem:notationYoneda}
	Let $\BBB$ be a small dg category and consider the associated dg Yoneda embedding $Y_\BBB:\BBB \ra \DD(\BBB)$. From this point on, we will abuse notations and write $Y_\BBB = h_{Y_\BBB} \in \RHom(\BBB,\DD(\BBB))$. In the same lines, if $\BBB$ is homotopically $\alpha$-cocomplete and $Y'_\BBB: \BBB \ra \DD_{\alpha}(\BBB)$ is the corestriction of the Yoneda embedding, we will write $Y'_\BBB = h_{Y'_{\BBB}} \in \RHom(\BBB, \DD_{\alpha}(\BBB))$.
\end{remark}

The following result extends To\"en's derived Morita theory (the case $\ccc = \DD(\CCC)$) from \cite[Thm 7.2]{homotopy-theory-dg-categories-derived-morita-theory} (see also \cite[Corollary 4.2]{internal-homs-via-extensions-dg-functors}). 
	
\begin{proposition}\label{kanextension} 
	Let $\BBB$ be a small dg category and $\ccc$ a well generated dg category. We have that the dg functor
	\begin{equation}\label{eqc}
	- \dertp Y_\BBB: \RHom_{\ccm}(\DD(\BBB), \ccc) \ra \RHom(\BBB, \ccc).
	\end{equation}
	is a quasi-equivalence, where $Y_{\BBB} \in \RHom(\BBB, \DD(\BBB))$ is the dg Yoneda embedding (see \Cref{rem:notationYoneda}). Therefore, $\RHom_{\ccm}(\DD(\BBB), \ccc) \cong \RHom(\BBB, \ccc)$ in $\hodgcat$.
	\begin{proof}
		If $\ccc = \DD(\CCC)$, the theorem reduces to derived Morita theory. In order to provide the proof for $\ccc$ an arbitrary well generated dg category, we will build upon the proof of \cite[Corollary 4.2]{internal-homs-via-extensions-dg-functors}.		
		Consider $\ccc$ a well generated dg category. In particular, by \Cref{thmenhGP} there exists a small dg category $\CCC$, a quasi-fully faithful dg functor $i: \ccc \lra \DD(\CCC)$ and a bimodule $a \in \RHom_{\ccm}(\DD(\CCC),\ccc)$ such that $a \dashv_{H^0} i$ and hence in particular, $a \dertp_{\DD(\CCC)} i \cong \id_{\ccc} \in H^0(\RHom(\ccc,\ccc))$ by \Cref{unit}. This implies that $[a]_{\iso} \circ [i] = \id_{\ccc}$ in $\Iso(H^0(\RHom(\ccc,\ccc))) \cong [\ccc,\ccc]$ (see \Cref{isoclasses}), where $[a]_{\iso}$ denotes the isomorphism class of $[a] \in H^0(\RHom(\ccc,\ccc))$.
		
		First, we prove that, for every well generated dg category $\ccc$, the map
		\begin{equation}\label{eq:bijectioninhodgcat}
			[\DD(\BBB), \ccc]_{\ccm} \lra [\BBB,\ccc]: f \mapsto f \circ [Y_{\BBB}]
		\end{equation}
		is a bijection, where $[-,-]_{\ccm}$ indicates the subset of morphisms in $\hodgcat$ such that the induced morphism between the homotopy categories preserves coproducts.
		 
		We first prove surjectivity. Consider $g \in [\BBB,\ccc]$. Then, $[i] \circ g \in [\BBB, \DD(\CCC)]$ and by derived Morita theory, there exists $f \in [\DD(\BBB), \DD(\CCC)]_{\ccm}$ such that $f \circ [Y_{\BBB}] = [i] \circ g$. Consider now $[a]_{\iso} \circ f$, which belongs to $[\DD(\BBB),\ccc]_{\ccm}$. Then, $[a]_{\iso} \circ f \circ [Y_{\BBB}] = [a]_{\iso} \circ [i] \circ g = g$, which proves surjectivity. 
		
		In order to prove injectivity, one can follow a very similar argument to that of the proof of \cite[Prop 3.10]{internal-homs-via-extensions-dg-functors} in which a first step towards the proof of derived Morita theory is provided. We provide the details here for convenience of the reader. Consider $f_1,f_2 \in [\DD(\BBB), \ccc]_{\ccm}$ such that $f_1 \circ [Y_{\BBB}] = f_2 \circ [Y_{\BBB}]$. By composing with $[i]$, we have $[i] \circ f_1 \circ [Y_{\BBB}] = [i] \circ f_2 \circ [Y_{\BBB}] \in [\BBB,\DD(\CCC)]$. It follows from \cite[Prop 2.11(3)]{internal-homs-via-extensions-dg-functors} that there exists a dg category $\aaa$ and a quasi-equivalence $I: \aaa \lra \DD(\BBB)$ such that $f_i = [F_i] \circ [I]^{-1}$ with $F_i: \aaa \lra \ccc$ a dg functor for $i=1,2$. Consequently, we have that $[i] \circ f_i = [i \circ F_i] \circ [I]^{-1}$ for $i = 1,2$. We denote by $\AAA$ the full dg subcategory of $\aaa$ such that $I' \coloneqq I_{|\AAA}$ induces a quasi-equivalence of $\AAA$ with the full dg subcategory $\qrep(\BBB)$ of $\DD(\BBB)$ and by $J: \AAA \hookrightarrow \aaa$ the inclusion. We write $G_i\coloneqq F_i \circ J: \AAA \lra \ccc$ for $i =1,2$. We hence have, for $i =1,2$, the following commutative diagram
		\begin{equation}
			\begin{tikzcd}
			\DD(\BBB) & \aaa \arrow[l, "\sim","I"'] \arrow[r, "i \circ F_i"] & \DD(\CCC). \\
			\qrep(\BBB) \arrow[u, hook,"J'"] & \AAA \arrow[l, "\sim","I'"'] \arrow[ru, "i \circ G_i"'] \arrow[u, "J", hook] & 
			\end{tikzcd}
		\end{equation}
		Following the notations of \cite[\S 3.1]{internal-homs-via-extensions-dg-functors}, we consider the extension of $i \circ G_i$:
		$$\widehat{i \circ G_i}: \dgmod(\AAA) \ra \dgmod(\CCC) : X \mapsto E_i \dertp X$$ 
		where $E_i \in \dgmod(\CCC \dertp \AAA^{\op})$ is the bimodule corresponding to the functor $i \circ G_i: \AAA \lra \DD(\CCC)$, and the restriction of $i \circ G_i$:
		$$\widetilde{i \circ G_i}: \dgmod(\CCC) \ra \dgmod(\AAA): X \mapsto \dgmod(\CCC)(i \circ G_i(-), X).$$
		By \cite[Prop 3.2]{internal-homs-via-extensions-dg-functors}, we have that $H^0(\widehat{i \circ G_i})$ is cocontinuous and it is easy to check that $\widehat{i \circ G_i}$ restricts to a dg functor $\DD(\AAA) \lra \DD(\CCC)$, because $i \circ G_i(\AAA) \subseteq \DD(\CCC)$. In addition, we have the adjunction $\widehat{i \circ G_i} \dashv \widetilde{i \circ G_i}$. 
		
		From the discussion above, we have that $$[i \circ G_1] \circ [I']^{-1}= [i \circ F_1] \circ [I]^{-1}\circ[J'] = [i \circ F_2] \circ [I]^{-1}\circ[J'] =[i \circ G_2] \circ [I']^{-1},$$ and hence $[ i \circ G_1 ] = [i \circ G_2]$. From \cite[Lem 3.9]{internal-homs-via-extensions-dg-functors} if follows that $[\widehat{i \circ G_1}]=[\widehat{i \circ G_2}]$.
		
		Consider now the restriction functor $J^*:\dgmod(\aaa) \ra \dgmod(\AAA): X \mapsto X\circ J$ and the composition $K \coloneqq J^* \circ Y_{\aaa}$. From the proof of \cite[Prop 1.17]{uniqueness-enhancement-triangulated-categories} it follows that $H^0(K)$ is cocontinuous and one has that $K(J(\AAA)) \cong Y_{\AAA}(\AAA) \subseteq \DD(\AAA)$. 
		
		Now observe that $\widetilde{i \circ G_i} \circ i \circ F_i (A) = \DD(\CCC)(i \circ G_i(-), i \circ F_i(A))$ for all $A \in \aaa$, and hence we have the following natural transformation $\gamma: K \lra \widetilde{i \circ G_i} \circ i \circ F_i$ induced by $i \circ F_i$:
		$$\gamma_A: K(A) = \aaa(J(-),A) \lra \DD(\CCC)(i \circ F_i \circ J(-), i \circ F_i (A)) = \widetilde{i \circ G_i} \circ i \circ F_i (A).$$
		By adjunction, we have a natural transformation $\beta: \widehat{i \circ G_i} \circ K \rightarrow i \circ F_i$ with the property that $H^0(\beta)_{|J(\AAA)}$ is an isomorphism, where $H^0(J(\AAA))$ forms a compact generator of the well generated triangulated category $H^0(\aaa)$. Consider the functor $\Phi_a: \DD(\CCC) \lra \qrep(\ccc)$ associated to $a \in \RHom_{\ccm}(\DD(\CCC),\ccc)$. By composing with $\Phi_a$ we obtain a natural transformation $\alpha: \Phi_a \circ \widehat{i\circ G_i} \circ K \lra \Phi_a \circ i \circ F_i$ such that $H^0(\alpha)_{|J(\AAA)}$ is an isomorphism. Then, we have that $H^0(\ccc)$ is well generated, $H^0(\Phi_a)$, $H^0( \widehat{i\circ G_i})$ and $H^0(K)$ are cocontinuous, and so is $H^0(\Phi_a) \circ H^0(i) \circ H^0(F_i) \cong H^0(F_i)$. Consequently, by the same argument of \cite[Rem 2.4]{internal-homs-via-extensions-dg-functors}, we have that $\alpha$ is a termwise homotopy equivalence, and hence $[a]_{\iso} \circ [\widehat{i\circ G_i}] \circ [K] = [a]_{\iso} \circ [i] \circ [F_i] = [F_i]$ for $i=1,2$. Now, as $[\widehat{i\circ G_1}] = [\widehat{i\circ G_2}]$, we obtain that $[F_1] = [F_2]$. This finally implies that $f_1 = f_2$ as desired. 
		
		Now, define $[ \DD(\BBB) \dertp \AAA, \ccc]'_{\ccm}$ as the subset of $[ \DD(\BBB) \dertp \AAA, \ccc]$ consisting of elements $f$ such that $H^0(f)(-,A)$ preserves coproducts for all $A \in \AAA$. Then, we have the following commutative diagram induced by the Yoneda embedding $Y_{\BBB}:\BBB \lra \DD(\BBB)$:
		\begin{equation}
		\begin{tikzcd}[column sep= 50pt]
		\left[ \AAA, \RHom_{\ccm}(\DD(\BBB),\ccc)\right]  \arrow[r,"{[-\dertp Y_\BBB]} \circ -"] &\left[ \AAA, \RHom(\BBB,\ccc)\right] \\	
		\left[ \DD(\BBB) \dertp \AAA, \ccc\right]'_{\ccm} \arrow[d,"\cong"'] \arrow[u,"\cong"]  &\left[ \AAA \dertp \BBB, \ccc\right] \arrow[d,"\cong"'] \arrow[u,"\cong"]\\
		\left[ \DD(\BBB), \RHom(\AAA,\ccc)\right]_{\ccm} \arrow[r,"- \circ {[Y_{\BBB}]}"]&\left[ \BBB, \RHom(\AAA,\ccc)\right], 
		\end{tikzcd}
		\end{equation}
		where the vertical arrows are induced by the $\dertp$-$\RHom$ adjunction in $\hodgcat$.
		As $\RHom(\AAA,\ccc)$ is well generated by \Cref{innerhomhalfway}, we have that the lower horizontal arrow is a bijection by the discussion above, and thus so is $[-\dertp Y_\BBB] \circ -$. Then, using \Cref{propunivyon} we can conclude that $[-\dertp Y_\BBB]$ is an isomorphism in $\hodgcat$, proving that $- \dertp Y_\BBB$ is a quasi-equivalence, as desired.
	\end{proof}
\end{proposition}
\begin{proposition}\label{derivedalphacocompletion} 
	Let $\BBB$ be a homotopically $\alpha$-cocomplete small dg category and $\ccc$ a well generated dg category. We have that the dg functor
	\begin{equation}\label{eqcbis}
	- \dertp Y'_\BBB: \RHom_{\ccm}(\DD_{\alpha}(\BBB), \ccc) \ra \RHom_{\alpha}(\BBB, \ccc)
	\end{equation}
	is a quasi-equivalence, where $Y'_\BBB \in \RHom(\BBB,\DD_\alpha(\BBB))$ is the corestriction of the dg Yoneda embedding (see \Cref{rem:notationYoneda}). Therefore, $\RHom_{\ccm}(\DD_{\alpha}(\BBB), \ccc) \cong \RHom_{\alpha}(\BBB, \ccc)$ in $\hodgcat$.
	\begin{proof}
		First recall that we have a dg Bousfield localisation $a \dashv_{H^0} i: \DD_{\alpha}(\BBB) \rightleftarrows \DD(\BBB)$. One can easily see that $a \dertp_{\DD(\BBB)} Y_{\BBB} \in \RHom_{\alpha}(\BBB,\DD_{\alpha}(\BBB))$ is $\alpha$-cocontinuous and it is isomorphic in $H^0(\RHom_{\alpha}(\BBB,\DD_{\alpha}(\BBB)))$ to the corestriction of the Yoneda embedding $Y_{\BBB}': \BBB \lra \DD_{\alpha}(\BBB)$, which can be easily deduced from the fact that $Y_{\BBB} = i \circ Y'_{\BBB}$. Hence, we have that $$i \dertp_{\DD_{\alpha}(\BBB)} a \dertp_{\DD(\BBB)} Y_{\BBB} \cong Y_{\BBB} \in H^0(\RHom(\BBB,\DD(\BBB))).$$ 
		
		We are going to show that, for any well generated dg category $\ccc$, the map
		\begin{equation}\label{eq4}
				[\DD_{\alpha}(\BBB),\ccc]_{\ccm} \lra [\BBB,\ccc]_{\alpha}: f \longmapsto f \circ [a]_{\iso} \circ [Y_{\BBB}]
		\end{equation}
		is a bijection, where $[-,-]_{\ccm}$ (resp. $[-,-]_{\alpha}$) indicates the subset of morphisms in $\hodgcat$ such that the induced morphism between the homotopy categories preserves small coproducts (resp. $\alpha$-small coproducts). Given $f \circ [a]_{\iso} \circ [Y_{\BBB}] = f'\circ [a]_{\iso}  \circ [Y_{\BBB}]$, then, by derived Morita theory, as both $f \circ [a]_{\iso}$ and $f' \circ [a]_{\iso}$ are cocontinuous,  we have that $f \circ [a]_{\iso} = f' \circ [a]_{\iso}$. Consequently,  $$f=f \circ [a]_{\iso} \circ [i]_{\iso} = f' \circ [a]_{\iso} \circ [i]_{\iso} = f',$$ which proves injectivity. 
		
		Next, consider $g \in [\BBB,\ccc]_{\alpha}$. Then, by derived Morita theory, there is an element $f \in [\DD(\BBB),\ccc]_{\ccm}$ such that $f \circ [Y_{\BBB}] = g$. We are going to show that $f$ factors through $[a]_{\iso} \in [\DD(\BBB), \DD_{\alpha}(\BBB)]_{\ccm}$. Indeed, by the description of the kernel of $H^0(a)$ provided in \Cref{remDalpha} and the universal property of the dg quotient \eqref{dgquotientrestricted}, $f$ factors through $[a]_{\iso}$ if and only if $$H^0(f)(\coprod_i Y_{\BBB}(B_i)) \cong H^0(f)(Y_{\BBB}(\coprod_i B_i)),$$ where $\coprod_i B_i$ is seen in $H^0(\BBB)$, for all $\alpha$-small coproducts. But this condition is readily seen to be satisfied taking into account that $f$ is cocontinuous and $f \circ [Y_{\BBB}] = g$ is $\alpha$-cocontinuous, and hence 
		\begin{equation*}
			f = t \circ [a]_{\iso}.
		\end{equation*}
		In addition, $t$ is also cocontinuous by the universal property of the dg quotient \eqref{dgquotientrestricted}, that is $t \in [\DD_{\alpha}(\BBB),\ccc]_{\ccm}$. 
		
		Now, observe that $$t \circ [a]_{\iso} \circ [Y_{\BBB}] = f \circ [Y_{\BBB}] = g,$$ which proves surjectivity. 
		
		Now, we define $[\DD_{\alpha}(\BBB) \dertp \AAA, \ccc]'_{\ccm}$ as the subset of $[\DD_{\alpha}(\BBB) \dertp \AAA, \ccc]$ consisting of elements $f$ such that $H^0(f)(-,A)$ preserves coproducts for all $A \in \AAA$. We define analogously $[\BBB \dertp \AAA, \ccc]'_{\alpha}$ as the subset of $\left[\BBB \dertp \AAA, \ccc\right]$ consisting of elements $f$ such that $H^0(f)(-,A)$ preserves $\alpha$-small coproducts for all $A \in \AAA$. Then, we have the following commutative diagram induced by the dg functor $Y'_\BBB: \BBB \lra \DD_{\alpha}(\BBB)$:
		\begin{equation}
		\begin{tikzcd}[column sep= 50pt]
		\left[ \AAA, \RHom_{\ccm}(\DD_{\alpha}(\BBB),\ccc)\right]  \arrow[r, " {[-\dertp Y'_{\BBB}]\circ-}"] &\left[ \AAA, \RHom_{\alpha}(\BBB,\ccc)\right] \\	
		\left[ \DD_{\alpha}(\BBB) \dertp \AAA, \ccc\right]'_{\ccm} \arrow[d,"\cong"'] \arrow[u,"\cong"]  &\left[ \BBB \dertp \AAA, \ccc\right]'_{\alpha} \arrow[d,"\cong"'] \arrow[u,"\cong"]\\
		\left[ \DD_{\alpha}(\BBB), \RHom(\AAA,\ccc)\right]_{\ccm} \arrow[r,"- \circ {[a]_{\iso}} \circ {[Y_{\BBB}]}","-\circ {[Y'_\BBB]}"']&\left[ \BBB, \RHom(\AAA,\ccc)\right]_{\alpha}, 
		\end{tikzcd}
		\end{equation}
		where the vertical arrows are induced by the $\dertp$-$\RHom$ adjunction in $\hodgcat$.
		As $\RHom(\AAA,\ccc)$ is well generated by \Cref{innerhomhalfway}, we have that the lower horizontal arrow is a bijection, and hence so is $[-\dertp Y'_{\BBB}]\circ-$. Therefore, using Proposition \ref{propunivyon} we can conclude that $[-\dertp Y'_{\BBB}]$ is an isomorphism in $\hodgcat$, as desired.		
	\end{proof}
\end{proposition}

Before proving the main result of the section, we will need the following lemma.

\begin{lemma}\label{intersectionwg}
	Let $\aaa$ be a well generated dg category and consider a small family of full well generated pretriangulated dg subcategories $\{\bbb_i\}_{i \in I}$ of $\aaa$ closed under homotopy coproducts. Then $\bigcap_i \bbb_i$ is a well generated pretriangulated dg subcategory of $\aaa$.
	\begin{proof}
		Observe that $H^0(\bigcap_i \bbb_i) = \bigcap_i H^0(\bbb_i)$ is a triangulated subcategory of $\aaa$, and hence $\bigcap_i \bbb_i$ is a pretriangulated dg subcategory of $\aaa$. It is thus sufficient to show that $H^0(\bigcap_i \bbb_i)$ is well generated. By hypothesis, we have that $H^0(\aaa)$ is well generated and that, for all $i \in I$,$H^0(\bbb_i) \subseteq H^0(\aaa)$ is a localising subcategory generated by a set of objects. Consequently, $\bigcap_{i \in I} H^0(\bbb_i) = H^0(\bigcap_{i \in I} \bbb_i)$ is also a localising subcategory of $H^0(\aaa)$ generated by a set of objects \cite[Lem 3.2]{bousfield-lattice-triangulated-category-stratification}. We can conclude by applying \cite[Thm 7.2.1]{localization-theory-triangulated-categories} that $H^0(\bigcap_i \bbb_i)$ is well generated.
	\end{proof}
\end{lemma}

\begin{remark}
	The proof of the following theorem is a dg parallel of the argument followed in \cite[Thm 2.60]{locally-presentable-accessible-categories} in order to prove that the category of models of a sketch taking values in an accessible category is again accessible.
\end{remark}	

\begin{theorem} \label{closed} \label{notnecwelgenerated}
	Let $\aaa,\bbb$ be two well generated dg categories. Then $\RHom_{\ccm}(\aaa,\bbb)$ is well generated.
	\begin{proof}
		By \Cref{thmenhGP}, we can choose a cardinal $\alpha$ such that $\aaa \cong \DD_{\alpha}(\AAA)$ for $\AAA$ a homotopically $\alpha$-cocomplete small dg category. We can further assume that $\AAA$ is cofibrant. By \Cref{derivedalphacocompletion}, it is enough to prove that $\RHom_{\alpha}(\AAA,\bbb)$ is well generated. 
		Consider the small family $\Lambda = \{(A_i)_{i\in I} \,\,|\,\, A_i \in \AAA, |I| < \alpha\}$ of all $\alpha$-small families of objects of $\AAA$. Given $\lambda=(A_i)_{i\in I} \in \Lambda$, denote by $\eee_{\lambda}$ the full dg subcategory of $\RHom(\AAA,\bbb)$ with objects $F$ such that the canonical morphism $\coprod_i^{H^0} H^0(F)(A_i) \ra H^0(F)(\coprod_i^{H^0} A_i)$ is an isomorphism in $H^0(\bbb)$. Observe that $\RHom_{\alpha}(\AAA,\bbb) = \bigcap_{\lambda \in \Lambda} \eee_{\lambda}$. We claim it is enough to prove that $\eee_{\lambda}$ is well generated for each $\lambda$. Indeed, we know by \Cref{innerhomhalfway} that $\RHom(\AAA,\bbb)$ is well generated, and one can readily check following the same argument of the proof of \Cref{cocompletecontinuousRHom} that $\eee_{\lambda}$ are pretriangulated dg subcategories of $\RHom(\AAA,\bbb)$ closed under homotopy coproducts. Hence, by \Cref{intersectionwg}, if  $\eee_{\lambda}$ is well generated for every $\lambda \in \Lambda$, we can conclude that $\bigcap_{\lambda \in \Lambda} \eee_{\lambda}$ is a well generated pretriangulated dg subcategory of $\RHom(\AAA,\bbb)$.
		
		It hence remains to prove that $\eee_{\lambda}$ is well generated for every $\lambda \in \Lambda$. In order to show this we will prove that $\eee_{\lambda}$ is a homotopy fiber product of a cospan diagram of well generated dg categories with cocontinuous dg functors. This allows us to conclude using the fact that the homotopy category of well generated dg categories (i.e. locally presentable dg categories) with cocontinuous morphisms is closed in $\hodgcat$ under homotopy limits (see the proof of \cite[Lem 3.3]{derived-azumaya-algebras-generators-twisted-derived-categories} or \cite[Rem 6.2.2]{motivic-homotopy-theory-noncommutative-spaces}).
		
		Fix $\lambda=(A_i)_{i\in I} \in \Lambda$ and consider the family of canonical morphisms $s_i: A_i \ra \coprod_i^{H^0}A_i$ in $H^0(\AAA)$. We fix a family $r_i: A_i \ra \coprod_i^{H^0}(A_i)$ in $Z^0(\AAA)$ lifting the $s_i$. Consider the dg category $\mathrm{Ar}_0$ with two objects $X, X'$ and morphisms $\mathrm{Ar}_0(X,X) = k1_X$, $\mathrm{Ar}_0(X',X') = k1_{X'}$, $\mathrm{Ar}_0(X',X) = 0$ and $\mathrm{Ar}_0(X,X') = kx$ with $x$ a morphism in degree 0. We introduce the following dg functors:
		\begin{itemize}
			\item We define the dg functor $C: \RHom_\ccm(\DD(\AAA),\bbb) \ra \RHom(\mathrm{Ar}_0,\bbb)$ as follows. For any $F \in \RHom_{\ccm}(\DD(\AAA),\bbb)$ we associate the quasi-functor $C(F)$, that as a dg functor $\Phi_{C(F)}:\mathrm{Ar}_0 \ra \qrep(\bbb)$ is given by assigning to $X$ the object $\Phi_F(\coprod_i h_{A_i}) \in \qrep(\bbb)$, to $X'$ the object $\Phi_F(h_{\coprod_i^{H^0} A_i}) \in \qrep(\bbb)$ and to $x$ the morphism $\Phi_F(\mathrm{can}_\lambda)$ where $\mathrm{can}_\lambda: \coprod_i h_{A_i} \ra h_{\coprod_i^{H^0}A_i}$ is the canonical morphism in $\DD(\AAA)$ induced by $h_{r_i}$. In order to lighten the notations, from this point on and for the rest of the proof we will not distinguish between right quasi-representable bimodules $F$ and their associated dg functor $\Phi_F$. Given a morphism $\gamma : F \ra G$, we associate to it the following natural morphism $C(\gamma)$ in $\RHom(\mathrm{Ar}_0,\bbb)$:
			\begin{equation*}\rule{\leftmargin}{0in}
				\begin{tikzcd}[column sep=16ex]
						F(\coprod_i h_{A_i}) \arrow[r,"C(F)(x) = F(\mathrm{can}_\lambda)"] \arrow[d,"\gamma_X = \gamma_{\coprod_i h_{A_i}}"'] &F(h_{\coprod_i^{H^0} A_i}) \arrow[d,"\gamma_{X'} = \gamma_{h_{\coprod_i^{H^0} A_i}} "]\\ 
						G(\coprod_i h_{A_i}) \arrow[r,"C(G)(x) = G(\mathrm{can}_\lambda)"] &G(h_{\coprod_i^{H^0} A_i}),
				\end{tikzcd}
			\end{equation*}
			already seen inside $\qrep(\bbb)$. We will denote this morphism by $(\gamma_{\coprod_i h_{A_i}},\gamma_{h_{\coprod_i^{H^0} A_i}})$, and from now on we will follow this notation for morphisms in $\RHom(\mathrm{Ar}_0, \bbb)$. In particular, if $\phi \in \RHom(\mathrm{Ar}_0, \bbb)$, we write $\phi = (\phi_1, \phi_2)$.
			\item We define the dg functor $I: \qrep(\bbb) \ra \RHom(\mathrm{Ar}_0,\bbb)$ given by associating to each $B \in \qrep(\bbb)$ the quasi-functor with constant value $B$ and such that $I(B)(x) = \id_{B}$. We define $I$ on morphisms in the natural way.
		\end{itemize}
		We are going to show that $\eee_{\lambda}$ is the homotopy limit of the following diagram
		\begin{equation*}
		\begin{tikzcd}
		&\qrep(\bbb) \arrow[d,"I"]\\
		\RHom_{\ccm}(\DD(\AAA),\bbb) \arrow[r,"C"] &\RHom(\mathrm{Ar}_0,\bbb).
		\end{tikzcd}
		\end{equation*}
		This will allow us to conclude. Indeed, $\bbb$ is a well generated dg category and hence so is $\qrep(\bbb)$ because they are isomorphic in $\hodgcat$. In addition, by \Cref{derivedalphacocompletion} we have that $\RHom_{\ccm}(\DD(\AAA),\bbb) \cong \RHom(\AAA,\bbb)$ in $\hodgcat$. Consequently, as a direct consequence of \Cref{innerhomhalfway}, we can conclude that $\RHom_{\ccm}(\DD(\AAA),\bbb)$ and $\RHom(\mathrm{Ar}_0,\bbb)$ are also well generated dg categories. Furthermore, both $I$ and $C$ are easily seen to be cocontinuous.
		
		In \cite[\S 4]{milnor-descent-cohesive-dg-categories} a model for the homotopy limit in $\hodgcat$ is described using path objects. In what follows we will construct a quasi-equivalence from this concrete model to $\eee_{\lambda}$. Let us begin with describing the model for $\ppp \coloneqq \RHom_{\ccm}(\DD(\AAA),\bbb) \times^{\text{h}}_{\RHom(\mathrm{Ar}_0,\bbb)} \qrep(\bbb)$. The objects of $\ppp$ are given by
		\begin{equation*}
		\begin{aligned}
	    \big\{ (F, B, \phi) \,\,|\,\, F \in& \RHom_{\ccm}(\DD(\AAA),\bbb), B \in \qrep(\bbb), \phi \in \RHom(\mathrm{Ar}_0,\bbb)^0(C(F),I(B)),\\
		&\phi \text{ is closed and becomes an isomorphism in } H^0(\RHom(\mathrm{Ar}_0,\bbb)) \big\}.
		\end{aligned}
		\end{equation*}
		The morphisms of degree $n$ are given by
		\begin{equation*}
		\begin{aligned}
		\ppp^n&((F_1,B_1,\phi_1),(F_2,B_2,\phi_2)) =\\ 
		&\RHom_{\ccm}(\DD(\AAA),\bbb)^n(F_1,F_2) \oplus \qrep(\bbb)^n (B_1,B_2) \oplus \RHom(\mathrm{Ar}_0,\bbb)^{n-1}(C(F_1),I(B_2)).
		\end{aligned}
		\end{equation*}
		Given a morphism $(\gamma,\mu,\nu): (F_1,B_1,\phi_1)\ra (F_2,B_2,\phi_2)$ of degree $n$ and a morphism $(\gamma',\mu',\nu'): (F_2,B_2,\phi_2)\ra (F_3,B_3,\phi_3)$ the composition is provided by
		$$(\gamma',\mu',\nu')(\gamma,\mu,\nu) = (\gamma'\gamma,\mu'\mu, (-1)^n \nu'C(\gamma) +  I(\mu') \nu) $$
		and the differential is given by
		$$d(\gamma,\mu,\nu) = (d\gamma,d\mu,d\nu + (-1)^n(\phi_2 C(\gamma) - I(\mu)\phi_1)).$$
		
		We first show that if $(F,B,\phi) \in \ppp$, then the quasi-representable bimodule $F \dertp Y \in \RHom(\AAA,\bbb)$ lies in $\eee_{\lambda}$. Indeed, because $\phi:C(F) \ra I(B)$ is a homotopy equivalence, we have that $H^0(F(\mathrm{can}_{\lambda}))$ is an isomorphism and, because $F$ is cocontinuous, the canonical morphism $\coprod_i^{H^0}(H^0(F)(h_{A_i})) \ra H^0(F)(\coprod_i h_{A_i})$ also is. Therefore, the composition
		\begin{equation*}
		\begin{tikzcd}
			\coprod_i^{H^0}(H^0(F)(h_{A_i})) \arrow[rr]  & & H^0(F)(\coprod_i h_{A_i}) \arrow[rr, "H^0(F(\mathrm{can}_\lambda))"] & & H^0(F)(h_{\coprod_i^{H^0} A_i})
		\end{tikzcd} 
		\end{equation*}	
		is an isomorphism, and thus we can conclude that the canonical morphism
		\begin{equation*}
			\coprod_i^{H^0}(H^0(F\dertp Y)(A_i)) \ra  H^0(F \dertp Y)(\coprod_i^{H^0} A_i)
		\end{equation*}
		is an isomorphism as well, proving the claim. 
		
		We define a dg functor $S: \ppp \ra \eee_{\lambda}$ as follows.
		To every $(F,B,\phi) \in \ppp$ we associate $F \dertp Y \in \eee_{\lambda}$ and to every morphism $(\gamma,\mu,\nu) \in  \ppp((F,B,\phi),(F',B',\phi'))$ we associate the morphism $\gamma \dertp Y \in \eee_{\lambda}(F\dertp Y,F'\dertp Y)$. It is readily seen that this is indeed a dg functor.
		To conclude, it is enough to show that $S$ is a quasi-equivalence.
		
		We first show that $S$ is quasi-essentially surjective. We know from the proof of \Cref{kanextension} that the dg functor
		\begin{equation*}
			\RHom_\ccm(\DD(\AAA),\bbb) \ra \RHom(\AAA,\bbb): F \mapsto F \dertp Y
		\end{equation*}
		is a quasi-equivalence. Consequently, given $G \in H^0(\eee_{\lambda})\subseteq H^0(\RHom(\AAA,\bbb))$, we can choose an $F \in H^0(\RHom_\ccm(\DD(\AAA),\bbb))$ such that there is an isomorphism $\psi: F \dertp Y \ra G$ in $H^0(\RHom(\AAA,\bbb))$. It is then easy to check that $F(\mathrm{can}_\lambda)$ induces an isomorphism in $H^0(\qrep(\bbb))$. 
		Denote by $\phi \in \RHom^0(\mathrm{Ar}_0,\bbb)(C(F),I(G(\coprod_i^{H^0}A_i)))$ the closed morphism of degree 0 given by $(\overline{\psi}\circ F(\mathrm{can}_\lambda), \overline{\psi})$, where $\overline{\psi}$ is a $0$-cycle lifting the isomorphism 
		\begin{equation*}
			\psi_{|\coprod_i^{H^0}A_i}: F\dertp Y(\coprod_i^{H^0}A_i) \ra G(\coprod_i^{H^0}A_i)
		\end{equation*}
		in $H^0(\qrep(\bbb))$. Observe that $\phi$ becomes an isomorphism in $H^0(\RHom(\mathrm{Ar}_0,\bbb))$.
		Therefore, we have that $(F, G(\coprod_i^{H^0}A_i),\phi)$ belongs to $\ppp$, and it is easy to see that, seen as an object in $H^0(\ppp)$, it is sent to $G \in H^0(\eee_{\lambda})$, proving that $S$ is quasi-essentially surjective as desired.
		
		We now show that $S$ is quasi-full. Consider $\sigma \in \zzz^n(\eee_{\lambda})(S(F,B,\phi),S(F',B',\phi'))=\zzz^n(\eee_{\lambda})(F\dertp Y,F'\dertp Y)$. As $\eee_{\lambda}$ is a full dg subcategory of $\RHom(\AAA,\bbb)$  and $-\dertp Y: \RHom_\ccm(\DD(\AAA),\bbb) \ra \RHom(\AAA,\bbb)$ is a quasi-equivalence by \Cref{kanextension}, we have that there exists a $\sigma' \in \zzz^n(\RHom_\ccm(\DD(\AAA),\bbb))(F,F')$ such that $H^n(-\dertp Y)([\sigma']) = [\sigma]$. 	Next observe that, because $[\phi] \in H^0(\mathrm{Ar}_0,\bbb)(C(F),I(B))$ is an isomorphism, we can consider a $0$-cycle $\psi \in \RHom(\mathrm{Ar}_0,\bbb)(I(B),C(F))$ such that $[\psi] = [\phi]^{-1}$. Therefore, there exists an $\alpha \in \RHom^{-1}(\mathrm{Ar}_0,\bbb)(C(F),C(F))$ such that $d(\alpha) = \id_{C(F)} - \psi \phi$.
		We define $\mu \in \qrep(\bbb)(B,B')$ as the composite
		\begin{equation*}
			\begin{tikzcd}[column sep= 1.2cm]
				B \arrow[r,"\psi_1"]& F(\coprod_i h_{A_i}) \arrow[r,"\sigma'_{\coprod_i h_{A_i}}"] &F'(\coprod_i h_{A_i}) \arrow[r,"\phi'_1"] &B',
			\end{tikzcd}
		\end{equation*}
		and $\nu \in \RHom^{n-1}(\mathrm{Ar}_0,\bbb)(C(F),I(B'))$ as the composite
		\begin{equation*}
			\begin{tikzcd}[column sep= 1.2cm]
				C(F) \arrow[r,"- \alpha"]& C(F) \arrow[r,"C(\sigma')"] &C(F') \arrow[r,"\phi'"] &I(B').
			\end{tikzcd}
		\end{equation*}
		We claim that $(\sigma',\mu,\nu) \in \zzz^n(\ppp)((F,B,\phi),(F',B',\phi'))$. Indeed, we have that
		\begin{equation*}
			\begin{aligned}
				d(\sigma',\mu,\nu) &=d \left( \sigma',\phi'_1 \sigma'_{\coprod_i h_{A_i}}  \psi_1, - (\phi' C(\sigma')  \alpha)\right) \\
				&= \left( 0,0, -d(\phi'  C(\sigma') \alpha) +(-1)^n \left( \phi'C(\sigma') - I(\phi'_1 \sigma'_{\coprod_i h_{A_i}}  \psi_1)\phi \right) \right)\\
				&= \left( 0,0, -(-1)^n \phi'C(\sigma') d(\alpha) +(-1)^n \left( \phi'C(\sigma') - \phi'C(\sigma') \psi \phi \right) \right)\\
				&= \left( 0,0, -(-1)^n \phi'C(\sigma') (\id_{C(F)} - \psi \phi) +(-1)^n \phi'C(\sigma') \left( \id_{C(F)} - \psi \phi \right) \right)\\ 
				&=(0,0,0),
			\end{aligned}
		\end{equation*} 
		where the third equality follows from
		\begin{equation*}
			\begin{aligned}
				I(\phi'_1 \sigma'_{\coprod_i h_{A_i}} \psi_1) \phi &= \left( \phi'_1 \sigma'_{\coprod_i h_{A_i}} \psi_1 \phi_1, \phi'_1 \sigma'_{\coprod_i h_{A_i}} \psi_1 \phi_2 \right)\\
				&= \left( \phi'_1 \sigma'_{\coprod_i h_{A_i}} \psi_1 \phi_1, \phi'_2 F'(\mathrm{can}_{\lambda})  \sigma'_{\coprod_i h_{A_i}} \psi_1 \phi_2 \right)\\
				&= \left( \phi'_1 \sigma'_{\coprod_i h_{A_i}} \psi_1 \phi_1, \phi'_2 \sigma'_{h_{\coprod_i^{H^0} A_i}} F(\mathrm{can}_\lambda) \psi_1 \phi_2 \right)\\
				&=\left( \phi'_1 \sigma'_{\coprod_i h_{A_i}} \psi_1 \phi_1, \phi'_2 \sigma'_{h_{\coprod_i^{H^0} A_i}} \psi_2 \phi_2 \right)\\
				&=\phi'C(\sigma')\psi\phi.
			\end{aligned}
		\end{equation*}
		By construction, one readily sees that $[(\sigma',\mu,\nu)] \in H^n(\ppp)((F,B,\phi),(F',B',\phi'))$ gets sent to $[\sigma] \in H^n(\eee_{\lambda})(F\dertp Y, F'\dertp Y)$ via $H^n(S)$, proving that $S$ is quasi-full as desired. 
%
		
		To finish the argument, it remains to show that $S$ is quasi-faithful. Consider $(\gamma,\mu,\nu) \in Z^n \ppp ((F,B,\phi),(F',B',\phi'))$ such that $[(\gamma,\mu,\nu)] \in H^n \ppp ((F,B,\phi),(F',B',\phi'))$ gets sent to 0 via
		$$H^n \ppp ((F,B,\phi),(F',B',\phi')) \ra H^n \eee_{\lambda}(F\dertp Y,F'\dertp Y).$$ 
		In what follows, we denote by $\phi_F \in \RHom^0(\mathrm{Ar}_0,\bbb)(C(F),F(h_{\coprod_i^{H^0}A_i}))$ the natural morphism $(F(\mathrm{can}_\lambda),\id_{F(h_{\coprod^{H^0}_i A_i})})$. Notice that $\phi_F$ is closed and induces an isomorphism in $H^0(\RHom(\mathrm{Ar}_0,\bbb))$.
		First, one observes that the morphism 
		$$(\id_F,\phi_2,0): (F,F(h_{\coprod^{H^0}_i A_i}),\phi_F) \ra (F,B,\phi)$$
	 	is a homotopy equivalence by using the characterization of homotopy equivalences in $\ppp$ provided in \cite[Lem 4.2]{milnor-descent-cohesive-dg-categories}. 
		Consequently, in order to conclude that $[(\gamma,\mu,\nu)] = 0$ is enough to show that $[(\gamma,\mu,\nu)(\id_F,\phi_2,0)] = 0$.
		We have that
		$$(\gamma,\mu,\nu)(\id_F,\phi_2,0) = (\gamma,  \mu \phi_2, (-1)^0 \nu C(\id_F) + I(\mu) 0) = (\gamma,  \mu \phi_2, \nu)$$ 
		in $Z^n \ppp ((F,F(h_{\coprod^{H^0}_i A_i}),\phi_F),(F',B',\phi'))$.
		Therefore, it suffices to show that there exists an $(\alpha,\beta,\delta) \in \ppp^{n-1}((F,F(\coprod^{H^0}_i A_i),\phi_F),(F',B',\phi'))$ such that $d(\alpha,\beta,\delta) = (\gamma, \mu \phi_2, \nu)$. First observe that $[\gamma \dertp Y] = 0$ in $H^n\eee_{\lambda}(F\dertp Y,F'\dertp Y)$ by hypothesis, and hence, $[\gamma] = 0$ in $H^n(\RHom_\ccm(\DD(\AAA),\bbb))(F,F')$. Thus, there exists an element $\alpha \in \RHom_{\ccm}(\DD(\AAA),\bbb))^{n-1}(F,F')$ such that $d\alpha = \gamma$. Our candidate $(\alpha,\beta,\delta)$ is going to be $(\alpha, \phi'_2 \alpha_{h_{\coprod_i^{ H^0}A_i}} + (-1)^n\nu_2, 0)$. First we compute $d\nu_2$. As $(\gamma,  \mu, \nu)$ is a $n$-cycle, we have that $d(\gamma,  \mu , \nu) = 0$, in particular, this implies that $0 = d\nu + (-1)^n (\phi' C(\gamma) -I(\mu) \phi)$, that is
		\begin{equation*}\label{dnu2}
			(d\nu_1,d\nu_2) = \left( (-1)^n (- \phi'_1 \gamma_{\coprod_i h_{A_i}} + \mu \phi_1), (-1)^n(- \phi'_2 \gamma_{h_{\coprod^{H^0}_i A_i}} + \mu \phi_2)  \right) . 
		\end{equation*}
		Making use of this, we can compute now:
		\begin{equation*}
		\begin{aligned}
		d\beta = d(\phi'_2 \alpha_{h_{\coprod_i^{H^0} A_i}} + (-1)^n\nu_2) &= \phi'_2\gamma_{h_{\coprod^{H^0}_i A_i}} + (-1)^n d\nu_2\\
		&= \phi'_2\gamma_{h_{\coprod^{H^0}_i A_i}} +(-1)^n ((-1)^n(- \phi'_2 \gamma_{h_{\coprod^{H^0}_i A_i}} + \mu \phi_2) )\\
		&= \mu \phi_2.
		\end{aligned}
		\end{equation*}
		Consequently, we have that
		\begin{equation*}
			\begin{aligned}
			d(\alpha,\beta,\delta) &= d(\alpha, \phi'_2 \alpha_{h_{\coprod_i^{H^0} A_i}} + (-1)^n\nu_2, 0)\\
			&= \left( \gamma,\mu \phi_2,(-1)^{n-1}( \phi' C(\alpha) - I (\phi'_2 \alpha_{h_{\coprod^{H^0}_i A_i}} + (-1)^n\nu_2) \phi_F) \right),
			\end{aligned}
		\end{equation*}
		where the last component is given by:
		\begin{equation*}
			\begin{aligned}
			(-1)^{n-1}&( \phi' C(\alpha) - I (\phi'_2 \alpha_{h_{\coprod^{H^0}_i A_i}} + (-1)^n\nu_2) \phi_F)=\\
			&= (-1)^{n-1}\left(  \phi'_1 \alpha_{\coprod_i h_{A_i}} - \phi'_2 \alpha_{h_{\coprod_i^{H^0} A_i}}F(\mathrm{can}_\lambda) - (-1)^n \nu_2 F(\mathrm{can}_\lambda), \right. \\ 
			& \,\,\,\,\,\, \left. \phi'_2 \alpha_{h_{\coprod_i^{H^0} A_i}} - \phi'_2 \alpha_{h_{\coprod_i^{H^0} A_i}} - (-1)^n\nu_2 \right)\\
			&= (-1)^{n-1}\left(  \phi'_1 \alpha_{\coprod_i h_{A_i}} - \phi'_2 \alpha_{h_{\coprod_i^{H^0} A_i}}F(\mathrm{can}_\lambda) - (-1)^n \nu_2 F(\mathrm{can}_\lambda), -(-1)^n \nu_2 \right)\\
			&= (-1)^{n-1}\left(  \phi'_1 \alpha_{\coprod_i h_{A_i}} - \phi'_1 \alpha_{\coprod_i h_{A_i}} - (-1)^n \nu_1,-(-1)^n \nu_2 \right)\\
			&=(-1)^{n-1} \left( - (-1)^n \nu_1,-(-1)^n \nu_2 \right) = (\nu_1,\nu_2) = \nu
			\end{aligned}
		\end{equation*}
		We hence have that $d(\alpha,\beta,\delta) = (\gamma,\mu \phi_2,\nu)$ as desired.
	\end{proof}
\end{theorem}

\section{The well generated tensor product}\label{quotandtpwg}

Let $\hodgcat_{\wg}$ denote the subcategory of $\frakv$-$\hodgcat$ given by the $\frakv$-small $\fraku$-well generated dg categories with cocontinuous quasi-functors. Up to equivalence, $\hodgcat_{\wg}$ is easily seen to be independent of the choice of $\frakv$.

\begin{definition}\label{defwellgeneratedtp}
	Let $\aaa$ and $\bbb$ be well generated dg categories. A \emph{well generated tensor product} of $\aaa$ and $\bbb$ is defined as a well generated dg category $\aaa \boxtimes \bbb$ such that for every well generated dg category $\ccc$, the following universal property holds:
	\begin{equation}\label{eqboxdgwg}
	\RHom_{\ccm}(\aaa \boxtimes \bbb, \ccc) \cong \RHom_{\ccm}(\aaa, \RHom_{\ccm}(\bbb, \ccc)).
	\end{equation}
\end{definition}

As a consequence, by \Cref{closed},  if we can show that the tensor product of well generated dg categories exists, the resulting monoidal structure on $\hodgcat_{\wg}$ is closed.

\begin{remark}
Note that the situation is different from the one for Grothendieck categories. As shown in \cite[Thm 5.4]{tensor-product-linear-sites-grothendieck-categories}, the tensor product of locally presentable $k$-linear categories is closed under Grothendieck categories, but the natural inner hom of cocontinuous functors between locally presentable categories is not (as follows for instance from \cite[Rem 6.5]{covers-envelopes-cotorsion-theories-locally-presentable-abelian-categories-contramodule-categories}). However, by Corollary \ref{corwglocpres}, the distinction between locally presentable categories and localisations of module categories does not exist on the derived level, whence this subtlety vanishes. An in depth study of the nature of morphisms categories between abelian categories is the topic of an ongoing joint project with Michel Van den Bergh.
\end{remark}

The rest of the paper is devoted to proving that the well generated tensor product exists (Theorem  \ref{thmexist}), and providing various constructions using localisation theory. In particular, \S \ref{partensdgquot} and \S \ref{partensorproductwg} discuss the relation between the tensor product and the dg quotient, in \S \ref{tplocsubcats} the tensor product is described in terms of localising subcategories of dg derived categories, and in \S \ref{tpstrictlocs} the tensor product is described in terms of their Bousfield localisations.

We start with some considerations regarding the internal hom in two variables in \S \ref{par2var}.

\subsection{Considerations in the two variable setting}\label{par2var}
We devote this section to prove that both ($\alpha$-)cocontinuity and annihilation of classes of objects behave suitably with respect to the monoidal structure. From now on, and for the rest of the paper, we will make implicit use of the fact that for every homotopically cocomplete small dg category, we can pick a cofibrant replacement in $\hodgcat$ which is also homotopically cocomplete (homotopically cocompleteness is preserved under quasi-equivalences) and this cofibrant replacement is the identity on objects (see \Cref{cofibrantreplacement} above). 

Let $\aaa, \bbb, \ccc$ be dg categories. Consider a right quasi-representable bimodule $F \in \RHom(\aaa \dertp \bbb, \ccc)$ and observe that the dg module $F \in \dgmod(\ccc \dertp \aaa^{\op} \dertp \bbb^{\op})$ with evaluations $F(C, A, B)$ gives rise on one hand to a bimodule $F_A = F(-, A, -) \in \dgmod(\ccc \dertp \bbb^{\op})$ for every $A \in \aaa$ and on the other hand to a bimodule $F_B = F(-, -, B) \in \dgmod(\ccc \dertp \aaa^{\op})$ for every $B \in \bbb$, and according to \eqref{eq2} these are all right quasi-representable. 

\begin{definition}
	We call $F \in \RHom(\aaa \dertp \bbb, \ccc)$ \emph{right cocontinuous} provided that every $F_B$ is cocontinuous, \emph{left cocontinuous } provided that every $F_A$ is cocontinuous, and \emph{bicocontinuous} provided that it is left and right cocontinuous. 
\end{definition}
	We denote by $\RHom_{\ccm,\ccm}(\aaa \dertp \bbb, \ccc) \subseteq \RHom(\aaa \dertp \bbb, \ccc)$ the full dg subcategory of bicocontinuous modules.
	
	Given a regular cardinal $\alpha$, the notions of \emph{left-, right- and bi-$\alpha$-cocontinuous} are defined similarly. In particular, we denote by $\RHom_{\alpha, \alpha}(\aaa \dertp \bbb, \ccc) \subseteq \RHom(\aaa \dertp \bbb, \ccc)$ the full dg subcategory of bi-$\alpha$-cocontinuous right quasi-representable bimodules.

\begin{definition}
	Consider $\nnn_{\aaa}$ a class of objects in $\aaa$ and $\nnn_{\bbb}$ a class of objects in $\bbb$. With the same notations as above, we say $F \in \RHom(\aaa \dertp \bbb, \ccc)$ \emph{biannihilates $(\nnn_{\aaa},\nnn_{\bbb})$} provided that every $F_A$ annihilates $\nnn_{\bbb}$ and every $F_B$ annihilates $\nnn_{\aaa}$. 
\end{definition}	
	We denote by $\RHom_{\nnn_{\aaa}, \nnn_{\bbb}}(\aaa \dertp \bbb, \ccc) \subseteq \RHom(\aaa \dertp \bbb, \ccc)$ the full dg subcategory of right quasi-representable bimodules that biannihilate $(\nnn_{\aaa},\nnn_{\bbb})$.
	
	Similarly, we denote by $\RHom_{(\ccm,\nnn_{\aaa}),(\ccm,\nnn_{\bbb}) }(\aaa \dertp \bbb, \ccc) \subseteq \RHom(\aaa \dertp \bbb, \ccc)$ the full dg subcategory of bicocontinuous right quasi-representable bimodules that biannihilate $(\nnn_{\aaa},\nnn_{\bbb})$.

We include the proof of the following statement for the convenience of the reader.
\begin{lemma}\label{innerhom} 
Let $\alpha \leq |\fraku|$ be a regular cardinal. The following hold:	
	\begin{enumerate}
		\item For homotopically cocomplete dg categories $\aaa$, $\bbb$ and $\ccc$, we have that the equivalence (\ref{eq2}) restricts to:
		\begin{equation}\label{innerhom1}
		\RHom_{\ccm,\ccm}(\aaa \dertp \bbb, \ccc) \cong \RHom_{\ccm}(\aaa, \RHom_{\ccm}(\bbb, \ccc))
		\end{equation}
		\item 
		For homotopically $\alpha$-cocomplete dg categories $\aaa$, $\bbb$ and $\ccc$, we have that the equivalence (\ref{eq2}) restricts to:
		\begin{equation}
		\RHom_{\alpha,\alpha}(\aaa \dertp \bbb, \ccc) \cong \RHom_{\alpha}(\aaa, \RHom_{\alpha}(\bbb, \ccc))
		\end{equation}
		\item For homotopically cocomplete dg categories $\aaa$, $\bbb$ and $\ccc$ and sets of objects $\nnn_{\aaa}$ in $\aaa$ and $\nnn_{\bbb}$ in $\bbb$, we have that the equivalence (\ref{eq2}) restricts to:
		\begin{equation} \label{innerhom3}
		\RHom_{(\ccm,\nnn_{\aaa}),(\ccm,\nnn_{\bbb}) }(\aaa \dertp \bbb, \ccc) \cong \RHom_{\ccm,\nnn_{\aaa}}(\aaa, \RHom_{\ccm, \nnn_{\bbb}}(\bbb, \ccc))
		\end{equation}
	\end{enumerate}
	\begin{proof} 
		Observe that (1) is the case $\alpha = |\fraku |$ of (2). We prove (2). First we show that for any $F \in \RHom_{\alpha,\alpha}(\aaa \dertp \bbb, \ccc) \subseteq \RHom(\aaa \dertp \bbb,\ccc)$ the image of $F$ via \eqref{eq2} is an element of $\RHom_{\alpha}(\aaa,\RHom_{\alpha}(\bbb,\ccc))$.
		If we denote by $\bar{F}$ the image of $F$ in $\RHom(\aaa,\RHom(\bbb,\ccc))$ via \eqref{eq2} we have that $\bar{F}(A) = F_A$ factors through $\RHom_{\alpha}(\bbb,\ccc) \subseteq \RHom(\bbb,\ccc)$ by hypothesis. We hence have that $\bar{F}$ belongs to $\RHom(\aaa,\RHom_{\alpha}(\bbb,\ccc))$. Let's now show that $\bar{F}$ actually belongs to $\RHom_{\alpha}(\aaa, \RHom_{\alpha}(\bbb,\ccc))$. Let $\{A_i\}_{i \in I}$ be an $\alpha$-small family of objects in $\aaa$. By definition, we have that
		\begin{equation}
			H^0(\bar{F})\left( \coprod_{i\in I}^{H^0(\aaa)} A_i\right)  = F_{\coprod_{i\in I}^{H^0(\aaa)} A_i} \in H^0(\RHom_{\alpha}(\bbb,\ccc)).
		\end{equation}
		For all $i \in I$ we have a natural morphism 
		\begin{equation}
		F_{\coprod_{i\in I}^{H^0(\aaa)} A_i} = H^0(\bar{F})(\coprod_{i\in I}^{H^0(\aaa)} A_i)  \longleftarrow H^0(\bar{F})(A_i)= F_{A_i},
		\end{equation}
		in $H^0(\RHom_{\alpha}(\bbb,\ccc))$, and hence we have the natural morphism
		\begin{equation} \label{naturalmorphismcoproduct}
		\coprod_{i\in I}^{H^0(\RHom_{\alpha}(\bbb,\ccc))}F_{A_i} \lra F_{\coprod_{i\in I}^{H^0(\aaa)} A_i}
		\end{equation}
		in $H^0(\RHom_{\alpha}(\bbb,\ccc))$, induced by the universal property of the coproduct. We claim that this morphism is an isomorphism. Indeed, observe that for all $B\in \bbb$, we have
		\begin{equation}
			 \begin{aligned}
			 	H^0 &\left(  \coprod_{i\in I}^{H^0(\RHom_{\alpha}(\bbb,\ccc))} F_{A_i}\right)  (B) =\\ &= \coprod_{i\in I}^{H^0(\ccc)} H^0(F_{A_i})(B) = \coprod_{i\in I}^{H^0(\ccc)} H^0(F(A_i,B,-)) = \coprod_{i\in I}^{H^0(\ccc)} H^0(F_{B})(A_i)=\\
			 	&= H^0(F_B)\left( \coprod_{i\in I}^{H^0(\aaa)}A_i\right)  = H^0\left( F\left( \coprod_{i\in I}^{H^0(\aaa)}A_i,B,-\right) \right)  = H^0\left( F_{\coprod_{i\in I}^{H^0(\aaa)}A_i}\right) (B),
			 \end{aligned}
		\end{equation}
		functorially in $B \in \bbb$, where the first equality follows from \Cref{cocompletecontinuousRHom} and the fourth from the fact that $F_B \in \RHom_{\alpha}(\aaa, \ccc)$. It follows that \eqref{naturalmorphismcoproduct} is an isomorphism. Consequently, we have that 	
		\begin{equation}
		H^0(\bar{F})(\coprod_{i\in I}^{H^0(\aaa)}A_i) = F_{\coprod_{i\in I}^{H^0(\aaa)}A_i} \cong \coprod_{i\in I}^{H^0(\RHom_{\alpha}(\bbb,\ccc))}F_{A_i} = \coprod_{i\in I}^{H^0(\RHom_{\alpha}(\bbb,\ccc))} H^0(\bar{F})(A_i)
		\end{equation}
		in $H^0(\RHom_{\alpha}(\bbb,\ccc))$, as desired. 
		
		To conclude it is enough to prove that for any $F \in \RHom(\aaa \dertp \bbb, \ccc)$, if its image $\bar{F}$ via \eqref{eq2} belongs to $\RHom_{\alpha}(\aaa,\RHom_{\alpha}(\bbb,\ccc))$, then $F$ lies in $\RHom_{\alpha,\alpha}(\aaa \dertp \bbb, \ccc)$. Take such an $F$. By definition, for every $A \in \aaa$ we have that
		$$\bar{F}(A) = F_A \in \RHom_{\alpha}(\bbb,\ccc),$$
		which proves the $\alpha$-cocontinuity of $F_A$ for all $A \in \aaa$. Let $\{A_i\}_{i \in I}$ an $\alpha$-small family of objects in $\aaa$. For every $B \in \bbb$ we have that 
		\begin{equation*}
			\begin{aligned}
			 H^0&(F_B) \left( \coprod_{i\in I}^{H^0(\aaa)}A_i \right) =\\
			 &= H^0 \left( F\left( \coprod_{i\in I}^{H^0(\aaa)}A_i,B,-\right) \right)  = H^0\left( F_{\coprod_{i\in I}^{H^0(\aaa)}A_i}\right) (B) = H^0\left( H^0(\bar{F})\left(\coprod_{i\in I}^{H^0(\aaa)}A_i\right) \right) (B) =\\
			 &= H^0\left( \coprod_{i\in I}^{H^0(\RHom_{\alpha}(\bbb,\ccc))} H^0(\bar{F})(A_i)\right) (B) = \coprod_{i\in I}^{H^0(\ccc)}(H^0(F_{A_i})(B)) = \coprod_{i\in I}^{H^0(\ccc)} H^0(F_B)(A_i), 
			\end{aligned}
		\end{equation*}
		where the fourth equality uses the fact that $\bar{F} \in \RHom_{\alpha}(\aaa,\RHom_{\alpha}(\bbb,\ccc))$ and the fith follows from \Cref{cocompletecontinuousRHom}. This proves the $\alpha$-cocontinuity of $F_B$ for all $B \in \bbb$. We can thus conclude that $F \in \RHom_{\alpha,\alpha}(\aaa \dertp \bbb,\ccc)$ as we wanted to show. 
		
		We prove (3). It is enough to see that the isomorphism \eqref{innerhom1} in $\hodgcat$ constructed above restricts to an isomorphism \eqref{innerhom3}. Let $F \in \RHom_{(\ccm, \nnn_{\aaa}), (\ccm, \nnn_{\bbb})}(\aaa \dertp \bbb,\ccc)$ and denote by $\bar{F}$ its image in $\RHom_{\ccm}(\aaa,\RHom_{\ccm}(\bbb,\ccc))$ via \eqref{innerhom1}. Then we have that 
		$$H^0\left( H^0(\bar{F})(A)\right) (B) = H^0(F_A)(B) = 0$$
		for all $B \in \nnn_{\bbb}$ and hence $\bar{F} \in \RHom_{\ccm}(\aaa,\RHom_{\ccm, \nnn_{\bbb}}(\bbb,\ccc))$. Now observe that, for all $B \in \bbb$, we have that
		$$H^0\left( H^0(\bar{F})(A)\right) (B) = H^0(F_A)(B) = H^0(F(A,B,-))= H^0(F_B)(A) = 0$$
		for all $A \in \nnn_{\aaa}$. Consequently, we have that $H^0(\bar{F})(A) = 0$ in $H^0(\RHom_{\ccm, \nnn_{\bbb}}(\bbb,\ccc))$ for all $A \in \nnn_{\aaa}$ and hence $\bar{F} \in \RHom_{\ccm,\nnn_{\aaa}}(\aaa,\RHom_{\ccm,\nnn_{\bbb}}(\bbb,\ccc))$ as desired.
		
		To conclude, it is enough to show that for all $F \in \RHom_{\ccm,\ccm}(\aaa \dertp \bbb, \ccc)$, if the image $\bar{F}$ of $F$ via \eqref{innerhom1} belongs to $\RHom_{\ccm,\nnn_{\aaa}}(\aaa,\RHom_{\ccm,\nnn_{\bbb}}(\bbb,\ccc))$, then $F$ is an element of $\RHom_{(\ccm, \nnn_{\aaa}), (\ccm, \nnn_{\bbb})}(\aaa \dertp \bbb,\ccc)$. For all $A \in \aaa$, we have that
		$$H^0(F_A) (B)= H^0\left( H^0(\bar{F})(A) \right) (B) = 0 $$
		for all $B\in \nnn_{\bbb}$, showing that, for all $A \in \aaa$, $F_A$ annihilates $\nnn_{\bbb}$. On the other hand, for all $B \in \bbb$, we have that
		$$H^0(F_B)(A) = H^0(F(A,B,-)) = H^0\left( H^0(\bar{F}) (A)\right) (B) = 0 $$
		for all $A \in \nnn_{\aaa}$, showing that, for all $B\in \bbb$, $F_B$ annihilates $\nnn_{\aaa}$ as desired.
	\end{proof}
	
\end{lemma}

\subsection{The tensor product of dg quotients}\label{partensdgquot}
Consider $\aaa, \bbb, \ccc \in$ $\hodgcat_{\wg}$ and suppose $\aaa \boxtimes \bbb$ exists. By \Cref{innerhom} above and the universal property of $\boxtimes$, we have an isomorphism in $\hodgcat$
\begin{equation}\label{twotensorproductswg}
\RHom_{\ccm}(\aaa \boxtimes \bbb,\ccc) \cong \RHom_{\ccm,\ccm}(\aaa \dertp \bbb,\ccc),
\end{equation}
for every well generated dg category $\ccc$.
Hence there exists, corresponding to the identity quasi-representable module on the left hand side by taking $\ccc= \aaa \boxtimes \bbb$, a canonical bi\-cocontinuous quasi-representable module $\otimes \in H^0(\RHom_{\ccm,\ccm}(\aaa \dertp \bbb,\aaa \boxtimes \bbb))$. We will denote the induced functor at the level of homotopy by $$\hztp: H^0(\aaa \dertp \bbb) \lra H^0(\aaa \boxtimes \bbb),$$  instead of our usual notation $H^0(\otimes)$.
Let $\xxx_{\aaa} \subseteq \aaa$ and $\xxx_{\bbb} \subseteq \bbb$ be classes of objects. We define the class 
\begin{equation}
\xxx_{\aaa} \hztp \xxx_{\bbb} = \{X_A \hztp X_B \,\,|\,\, X_A \in \xxx_A, X_B \in \xxx_B \}
\end{equation} 
of objects in $H^0(\aaa \boxtimes \bbb)$. 

\begin{remark}
	Let $\ccc$ be a dg category. Observe that taking a class of objects in $\ccc$ is the same as taking a class of objects in $H^0(\ccc)$ as  $\Obj(H^0(\ccc)) = \Obj(\ccc)$. 
\end{remark}

In first place, let's analyse the relation of the well generated tensor product and the annihilation of classes of objects in $\hodgcat_{\wg}$.
\begin{proposition}\label{proptwosidedobjectswg} 

	Consider classes $\nnn_{\aaa} \subseteq H^0(\aaa)$ and $\nnn_{\bbb} \subseteq H^0(\bbb)$ of objects. The class
	\begin{equation}
	\nnn_{\aaa} \boxtimes_{\mathrm{Cl}} \nnn_{\bbb} = (\nnn_{\aaa} \hztp \bbb) \cup (\aaa \hztp \nnn_{\bbb}) \subseteq H^0(\aaa \boxtimes \bbb)
	\end{equation}
	is such that
	\begin{equation}\label{tensorandquotientwg}
	\RHom_{\ccm, \nnn_{\aaa} \boxtimes_{\mathrm{Cl}} \nnn_{\bbb}}(\aaa \boxtimes \bbb, \ccc) \cong \RHom_{(\ccm,\nnn_{\aaa}),(\ccm,\nnn_{\bbb}) }(\aaa \dertp \bbb, \ccc).
	\end{equation}
	\begin{proof}
		We have the isomorphism in $\hodgcat$
		\begin{equation*}
			\RHom_{\ccm}(\aaa \boxtimes \bbb,\ccc) \cong \RHom_{\ccm,\ccm}(\aaa \dertp \bbb,\ccc)
		\end{equation*}
		from \eqref{twotensorproductswg} given at the $H^0$-level by composition with the canonical bicocontinuous quasi-rep\-re\-sent\-able bimodule $\otimes$ between $\aaa \dertp \bbb$ and $\aaa \boxtimes \bbb$. Then it is enough to see that this isomorphism restricts to an isomorphism \eqref{tensorandquotientwg} in $\hodgcat$. 
		
		Consider $F \in \RHom_{\ccm, \nnn_{\aaa} \boxtimes_{\mathrm{Cl}} \nnn_{\bbb}}(\aaa \boxtimes \bbb, \ccc)$. Then $F \dertp_{\aaa \boxtimes \bbb} \otimes \in H^0(\RHom_{\ccm,\ccm}(\aaa \dertp \bbb,\ccc))$ is trivially seen to biannihilate ($\nnn_{\aaa}$-$\nnn_{\bbb}$).
		
		On the other hand, given any $G\in \RHom_{(\ccm, \nnn_{\aaa}), (\ccm, \nnn_{\bbb})}(\aaa \dertp \bbb, \ccc)$, we have that $$G \cong F \dertp_{\aaa \boxtimes \bbb} \otimes \in H^0(\RHom_{(\ccm, \nnn_{\aaa}), (\ccm, \nnn_{\bbb})}(\aaa \dertp \bbb, \ccc))$$ for some $F \in \RHom_{\ccm}(\aaa \boxtimes \bbb,\ccc)$. Consequently, for every object $B \in \bbb$, we have that $H^0(F)(\nnn_A \hztp B) \cong H^0(G)(\nnn_{\aaa},B) = 0$ in $H^0(\ccc)$ and, similarly, for every object $A \in \aaa$, $H^0(F)(A \hztp \nnn_{\bbb}) \cong H^0(G)(A,\nnn_{\bbb})= 0$ in $H^0(\ccc)$. Thus we have that $H^0(F)$ annihilates ${\nnn_{\aaa}\boxtimes_{\mathrm{Cl}} \nnn_{\bbb}}$, therefore $F \in \RHom_{\ccm,\nnn_{\aaa} \boxtimes_{\mathrm{Cl}} \nnn_{\bbb}}(\aaa \boxtimes \bbb, \ccc)$ as desired.
	\end{proof}	
\end{proposition}
\begin{definition}
	We will call $\nnn_{\aaa} \boxtimes_{\mathrm{Cl}} \nnn_{\bbb}$ the \emph{tensor product of classes of objects} $\nnn_{\aaa}$ and $\nnn_{\bbb}$.
\end{definition}

\begin{remark}
	Let $\bbb$, $\ccc$ be well generated dg categories and let $\nnn$ be a class of objects in $\bbb$. Let $\langle \nnn \rangle \subseteq H^0(\bbb)$ be the smallest localising subcategory containing $\nnn$. Then, given $F \in \RHom_{\ccm}(\bbb, \ccc)$, the induced $H^0(F): H^0(\bbb) \lra H^0(\ccc)$ is exact and cocontinuous.
	As a consequence, $\Kern(H^0(F))$ is a localising subcategory of $H^0(\bbb)$. It follows that 
	\begin{equation}\label{localisingandnot}
	\RHom_{\ccm,\nnn}(\bbb, \ccc) = \RHom_{\ccm,\langle \nnn \rangle}(\bbb, \ccc).
	\end{equation}
\end{remark}

\begin{lemma}\label{lemgen}
	Let $\aaa$, $\bbb$ be two well generated dg categories and $\www_{\aaa} \subseteq H^0(\aaa)$ and $\www_{\bbb} \subseteq H^0(\bbb)$ localising subcategories generated by sets. Let $\GGG_{\aaa}$ (resp. $\GGG_{\bbb}$) be a set of generators of $H^0(\aaa)$ (resp. $H^0(\bbb)$) and $\nnn_{\aaa}$ (resp. $\nnn_{\bbb}$) be a set of generators of $\www_{\aaa}$ (resp. $\www_{\bbb}$). We have that:
	$$\langle \www_{\aaa} \boxtimes_{\mathrm{Cl}} \www_{\bbb} \rangle =  \langle (\nnn_{\aaa} \hztp \GGG_{\bbb}) \cup (\GGG_{\aaa} \hztp \nnn_{\bbb}) \rangle.$$
	Hence $\langle \www_{\aaa} \boxtimes_{\mathrm{Cl}} \www_{\bbb} \rangle$ is generated by a set of objects.
	\begin{proof}
		By definition we have that
		$$\langle \www_{\aaa} \boxtimes_{\mathrm{Cl}} \www_{\bbb} \rangle = \langle (\langle \nnn_{\aaa} \rangle \hztp \bbb) \cup (\aaa \hztp \langle \nnn_{\bbb} \rangle) \rangle.$$
		As it is a localising subcategory and it trivially contains $\nnn_{\aaa} \hztp \GGG_{\bbb} \cup \GGG_{\aaa} \hztp \nnn_{\bbb}$, we have that
		$\langle \nnn_{\aaa} \hztp \GGG_{\bbb} \cup \GGG_{\aaa} \hztp \nnn_{\bbb} \rangle \subseteq \langle \www_{\aaa} \boxtimes_{\mathrm{Cl}} \www_{\bbb} \rangle$. 
		
		In order to prove the other inclusion, we consider an element $X \in \langle \nnn_{\aaa} \rangle \hztp \bbb$. If it belonged to $\aaa \hztp \langle \nnn_{\bbb} \rangle$, we argue analogously. We know we can choose regular cardinals $\alpha$ and $\beta$ such that the generators $\nnn_{\aaa}$ are all $\alpha$-compact in $\www_{\aaa}$ and the generators $\GGG_{\bbb}$ are all $\beta$-compact in $\bbb$.   
		Combining \cite[Lem 4.4.5 \& Lem B.1.3]{triangulated-categories}, we have that 
		$$X \cong W \hztp B,$$
		where $W$ (resp. $B$) can be written in terms of objects in $\nnn_{\aaa}$ (resp. in $\GGG_{\bbb}$) by using coproducts, cones, direct summands and shifts.
		As $\hztp$ is bicocontinuous and an exact functor in each variable, we have that $X$ can be written using coproducts, cones, direct summands and shifts in terms of elements of the form $N \hztp G$ where $N \in \nnn_{\aaa}$ and $G \in \GGG_{\bbb}$. Therefore, $X$ is an element of $\langle \nnn_{\aaa} \hztp \GGG_{\bbb} \cup \GGG_{\aaa} \hztp \nnn_{\bbb} \rangle$. Consequently, we also have an inclusion $\langle \www_{\aaa} \boxtimes_{\mathrm{Cl}} \www_{\bbb} \rangle \subseteq  \langle \nnn_{\aaa} \hztp \GGG_{\bbb} \cup \GGG_{\aaa} \hztp \nnn_{\bbb} \rangle$ which concludes the proof.
	\end{proof}
\end{lemma}

\begin{theorem}\label{thmboxquotwg}
Let $\aaa$, $\bbb$ be two  well generated dg categories such that $\aaa \boxtimes \bbb$ exists, and consider $\www_{\aaa} \subseteq H^0(\aaa)$ and $\www_{\bbb} \subseteq H^0(\bbb)$ localising subcategories generated by sets.
	We have
	\begin{equation}
	\frac{\aaa}{\www_{\aaa}} \boxtimes \frac{\bbb}{\www_{\bbb}} \cong \frac{\aaa \boxtimes \bbb}{\langle\www_{\aaa} \boxtimes_{\mathrm{Cl}} \www_{\bbb}\rangle}
	\end{equation}
in $\hodgcat_{\wg}$.
	\begin{proof}
		The subcategory $\langle\www_{\aaa} \boxtimes_{\mathrm{Cl}} \www_{\bbb}\rangle \subseteq H^0(\aaa \boxtimes \bbb)$ is a localising subcategory generated by a set of objects as proved in \Cref{lemgen}. Hence, $\aaa \boxtimes \bbb/\langle\www_{\aaa} \boxtimes_{\mathrm{Cl}} \www_{\bbb}\rangle$ is a well generated dg category. If we show that it satisfies the universal property \eqref{eqboxdgwg}, we conclude our argument. For any well generated dg category $\ccc$, we have:
		$$\begin{aligned}
		\RHom_{\ccm}(\frac{\aaa \boxtimes \bbb}{\langle\www_{\aaa} \boxtimes_{\mathrm{Cl}} \www_{\bbb}\rangle}, \ccc) &\cong \RHom_{\ccm, \langle\www_{\aaa} \boxtimes_{\mathrm{Cl}} \www_{\bbb}\rangle}(\aaa \boxtimes \bbb, \ccc)\\
		&\cong \RHom_{\ccm, \www_{\aaa} \boxtimes_{\mathrm{Cl}} \www_{\bbb}}(\aaa \boxtimes \bbb, \ccc)\\
		& \cong  \RHom_{(\ccm,\www_{\aaa}),(\ccm,\www_{\bbb}) }(\aaa \dertp \bbb, \ccc) \\
		& \cong \RHom_{\ccm,\www_{\aaa}}(\aaa, \RHom_{\ccm,\www_{\bbb}}(\bbb, \ccc)) \\
		& \cong \RHom_{\ccm}(\frac{\aaa}{\www_{\aaa}}, \RHom_{\ccm}(\frac{\bbb}{\www_{\bbb}}, \ccc)),
		\end{aligned}$$
		where the first and last isomorphisms come from the universal property of the dg quotient in $\hodgcat_{\wg}$ (see \eqref{dgquotientrestricted}), the second follows from \eqref{localisingandnot} above, the third one is given by \Cref{proptwosidedobjectswg} and the fourth one by \Cref{innerhom}.
	\end{proof}
\end{theorem}

\begin{corollary}\label{corexisttpquotwg}
	Let $\aaa, \bbb$ be two well generated dg categories. If the tensor product $\aaa \boxtimes \bbb$ exists, so does the well generated tensor product between any two dg quotients of $\aaa$, $\bbb$ with respect to localising subcategories generated by a set of objects.
	\begin{proof}
		This is a direct consequence of \Cref{thmboxquotwg}.
	\end{proof}
\end{corollary}

\subsection{Tensor product of well generated dg categories}\label{partensorproductwg}
In this section we show that the well generated tensor product exists and we provide a construction.

We will proceed as follows. We will show that the well generated tensor product of derived dg categories exists and it is again a derived dg category. This result will allow us, using Theorem \ref{thmenhGP}, to approach the construction of the tensor product for arbitrary well generated dg categories making essential use of \Cref{corexisttpquotwg} above.

\begin{proposition}\label{proptpderivedcats}
	Consider small dg categories $\AAA$ and $\BBB$. In $\hodgcat_{\wg}$, we have
	\begin{equation}
	\DD(\AAA) \boxtimes \DD(\BBB) \cong \DD(\AAA \dertp \BBB).
	\end{equation}
	\begin{proof}
		For a well generated dg category $\ccc$, we have
		$$\begin{aligned}
		\RHom_{\ccm}(\DD(\AAA \dertp \BBB), \ccc) & \cong \RHom(\AAA \dertp \BBB, \ccc) \\
		& \cong \RHom(\AAA, \RHom(\BBB, \ccc)) \\
		& \cong \RHom_{\ccm}(\DD(\AAA), \RHom_{\ccm}(\DD(\BBB), \ccc))
		\end{aligned}$$
		where the first and the last isomorphisms are given by \eqref{eqc} and the second one is by the $\dertp$ - $\RHom$-adjunction in $\hodgcat$.
	\end{proof}
\end{proposition}

We are finally in the position to prove the existence of the well generated tensor product. 
\begin{theorem}\label{thmexist}
	Let $\aaa$, $\bbb$ be two well generated dg categories such that $\aaa \cong \DD(\AAA)/\www_{\AAA}$ and $\bbb \cong \DD(\BBB)/\www_{\BBB}$ in $\hodgcat$ for small dg categories $\AAA$, $\BBB$ with $\www_{\AAA}\subseteq \D(\AAA)$ and $\www_{\BBB}\subseteq \D(\BBB)$ localising subcategories generated by a set of objects. Then, the well generated tensor product of $\aaa$ and $\bbb$ exists and it is given by
	\begin{equation}
		\aaa \boxtimes \bbb \cong \DD(\AAA \dertp \BBB)/\langle\www_{\AAA} \boxtimes_{\mathrm{Cl}} \www_{\BBB} \rangle.
	\end{equation}
	In particular, $\aaa \boxtimes \bbb$ is independent of the chosen realisations of $\aaa$ and $\bbb$.
	\begin{proof}
		We have $\aaa \cong \DD(\AAA)/\www_{\AAA}$ and $\bbb \cong \DD(\BBB)/\www_{\BBB}$ with $\www_{\AAA}$ and $\www_{\BBB}$ localising subcategories generated by a set of objects. By \Cref{proptpderivedcats} we know that $\DD(\AAA) \boxtimes \DD(\BBB)$ exists and equals $\DD(\AAA \dertp \BBB)$. Then, by \Cref{thmboxquotwg}, we have that $\aaa \boxtimes \bbb \cong \DD(\AAA)/\www_{\AAA} \boxtimes \DD(\BBB)/\www_{\BBB}$ exists and it is given by $\DD(\AAA \dertp \BBB)/\langle\www_{\AAA} \boxtimes_{\mathrm{Cl}} \www_{\BBB} \rangle$, and it is obviously independent of the realizations chosen, as it fulfils the universal property.
	\end{proof}
\end{theorem}

\begin{corollary}
The homotopy category $\hodgcat_{\wg}$ of well generated dg categories with cocontinuous quasi-functors is symmetric monoidal closed.
\end{corollary} 

\begin{proof}
This follows from \Cref{thmexist} and \Cref{closed}.
\end{proof}

\subsection{Tensor product of localising subcategories}\label{tplocsubcats}
In this section we provide an alternative description of the tensor product from \S\ref{partensorproductwg}, in the spirit of \cite[\S2.5]{tensor-product-linear-sites-grothendieck-categories}, which does not appeal to choices of generators of localising subcategories. In the next section, this construction will lead, via the equivalent approaches to localisation theory described in \S\ref{equivapproaches}, to a description of the tensor product in terms of Bousfield localisations (in the spirit of \cite[\S 2.6]{tensor-product-linear-sites-grothendieck-categories}), which will be used in \S\ref{tpalphacocontinuous}.

Let $\AAA$, $\BBB$ be two small dg categories and consider the derived dg categories $\DD(\AAA)$ and $\DD(\BBB)$. Let $\www_{\AAA} \subseteq \D(\AAA)$ and $\www_{\BBB} \subseteq \D(\BBB)$ be localising subcategories generated by sets of objects. 
Inspired upon the construction of $\boxtimes$ above, we can define a \emph{tensor product of localising subcategories generated by a set} as follows.
\begin{definition}
	With the notations above, we put
	\begin{equation}\label{deftplocalisingsubcats}
		\www_{\AAA} \boxtimes \www_{\BBB} = \langle \www_{\AAA} \boxtimes_{\mathrm{Cl}} \www_{\BBB}\rangle.
	\end{equation}
\end{definition}
We define one-sided localising subcategories of $\D(\AAA \dertp \BBB)$ as follows:
\begin{equation}\label{onesidedlocsubcats}
\begin{aligned}
\www_1 &\coloneqq \{F \in \D(\AAA \dertp \BBB) | F(-,B) \in \www_{\AAA} \text{ for all } B\in \BBB \}\\
\www_2 &\coloneqq \{F \in \D(\AAA \dertp \BBB) | F(A,-) \in \www_{\BBB} \text{ for all } A\in \AAA \}
\end{aligned}
\end{equation}

\begin{theorem}\label{tplocalisingsubcats}
	The \emph{tensor product of localising subcategories generated by a set} $\www_{\AAA} \boxtimes \www_{\BBB}$ is given by
	$$\www_1 \vee \www_2 = \langle \www_1 \cup \www_2 \rangle$$
	in the poset $W_{\mathsf{dg}}$ of localising subcategories of $\D(\AAA \dertp \BBB)$ generated by a set of objects.
\end{theorem}	
	In order to prove this result, we first provide an explicit description of the quasi-representable bimodule $\otimes$ between $\DD(\AAA) \dertp \DD(\BBB)$ and $\DD(\AAA) \boxtimes \DD(\BBB) \cong \DD(\AAA  \dertp \BBB)$ (see \S \ref{quotandtpwg}).
	\begin{lemma}\label{lemcantpderdg}
		Let $\AAA$ and $\BBB$ be small dg categories and consider the canonical bimodule $\otimes \in \RHom_{\ccm,\ccm}(\DD(\AAA) \dertp \DD(\BBB), \DD(\AAA \dertp \BBB))$. Then, given $F \in \DD(\AAA), G \in \DD(\BBB)$, we have that:
		\begin{equation}
		(F \hztp G) (A,B) = F(A) \dertp G(B)
		\end{equation}
		in $\D(k)$ for all  and $A \in \AAA, B \in \BBB$.
		\begin{proof}
			Recall that given $\CCC$ a small dg category, representables $\{\CCC(-,C)\}_{C \in \CCC}$ form a set of compact generators of $\D(\CCC)$. Consequently, we have that $F$ (resp. $G$) can be written in terms of representables in $\D(\AAA)$ (resp. in $\D(\BBB)$) by using coproducts, cones and shifts. Because $\hztp$ is bicocontinuous and exact in each variable, and thus commutes with coproducts, cones and shifts in both variables, $F \hztp G$ can be also written in $\D(\AAA \dertp \BBB)$ in terms of elements of the form $\AAA(-,A) \hztp \BBB(-,B)$ using direct sums, cones and shifts.
			
			Now recall that $\otimes$ is just the image of the identity in $H^0(\RHom_{\ccm}(\DD(\AAA \dertp \BBB),\DD(\AAA \dertp \BBB)))$ via the chain of isomorphisms
			\begin{equation*}
			\begin{aligned}
			\RHom_{\ccm}(\DD(\AAA \dertp \BBB),\DD(\AAA \dertp \BBB))) &\cong \RHom(\AAA \dertp \BBB,\DD(\AAA \dertp \BBB))\\ 
			&\cong (\RHom(\AAA,\RHom(\BBB,\DD(\AAA \dertp \BBB)))\\
			&\cong \RHom_{\ccm}(\DD(\AAA), \RHom_{\ccm}(\DD(\BBB), \DD(\AAA \dertp \BBB))\\
			&\cong \RHom_{\ccm,\ccm}(\DD(\AAA) \dertp \DD(\BBB), \DD(\AAA \dertp \BBB))
			\end{aligned}
			\end{equation*}
			defined above (see \eqref{eqc}, \eqref{eq1} and \eqref{innerhom1}). On the other hand, observe that the identity in $\RHom_{\ccm}(\DD(\AAA \dertp \BBB),\DD(\AAA \dertp \BBB))$ gets mapped under the first quasi-equivalence \eqref{eqc} to the Yoneda embedding $\AAA \dertp \BBB \lra \DD(\AAA \dertp \BBB)$. Therefore, when restricted to the representables, one just has that
			$$\AAA(-,A) \hztp \BBB(-,B) = (\AAA \dertp \BBB) (-, (A,B)).$$ 
			Now observe that 
			\begin{equation}
			\begin{aligned}
			\left( ( \AAA \dertp \BBB) (-, (A,B))\right)  (A',B') &= (Q(\AAA) \otimes \BBB) ((A',B'),(A,B))\\ 
			&= Q(\AAA) (A',A) \otimes \BBB(B',B)\\
			&= \AAA(A',A) \dertp \BBB(B',B),
			\end{aligned}
			\end{equation}
			where $Q$ denotes the cofibrant replacement functor in $\dgcat$, which can be chosen such that $Q(\AAA) \lra \AAA$ is the identity on objects (\Cref{cofibrantreplacement}). In addition, also by \Cref{cofibrantreplacement}, we have that the induced $Q(\AAA)(A',A) \lra \AAA(A',A)$ is a cofibrant replacement for $\AAA(A',A)$ in $\C(k)$.

			Recall that coproducts, cones and shifts are point-wise in $\D(\AAA \dertp \BBB)$, and hence the evaluation of $F \hztp G$ at any point $(A',B')$ can be written in terms of elements of the form $\AAA(A',A) \dertp \BBB(B',B)$ using coproducts, cones and shifts. But as $\dertp$ is bicocontinuous in $\D(k)$ and applying again that coproducts, cones and shifts are point-wise in $\D(\AAA)$ and $\D(\BBB)$, we obtain that 
			$$(F \hztp G) (A',B') = F(A') \dertp G(B')$$
			for all $(A',B') \in \AAA \dertp \BBB$ and we conclude. 
		\end{proof}
	\end{lemma}

We proceed now to prove \Cref{tplocalisingsubcats}:
\begin{proof}
	Let $\nnn_{\aaa}$ and $\nnn_{\bbb}$ be sets of generators of $\www_{\AAA}$ and $\www_{\BBB} $ respectively. Consider the set of compact objects $\GGG_{\DD(\AAA)} = \{\AAA(-,A)\}_{A \in \AAA}$ as a set of generators of $\D(\AAA)$ and respectively the set of compact objects $\GGG_{\DD(\BBB)} = \{ \BBB(-,B)\}_{B \in \BBB}$ as a set of generators of $\D(\BBB)$. By \Cref{generatedbyaset2variables} we know that $\www_1$ and $\www_2$ are localising subcategories in $\D(\AAA \dertp \BBB)$ generated by a set of objects. More concretely, it follows from \Cref{adjunctiontwovariables} and \Cref{generatedbyaset2variables} combined with \Cref{lemcantpderdg} that
	\begin{equation}
	\begin{aligned}
	\www_1 &= \langle \nnn_{\aaa} \hztp \GGG_{\DD(\BBB)} \rangle;\\
	\www_2 &= \langle \GGG_{\DD(\AAA)} \hztp \nnn_{\bbb} \rangle.
	\end{aligned}
	\end{equation} 
Hence, we can conclude that $$\www_{\AAA} \boxtimes \www_{\BBB}= \langle \www_{\AAA} \boxtimes_{\mathrm{Cl}} \www_{\BBB}\rangle = \langle \langle \nnn_{\aaa} \hztp \GGG_{\DD(\BBB)} \rangle \cup \langle \GGG_{\DD(\AAA)} \hztp \nnn_{\bbb} \rangle \rangle = \langle \www_1 \cup \www_2 \rangle = \www_1 \vee \www_2,$$ where the second equality is a direct consequence of \Cref{lemgen}.
\end{proof}

\subsection{Tensor product of dg Bousfield localisations}\label{tpstrictlocs}
Let $\AAA$, $\BBB$ be two small dg categories and consider the derived dg categories $\DD(\AAA)$ and $\DD(\BBB)$. Consider respective Bousfield localisations with kernels generated by a set of objects given by the dg subcategories $\LLL_{\AAA} \subseteq \DD(\AAA)$ and $\LLL_{\BBB} \subseteq \DD(\BBB)$ with respective quasi-left adjoints $F_{\AAA}$ and $F_{\BBB}$. Denote by $\www_{\AAA} = \Kern(H^0(F_{\AAA}))$ and $\www_{\BBB} = \Kern(H^0(F_{\BBB}))$ the corresponding localising subcategories generated by a set.

Consider the following full dg subcategories of $\DD(\AAA \dertp \BBB)$:
\begin{itemize}
	\item $\LLL_1 = \{F \in \DD(\AAA \dertp \BBB) \text{ | } F(-,B) \in \LLL_{\AAA} \text{ for all } B\in \BBB \} \subseteq \DD(\AAA \dertp \BBB)$;
	\item $\LLL_2 = \{F \in \DD(\AAA \dertp \BBB) \text{ | } F(A,-) \in \LLL_{\BBB} \text{ for all }A\in \AAA \} \subseteq \DD(\AAA \dertp \BBB)$.
\end{itemize}
The natural functors
\begin{itemize}
	\item $F_1: \D(\AAA \dertp \BBB) \lra H^0(\LLL_1): X \longmapsto \left( F_1(X):(A,B) \longmapsto H^0(F_{\AAA})(X(-,B))(A) \right) $;
	\item $F_2: \D(\AAA \dertp \BBB) \lra H^0(\LLL_2): X \longmapsto \left( F_2(X):(A,B) \longmapsto H^0(F_{\BBB})(X(A,-))(B) \right) $;
\end{itemize}	
can be easily seen to be the left adjoints for the inclusions $H^0(i_1):H^0(\LLL_1) \lra \D(\AAA \dertp \BBB)$ and $H^0(i_2): H^0(\LLL_2) \lra \D(\AAA \dertp \BBB)$ respectively. We have thus that $\LLL_1$and $\LLL_2$ are Bousfield localisations of $\DD(\AAA \dertp \BBB)$. Additionally, following the notations from \eqref{onesidedlocsubcats} above, one can observe that 
$$\Kern(F_1) = \{ F \in \D(\AAA \dertp \BBB) \text{ | } F(-,B) \in \www_{\AAA} \text{ for all }B\in \BBB\} = \www_1,$$
and analogously
$$\Kern(F_2) = \{ F \in \D(\AAA \dertp \BBB) \text{ | } F(A,-) \in \www_{\BBB} \text{ for all }A\in \AAA\} = \www_2.$$
As $\www_{\AAA}$ and $\www_{\BBB}$ are by hypothesis generated by a set, we have, as a consequence of \Cref{generatedbyaset2variables} above, that $\www_1 = \Kern(F_1)$ and $\www_2 = \Kern(F_2)$ are also generated by a set of objects. Hence $i_1$ and $i_2$ are Bousfield localisations of $\DD(\AAA \dertp \BBB)$ with kernel of the left adjoint at the $0^{\text{th}}$-cohomology level generated by a set of objects and we have the following:
\begin{proposition}\label{propcomparison}
	The localising subcategory $\www_1$ (resp. $\www_2$) and the well generated Bousfield localisation $\LLL_1$ (resp. $\LLL_2$) correspond under the isomorphism between $W_{\mathsf{dg}}$ and $L_{\mathsf{dg}}^{\op}$.
\end{proposition}

\begin{theorem}\label{deftpstrictloc}
	The tensor product $\LLL_{\AAA} \boxtimes \LLL_{\BBB}$ is given by
	$$\LLL_1 \wedge \LLL_2 = \LLL_1 \cap \LLL_2$$
	in the poset $L_{\mathsf{dg}}$ of dg Bousfield localisations of $\DD(\AAA \dertp \BBB)$ with kernel of the left adjoint at the $0^{\text{th}}$-cohomology level generated by a set of objects.
	\begin{proof}
		We have that: 
		\begin{equation*}
		\begin{aligned}
		\LLL_1 \cap \LLL_2 &= \LLL_1 \wedge \LLL_2 \\
		&= (\www_{1} \vee \www_{2})^{\perp}\\ 
		&= \langle \www_{1} \cup \www_{2} \rangle^{\perp}\\
		&= (\www_{\AAA} \boxtimes \www_{\BBB})^{\perp}\\
		&\cong  \DD(\AAA \dertp \BBB)/\www_{\AAA} \boxtimes \www_{\BBB}\\
		&\cong \DD(\AAA)/\www_{\AAA} \boxtimes \DD(\BBB)/\www_{\BBB}\\
		&= \LLL_{\AAA} \boxtimes \LLL_{\BBB}
		\end{aligned}
		\end{equation*}  
		where the first equality follows from \Cref{descriptionmeet} below, and the fourth is given by \Cref{tplocalisingsubcats}. 
	\end{proof}
\end{theorem}
\begin{lemma} \label{descriptionmeet}
	Let $\ccc$ be a well generated dg category. Given $\LLL$ and $\LLL'$ two dg Bousfield localisations of $\ccc$, we have that
	\begin{equation}
	\LLL \wedge \LLL' = \LLL \cap \LLL'
	\end{equation}
	in the poset $L_{\mathsf{dg}}$ of dg Bousfield localisations of $\ccc$ with kernel of the left adjoint at the $0^{\text{th}}$-cohomology level generated by a set of objects.
	\begin{proof}
		Observe we have that: 
		\begin{equation*}
		\begin{aligned}
		\LLL \wedge \LLL' &= (\www_{\LLL} \vee \www_{\LLL'})^{\perp}\\ 
		&= \langle \www_{\LLL} \cup \www_{\LLL'} \rangle^{\perp}\\
		&= \www_{\LLL}^{\perp} \cap \www_{\LLL'}^{\perp}\\
		&=  \LLL \cap \LLL'
		\end{aligned}
		\end{equation*}  
		where the first and last equalities are given by the isomorphism of posets described in \S\ref{posetisom}, the second by the description of the poset of localising subcategories generated by a set and the third by \Cref{leftorth}.
	\end{proof}
\end{lemma}

\section{Tensor product in terms of \texorpdfstring{$\alpha$}{alpha}-cocontinuous derived categories}\label{tpalphacocontinuous}

In this section we provide the description of the tensor product of well generated dg categories when we realise them as $\alpha$-cocontinuous dg categories. We make use of the description of the tensor product of Bousfield localisations of dg derived categories provided in \S\ref{tpstrictlocs}.
\begin{proposition}\label{alphaalphaderived}
	Let $\AAA, \BBB$ be two homotopically $\alpha$-cocomplete small dg categories and consider their respective $\alpha$-cocontinuous dg derived categories $\DD_{\alpha}(\AAA), \DD_{\alpha}(\BBB)$. Then we have that
	$$\DD_{\alpha}(\AAA) \boxtimes \DD_{\alpha}(\BBB) = \DD_{\alpha, \alpha}(\AAA \dertp \BBB),$$
	where $\DD_{\alpha, \alpha}(\AAA \dertp \BBB)$ denotes the full dg subcategory of $\DD(\AAA) \boxtimes \DD(\BBB) = \DD(\AAA \dertp \BBB)$ formed by the bimodules $F$ such that $F(A,-) \in \DD_{\alpha}(\BBB)$ for all $A \in \AAA$ and $F(-,B) \in \DD_{\alpha}(\AAA)$ for all $B \in \BBB$.
	\begin{proof}
		This follows from \Cref{deftpstrictloc}. 
	\end{proof}
\end{proposition}

Consider $\AAA, \BBB$ two homotopically $\alpha$-cocomplete small dg categories. We know that $\DD_{\alpha}(\AAA) \boxtimes \DD_{\alpha}(\BBB) = \DD_{\alpha, \alpha}(\AAA \dertp \BBB)$ is a well generated dg category, and hence, there exists a regular cardinal $\beta$ and a homotopically $\beta$-cocomplete small dg category $\CCC$ such that $\DD_{\alpha, \alpha}(\AAA \dertp \BBB) \cong \DD_{\beta}(\CCC)$. It is reasonable to ask the following questions:
\begin{itemize}
	\item Can we find such a $\CCC$ with $\beta = \alpha$? Or in other words, is the tensor product of $\alpha$-compactly generated dg categories again $\alpha$-compactly generated?
	\item Can $\CCC$ be found in terms of the provided $\AAA$ and $\BBB$?
\end{itemize}
The answer to both questions is affirmative (see \Cref{tpsmalltobig} and \Cref{alphapreserved} below). Showing this will be the main goal of this chapter.

\subsection{Tensor product of homotopically \texorpdfstring{$\alpha$}{alpha}-cocomplete dg categories} \label{secttpalphacocomplete}
Fixed a $\fraku$-small regular cardinal $\alpha$, we can define a homotopically $\alpha$-cocomplete tensor product in the full subcategory $\hodgcat_{\alpha}$ of $\hodgcat$ given by the homotopically $\alpha$-cocomplete $\fraku$-small dg categories. 
\begin{definition}
	Let $\AAA$ and $\BBB$ be homotopically $\alpha$-cocomplete dg categories. A \emph{homotopically $\alpha$-cocomplete tensor product} of $\AAA$ and $\BBB$ is defined as a homotopically $\alpha$-co\-com\-plete small dg category $\AAA \dertp_{\alpha} \BBB$ such that the following universal property holds in $\hodgcat_{\alpha}$:
	\begin{equation}
	\RHom_{\alpha}(\AAA \dertp_{\alpha} \BBB,\CCC) \cong \RHom_{\alpha}(\AAA,\RHom_{\alpha}(\BBB,\CCC)).
	\end{equation}
\end{definition}
\begin{remark}
	Observe that for $\alpha = \aleph_0$, as the homotopy category of a dg category is in particular $\Ab$-enriched, we have that for $\AAA,\BBB \in \hodgcat_{\aleph_0}$:
	\begin{itemize}
		\item $\RHom_{\aleph_0}(\AAA,\BBB) = \RHom(\AAA,\BBB)$;
		\item and hence $\AAA \dertp_{\aleph_0}\BBB = \AAA \dertp \BBB$.
	\end{itemize}
\end{remark}

\begin{remark}\label{topostheory}
	The following theorem, together with \Cref{tpsmalltobig}, constructs a homotopically $\alpha$-cocomplete dg category $\DDD$ such that $\DD_{\alpha}(\DDD) \cong \DD_{\alpha, \alpha}(\AAA \dertp \BBB)$ in $\hodgcat$, and shows that $\DDD$ is actually the homotopically $\alpha$-cocomplete tensor product of $\AAA$ and $\BBB$. The argument, despite the technicalities intrinsic to this setup, is essentially of topos theoretic nature. Let us describe here the outline of the proof roughly, ignoring the fact that we are working with cofibrant objects, and not just categories of dg modules, and that we are working with quasi-functors, instead of with dg functors. We first construct a candidate $\DDD$ for the homotopically $\alpha$-cocomplete tensor product of $\AAA$ and $\BBB$ together with a dg functor $F: \AAA \dertp \BBB \lra \DDD$ which is $\alpha$-cocontinuous in each variable. Intuitively, one can think of these small dg categories as ``dg Grothendieck sites''. Then, the fact that $F$ is $\alpha$-cocontinuous in each variable allows to observe that the restriction of scalars $F^*: \DD(\DDD) \lra \DD(\AAA \dertp \BBB)$ restricts to a map $F_s: \DD_{\alpha}(\DDD) \lra \DD_{\alpha, \alpha}(\AAA \dertp \BBB)$ between the ``categories of sheaves''. This is, in topos theoretical language, saying that $F$ is a ``continuous morphism of sites''. Then, using a parallel argument to that of classical topos theory, one has that $F_s$ has a left adjoint, that we will denote in the proof by $\Ind^{\alpha}_F$ such that 
	\begin{equation*}
	\begin{tikzcd}
	\AAA \dertp \BBB \arrow[d, "Y_{\AAA \dertp \BBB}"] \arrow[r, "F"] & \DDD \arrow[d, "Y_{\DDD}"] \\
	\DD(\AAA \dertp \BBB) \arrow[d, "a_{\AAA \dertp \BBB}"] \arrow[r, "F_!"] & \DD(\DDD) \arrow[d, "a_{\DDD}"] \\
	{\DD_{\alpha,\alpha}(\AAA \dertp \BBB)} \arrow[r, "\Ind^{\alpha}_F"] & \DD_{\alpha}(\DDD),
	\end{tikzcd}
	\end{equation*}
	is a commutative diagram. In particular, one has that $\Ind_{\alpha}F = a_{\DDD} \circ F_! \circ i_{\AAA \dertp \BBB}$, where $i_{\AAA \dertp \BBB}: \DD_{\alpha, \alpha}(\AAA \dertp \BBB) \subseteq \DD(\AAA \dertp \BBB)$ denotes the dg embedding. Then, by means of the concrete construction of $\DDD$, one can conclude, and we will do so combining \Cref{thmalphatp} and \Cref{tpsmalltobig}, that $\Ind_{\alpha}F$ is an isomorphism in $\hodgcat$.
\end{remark}

\begin{theorem}\label{thmalphatp} 

	Let $\alpha$ be a regular cardinal and $\AAA, \BBB$ homotopically $\alpha$-cocomplete $\fraku$-small dg categories. Then, there exists a homotopically $\alpha$-cocomplete $\fraku$-small dg category $\DDD$ such that 
	\begin{equation}
	\RHom_{\alpha}(\AAA, \RHom_{\alpha}(\BBB,\ccc)) \cong \RHom_{\alpha}(\DDD,\ccc)
	\end{equation}
	for all $\fraku$-well generated $\frakv$-small dg category $\ccc$. Moreover, we have that $\DDD = \AAA \dertp_{\alpha} \BBB$.
	\begin{proof}
		The construction of $\DDD$ will be obtained by mimicking the construction of the tensor product of $\alpha$-cocomplete $k$-linear categories following \cite[\S 6.5]{basic-concepts-enriched-category-theory}, or \cite[\S 10]{structures-defined-finite-limits-enriched-context} and \cite[\S 2.4]{tensor-product-finitely-cocomplete-abelian-categories} for the concrete case of $\alpha = \aleph_0$.
		
		Consider the Yoneda embedding $Y_{\AAA \dertp \BBB}:\AAA \dertp \BBB \lra \DD(\AAA \dertp \BBB)$ and the quasi-adjunction $a_{\AAA\dertp \BBB} \dashv_{H^0} i_{\AAA \dertp \BBB}$ where $i_{\AAA \dertp \BBB}: \DD_{\alpha,\alpha}(\AAA \dertp \BBB) \subseteq \DD(\AAA \dertp \BBB)$ is the natural inclusion. Recall that $a_{\AAA\dertp \BBB} \in \RHom_{\ccm}(\DD(\AAA \dertp \BBB),\DD_{\alpha,\alpha}(\AAA \dertp \BBB))$. 
		
		Consider the bimodule $a_{\AAA\dertp \BBB} \dertp_{\DD(\AAA \dertp \BBB)} Y_{\AAA\dertp \BBB} \in \RHom(\AAA \dertp \BBB,\DD_{\alpha,\alpha}(\AAA \dertp \BBB))$. We prove that $a_{\AAA\dertp \BBB} \dertp_{\DD(\AAA \dertp \BBB)} Y_{\AAA\dertp \BBB} $ is bi-$\alpha$-cocontinuous. Observe that that is the case if and only if 
		$$\D_{\alpha,\alpha}(\AAA \dertp \BBB)(H^0(a_{\AAA\dertp \BBB}) \circ H^0(Y_{\AAA\dertp \BBB})(-,-),X):( \AAA \dertp \BBB)^{\op} \lra \D(k)$$
		sends $\alpha$-small coproducts in both variables to $\alpha$-small products for all $X \in \D_{\alpha,\alpha}(\AAA \dertp \BBB)$, where we put $\D_{\alpha,\alpha}(\AAA \dertp \BBB) = H^0(\DD_{\alpha,\alpha}(\AAA \dertp \BBB))$. We have that 
		\begin{equation*}
		\begin{aligned}
		\D_{\alpha,\alpha}(\AAA \dertp \BBB)(H^0(a_{\AAA\dertp \BBB}) \circ H^0(Y_{\AAA\dertp \BBB})(-,-),X) &=\D(\AAA \dertp \BBB)(H^0(Y_{\AAA\dertp \BBB})(-,-),X) =\\
		&= H^0(X)(-,-),
		\end{aligned}
		\end{equation*}
		which, because $X \in \D_{\alpha,\alpha}(\AAA \dertp \BBB)$, sends $\alpha$-small coproducts in both variables to $\alpha$-small products. Consequently, $a_{\AAA\dertp \BBB} \dertp_{\DD(\AAA \dertp \BBB)} Y_{\AAA\dertp \BBB} $ is bi-$\alpha$-cocontinuous.   
				
		Denote by $\GGG$ the set of representables in $\DD(\AAA \dertp \BBB)$ and consider $\ttt \subseteq H^0(\DD_{\alpha, \alpha}(\AAA \dertp \BBB)) = \D_{\alpha,\alpha}(\AAA \dertp \BBB)$ the closure of $H^0(a_{\AAA \dertp \BBB})(\GGG)$ under $\alpha$-small coproducts. Denote by $\DDD \subseteq \qrep(\DD_{\alpha,\alpha}(\AAA \dertp \BBB))$ the enhancement of $\ttt$ via the natural enhancement $\qrep(\DD_{\alpha,\alpha}(\AAA \dertp \BBB))$ of $\D_{\alpha, \alpha}(\AAA \dertp \BBB)$. In particular, observe that $\DDD$ is an essentially small dg category which is homotopically $\alpha$-cocomplete.
		
		Consider the functor  
		$$F: \AAA \dertp \BBB \lra \DDD$$ 
		induced by the bimodule $a_{\AAA\dertp \BBB} \dertp_{\DD(\AAA \dertp \BBB)} Y_{\AAA\dertp \BBB} $, which remains bi-$\alpha$-cocontinuous.

		Consider $\ccc$ a well generated dg category. We are going to show that 
		\begin{equation}\label{alphabijectionlarge}
		\phi: [\DDD,\ccc]_{\alpha} \lra [\AAA \dertp \BBB, \ccc]_{\alpha,\alpha} : f \longmapsto f \circ [F]
		\end{equation}
		is a bijection. 

		Observe that we have
		\begin{equation}\label{chaineqvs}
		\begin{aligned}
		[\DD_{\alpha,\alpha}(\AAA \dertp \BBB),\ccc]_{\ccm} &\cong [\DD_{\alpha}(\AAA),\RHom_{\ccm}(\DD_{\alpha}(\BBB),\ccc)]_{\ccm} \\
		&\cong [\DD_\alpha(\AAA),\RHom_{\alpha}(\BBB,\ccc)]_{\ccm} \\
		&\cong  [\AAA,\RHom_{\alpha}(\BBB,\ccc)]_{\alpha} \\
		&\cong [\AAA \dertp \BBB, \ccc]_{\alpha,\alpha}
		\end{aligned}
		\end{equation}
		where the first bijection follows from the definition of the tensor product of well generated dg categories together with \Cref{alphaalphaderived}, the second from \Cref{derivedalphacocompletion}, the third from \Cref{notnecwelgenerated} together with \Cref{derivedalphacocompletion} and the last one from the $\dertp - \RHom$ adjunction.  Observe that an element $g \in [\DD_{\alpha,\alpha}(\AAA \dertp \BBB),\ccc]_{\ccm}$ gets sent to $g \circ  [Y'_{\AAA} \otimes_{H^0} Y'_{\BBB}]_{\iso} \in [\AAA \dertp \BBB, \ccc]_{\alpha,\alpha}$. If we denote by $a_{\AAA}$ (resp. $a_{\BBB}$) the quasi-left adjoint of the inclusion $i_{\AAA}: \DD_{\alpha}(\AAA) \subseteq \DD(\AAA)$ (resp. $\DD_{\alpha}(\BBB) \subseteq \DD(\BBB)$), it is easy to see, using the construction of the tensor product in terms of quotients as exposed in \Cref{thmboxquotwg}, that $g \circ  [Y'_{\AAA} \otimes_{H^0} Y'_{\BBB}]_{\iso} = g \circ [(a_{\AAA}\dertp_{\DD(\AAA)} Y_{\AAA}) \otimes_{H^0} (a_{\BBB} \dertp_{\DD(\BBB)} Y_{\BBB})]_{\iso} = g \circ [a_{\AAA\dertp \BBB}]_{\iso} \circ [Y_{\AAA \dertp \BBB}]$. We denote by $t_{\AAA \dertp \BBB, \ccc}: [\AAA \dertp \BBB, \ccc]_{\alpha,\alpha} \lra [\DD_{\alpha,\alpha}(\AAA \dertp \BBB),\ccc]_{\ccm}$ the inverse of this bijection.	
		
		We have a map
		\begin{equation}\label{inversealphasmall}
			[\AAA \dertp \BBB,\ccc]_{\alpha,\alpha} \overset{t_{\AAA\dertp \BBB, \ccc}}{\underset{\cong}{\lra}} [\DD_{\alpha, \alpha}(\AAA \dertp \BBB),\ccc]_{\ccm} \lra [\DDD,\ccc]_{\alpha} : f \longmapsto t_{\AAA \dertp \BBB, \ccc}(f) \circ j,
		\end{equation}
		where $j = [\bar{Y}_{\DD_{\alpha, \alpha}(\AAA \dertp \BBB)}]^{-1} \circ [i] \in [\DDD,\DD_{\alpha, \alpha}(\AAA \dertp \BBB)]$, with $i: \DDD \subseteq \qrep(\DD_{\alpha, \alpha}(\AAA \dertp \BBB))$ and $\bar{Y}_{\DD_{\alpha, \alpha}(\AAA \dertp \BBB)}: \DD_{\alpha, \alpha}(\AAA \dertp \BBB) \lra \qrep(\DD_{\alpha, \alpha}(\AAA \dertp \BBB))$ the natural quasi-equivalence provided by the Yoneda embedding.
		We are going to show that this is an inverse map of \eqref{alphabijectionlarge}.
		
		We have that $t_{\AAA \dertp \BBB}(f)\circ j \circ [F] = t_{\AAA \dertp \BBB}(f) \circ [a_{\AAA\dertp \BBB}]_{\iso} \circ [Y_{\AAA\dertp \BBB}] = f$ for any element $f \in [\AAA \dertp \BBB,\ccc]_{\alpha,\alpha}$. Hence \eqref{inversealphasmall} is a right inverse of \eqref{alphabijectionlarge}. 
		
		Now we want to show that $t_{\AAA \dertp \BBB, \ccc}(g \circ [F]) \circ j = g$. This equality is more involved and in order to prove it we will use the topos theoretical argument mentioned in \Cref{topostheory} above, which can also be seen as an $\alpha$-version of the usual extensions of dg functors. Denote by $\Ind_{F}^{\alpha} \coloneqq t_{\AAA \dertp \BBB, \DD_{\alpha}(\DDD)}([a_{\DDD}]_{\iso} \circ [Y_{\DDD}] \circ [F]) \in [\DD_{\alpha, \alpha}(\AAA \dertp \BBB), \DD_{\alpha}(\DDD)]_{\ccm}$. We hence have that
		\begin{equation} \label{firstequality}
		\Ind^{\alpha}_{F} \circ [a_{\AAA \dertp \BBB}]_{\iso} \circ [Y_{\AAA \dertp \BBB}] = [a_{\DDD}]_{\iso} \circ [Y_{\DDD}] \circ [F].
		\end{equation}
		Observe that $\Ind^{\alpha}_{F} \circ [a_{\AAA \dertp \BBB}]_{\iso} = [a_{\DDD}]_{\iso} \circ [F_!]$ and hence $\Ind^{\alpha}_{F} = [a_{\DDD}]_{\iso} \circ [F_!] \circ [i_{\AAA \dertp \BBB}]$.
		We claim that
		\begin{equation}\label{secondequality}
		\Ind^{\alpha}_{F} \circ j = [a_{\DDD}]_{\iso} \circ [Y_{\DDD}].
		\end{equation}
		Observe this will be enough to conclude. Indeed, as $\Ind^{\alpha}_{F}$ is cocontinuous, we have a diagram 
		\begin{equation}
		\begin{tikzcd}[column sep= large]
		{[\DDD,\ccc]_{\alpha}} \arrow[r, "\cong"',"s"] \arrow[d, "{(-)\circ [F]}"'] & {[\DD_{\alpha}(\DDD),\ccc]_{\ccm}} \arrow[d, "(-) \circ \Ind^{\alpha}_{F}"] \\
		{[\AAA \dertp \BBB,\ccc]_{\alpha,\alpha}} \arrow[r, "\cong","t_{\AAA \dertp \BBB, \ccc}"'] & {[\DD_{\alpha,\alpha}(\AAA \dertp \BBB),\ccc]_{\ccm}},
		\end{tikzcd}
		\end{equation}
		which is commutative as a direct consequence of \eqref{firstequality}, where $s$ denotes the inverse of the bijection $[\DD_{\alpha}(\DDD),\ccc]_{\ccm} \ra [\DDD,\ccc]_{\alpha}: f \mapsto f \circ [a_{\DDD}]_{\iso} \circ [Y_{\DDD}]$ from \Cref{derivedalphacocompletion}. Then, we have that 
		$$t_{\AAA \dertp \BBB, \ccc}(g \circ [F]) \circ j= s(f) \circ \Ind^{\alpha}_{F} \circ j = s(f) \circ [a_{\DDD}]_{\iso} \circ [Y_{\DDD}] = f,$$ 
		where the second equality comes from \eqref{secondequality}. Consequently, \eqref{inversealphasmall} is also a left inverse of \eqref{alphabijectionlarge}, which concludes the argument. 
		
		It hence only remains to prove that \eqref{secondequality} holds. Consider the dg functor $F: \AAA \dertp \BBB \lra \DDD$ and the associated restriction $F^*: \dgmod(\DDD) \lra \dgmod(\AAA \dertp \BBB)$ and extension $F_!: \DD(\AAA \dertp \BBB) \lra \DD(\DDD)$. Denote by $\DDD'$ the full dg subcategory of $\DD_{\alpha, \alpha}(\AAA \dertp \BBB)$ quasi-equivalent to $\DDD$ via the quasi-equivalence $\bar{Y}_{\AAA \dertp \BBB}$:
		\begin{equation*}
		\begin{tikzcd}
		{\DD_{\alpha,\alpha}(\AAA \dertp \BBB)} \arrow[r, "\bar{Y}_{\AAA \dertp \BBB}", "\sim"'] & {\qrep(\DD_{\alpha,\alpha}(\AAA \dertp \BBB))} \\
		\DDD' \arrow[u, "I", hook] \arrow[r, "G", "\sim"'] & \DDD \arrow[u, "i", hook]
		\end{tikzcd}
			\end{equation*}
		Observe that, for all $D \in \DDD'$, 
		\begin{equation*}
			\begin{aligned}
			F^* \circ Y_{\DDD} \circ G(D) &= F^*(h_{G(D)}) =\\
			&= \DDD(F(-),G(D)) =\\
			&= \qrep(\DD_{\alpha, \alpha}(\AAA \dertp \BBB))(i \circ F(-),i\circ G(D)) =\\
			&= \qrep(\DD_{\alpha, \alpha}(\AAA \dertp \BBB))(\Phi_{a_{\AAA \dertp \BBB}} \circ Y_{\AAA \dertp \BBB}(-),i\circ G(D)).
			\end{aligned}
		\end{equation*}		
		We hence have that
		\begin{equation*}
			\begin{aligned}
			i_{\AAA \dertp \BBB} \circ I(D)  &= \DD(\AAA \dertp \BBB)(Y_{\AAA \dertp \BBB}(-), i_{\AAA \dertp \BBB} \circ I(D)) \lra \\
			&\lra \qrep(\DD_{\alpha, \alpha}(\AAA \dertp \BBB))(\Phi_{a_{\AAA \dertp \BBB}} \circ Y_{\AAA \dertp \BBB} (-), \Phi_{a_{\AAA \dertp \BBB}} \circ i_{\AAA \dertp \BBB} \circ I(D)) =\\
			&= \qrep(\DD_{\alpha, \alpha}(\AAA \dertp \BBB))(\Phi_{a_{\AAA \dertp \BBB}} \circ Y_{\AAA \dertp \BBB} (-), i \circ G(D)).
			\end{aligned}
		\end{equation*}
		Consequently, we have a natural transformation $i_{\AAA \dertp \BBB} \circ I \lra F^* \circ Y_{\DDD} \circ G$. By adjunction, we have a natural transformation $F_! \circ i_{\AAA \dertp \BBB} \circ I \lra Y_{\DDD} \circ G$ and by composition a natural transformation
		\begin{equation}\label{naturaltransformationalpha}
			\alpha: \Phi_{a_{\DDD}} \circ F_! \circ i_{\AAA \dertp \BBB} \circ I \lra \Phi_{a_{\DDD}} \circ Y_{\DDD} \circ G.
		\end{equation} 
		Now, observe that every object $D \in H^0(\DDD')$ is isomorphic to $\coprod_{i\in I} H^0(G)^{-1} H^0(F)(A_i,B_i)$ where the coproduct is $\alpha$-small. Then we have that 
		\begin{equation}\label{homotopyisomorphism}
			\begin{aligned}
			& H^0(a _{\DDD}) \circ H^0(F_! \circ i_{\AAA \dertp \BBB} \circ I)(D) =\\
			&= H^0(a_{\DDD}) \circ H^0(F_! \circ i_{\AAA \dertp \BBB}) \left( \coprod_{i\in I} H^0(a_{\AAA \dertp \BBB}) \circ H^0(Y_{\AAA \dertp \BBB}) (A_i,B_i) \right) =\\
			&= H^0(a_{\DDD}) \circ H^0(F_! \circ i_{\AAA \dertp \BBB}) \circ H^0(a_{\AAA \dertp \BBB}) \left( \coprod_{i\in I} H^0(Y_{\AAA \dertp \BBB}) (A_i,B_i) \right) =\\
			&= H^0(a_{\DDD}) \circ H^0(F_!) \left( \coprod_{i\in I} H^0(Y_{\AAA \dertp \BBB})(A_i, B_i) \right)  =\\
			&= \coprod_{i\in I} H^0(a_{\DDD}) \circ H^0(F_!) \circ H^0(Y_{\AAA \dertp \BBB})(A_i, B_i) =\\
			&= \coprod_{i\in I} H^0(a_{\DDD}) \circ H^0(Y_{\DDD} \circ F) (A_i,B_i) =\\
			&= H^0(Y'_{\DDD}) \left( \coprod_{i\in I} H^0(F)(A_i,B_i)\right)  =\\
			&= H^0(a_{\DDD}) \circ H^0(Y_{\DDD} \circ G) (D),
			\end{aligned}
		\end{equation}
		where the only non-trivial equality is the third one. It follows from the fact that 
		$$H^0(a_{\DDD}) \circ H^0(F_! \circ i_{\AAA \dertp \BBB}) \circ H^0(a_{\AAA \dertp \BBB}) = H^0(a_{\DDD}) \circ H^0(F_!),$$
		which can be deduced by using the adjunctions $H^0(a_{\DDD}) \dashv H^0(i_{\DDD})$, $H^0(a_{\AAA \dertp \BBB}) \dashv H^0(i_{\AAA \dertp \BBB})$ and $H^0(F_!) \dashv H^0(F^*)$ together with the fact that the image of $H^0(F^*)(\D_{\alpha}(\DDD))$ lies in $\D_{\alpha,\alpha}(\AAA \dertp \BBB) \subseteq \D(\AAA \dertp \BBB)$. 
		From \eqref{homotopyisomorphism}, one can conclude that the natural transformation $\alpha$ from \eqref{naturaltransformationalpha} is a termwise homotopy equivalence and consequently, we have that $[a_{\DDD}]_{\iso} \circ [F_!] \circ [i_{\AAA \dertp \BBB}] \circ [I] = [a_{\DDD}]_{\iso} \circ [Y_{\DDD}] \circ [G]$.
		Thus, we have that
		\begin{equation*}
			\begin{aligned}
			\Ind^{\alpha}_{F} \circ j &= [a_{\DDD}]_{\iso} \circ [F_!] \circ [i_{\AAA \dertp \BBB}] \circ [\bar{Y}_{\DD_{\alpha, \alpha}(\AAA \dertp \BBB)}]^{-1} \circ [i] =\\
			&= [a_{\DDD}]_{\iso} \circ [F_!] \circ [i_{\AAA \dertp \BBB}] \circ [I] \circ [G]^{-1} = \\
			&= [a_{\DDD}]_{\iso} \circ [Y_{\DDD}] \circ [G] \circ [G]^{-1} =\\
			&= [a_{\DDD}]_{\iso} \circ [Y_{\DDD}],
			\end{aligned}
		\end{equation*}
		as we wanted to show.
		
		We hence have that $\phi: [\DDD,\ccc]_{\alpha} \lra [\AAA \dertp \BBB, \ccc]_{\alpha,\alpha} : f \longmapsto f \circ [F]$ is a bijection. Given another small dg category $\EEE$, we denote by $[\DDD \dertp \EEE,\ccc]'_{\alpha}$ the subset of $[\DDD \dertp \EEE,\ccc]$ of $\alpha$-cocontinuous morphisms in the first variable, and by $[(\AAA \dertp \BBB) \dertp \EEE,\ccc]'_{\alpha,\alpha}$ the subset of $[(\AAA \dertp \BBB) \dertp \EEE,\ccc]$ of $\alpha$-cocontinuous morphisms in both the first and second variables. We have the following diagram
		\begin{equation}
		\begin{tikzcd}[column sep= 46pt]
		\left[ \EEE, \RHom_{\alpha}(\DDD,\ccc)\right]  &\left[ \EEE, \RHom_{\alpha} (\AAA, \RHom_{\alpha}(\BBB,\ccc))\right] \\	
		\left[ \DDD \dertp \EEE, \ccc \right]'_{\alpha} \arrow[d,"\cong"'] \arrow[u,"\cong"]  &\left[(\AAA \dertp \BBB) \dertp \EEE, \ccc \right]'_{\alpha,\alpha} \arrow[d,"\cong"'] \arrow[u,"\cong"]\\
		\left[ \DDD, \RHom(\EEE,\ccc)\right]_{\alpha} \arrow[r,"- \circ {[F]}"]&\left[ \AAA \dertp \BBB, \RHom(\EEE,\ccc)\right]_{\alpha,\alpha}. 
		\end{tikzcd}
		\end{equation}
		Observe that $\RHom(\EEE,\ccc)$ is well generated as a direct consequence of \Cref{innerhomhalfway}, and hence the horizontal arrow is a bijection by \eqref{alphabijectionlarge}. Thus, as a direct consequence of Yoneda lemma, we have that $\RHom_{\alpha}(\DDD,\ccc) \cong \RHom_{\alpha} (\AAA, \RHom_{\alpha}(\BBB,\ccc))$ in $\hodgcat$ as we wanted to show.
		
		Now, given any homotopically $\alpha$-cocomplete small dg category $\CCC$, we want to show that $\phi': [\DDD,\CCC]_{\alpha} \lra [\AAA \dertp \BBB, \CCC]_{\alpha,\alpha} : f \longmapsto f \circ [F]$ is a bijection. From the argument above, we have that $\phi: [\DDD,\DD_{\alpha}(\CCC)]_{\alpha} \lra [\AAA \dertp \BBB, \DD_{\alpha}(\CCC)]_{\alpha,\alpha} : f \longmapsto f \circ [F]$ is a bijection. Observe that the corestriction $Y'_{\CCC}: \CCC \lra \DD_{\alpha}(\CCC)$ of the Yoneda embedding induces injections
		$$[\DDD,\CCC]_{\alpha} \subseteq [\DDD,\DD_{\alpha}(\CCC)]_{\alpha}$$
		and 
		$$[\AAA \dertp \BBB,\CCC]_{\alpha, \alpha} \subseteq [\AAA \dertp \BBB,\DD_{\alpha}(\CCC)]_{\alpha,\alpha}.$$
		It is then easy to check that $\phi'$ can be obtained as the restriction of $\phi$ to $[\DDD,\CCC]_{\alpha}$, and hence we have that $\phi'$ is injective. As the elements $H^0(F)(\AAA \dertp \BBB)$ generate $H^0(\DDD)$ under $\alpha$-small coproducts and $Y'_{\CCC}$ is $\alpha$-cocontinuous, we can conclude that it is also surjective. Then, a similar argument as above using the universal property of the internal hom and Yoneda lemma allows us to prove that $\RHom_{\alpha}(\AAA, \RHom_{\alpha}(\BBB,\CCC)) \cong \RHom_{\alpha}(\DDD,\CCC)$, showing that $\DDD = \AAA \dertp_{\alpha} \BBB$ as desired.
	\end{proof}
\end{theorem}

\subsection{Tensor product of \texorpdfstring{$\alpha$}{alpha}-cocontinuous derived dg categories} \label{smalltobig}
\begin{proposition}\label{tpsmalltobig}
	Let $\AAA$, $\BBB$ be two homotopically $\alpha$-cocomplete small dg categories. Then, we have that
	\begin{equation}
	\DD_{\alpha}(\AAA) \boxtimes \DD_{\alpha}(\BBB) \cong \DD_{\alpha}(\AAA \dertp_{\alpha} \BBB)
	\end{equation}
	in $\hodgcat_{\wg}$.
	\begin{proof}
		We have that: 
		\begin{equation*}
			\begin{aligned}
				\RHom_{\ccm}(\DD_{\alpha}(\AAA),\RHom_{\ccm}(\DD_{\alpha}(\BBB), \ccc)) &\cong \RHom_{\alpha}(\AAA,\RHom_{\alpha}(\BBB,\ccc))\\
				&\cong \RHom_{\alpha}(\AAA \dertp_{\alpha} \BBB,\ccc)\\
				&\cong \RHom_{\ccm}(\DD_{\alpha}(\AAA \dertp_{\alpha}\BBB),\ccc)
			\end{aligned}
		\end{equation*}
		for every well generated dg category $\ccc$, where the first isomorphism comes from \Cref{derivedalphacocompletion} together with \Cref{closed}, the second isomorphism follows from \Cref{thmalphatp} and the last isomorphism from \Cref{derivedalphacocompletion}. This concludes the argument.
	\end{proof}	
\end{proposition}
\begin{corollary}\label{alphapreserved}
	The tensor product of two $\alpha$-compactly generated dg categories is again $\alpha$-compactly generated.
	\begin{proof}
		The theorem follows from the enhanced Gabriel-Popescu theorem (Theorem \ref{thmenhGP}) and \Cref{tpsmalltobig} above.
	\end{proof}
\end{corollary}

\def\cprime{$'$}
\providecommand{\bysame}{\leavevmode\hbox to3em{\hrulefill}\thinspace}
\bibliography{bibliography}
\bibliographystyle{amsplain}
\end{document}